\documentclass[12pt]{amsart}
\usepackage{amsmath}
\usepackage{amssymb}
\usepackage{amsthm}
\usepackage{psfig}
\newcommand{\gap}{\hspace{12pt}}
\newcommand{\ignore}[1]{\relax}
\newcommand{\IGNORE}[1]{\relax}

\newcommand{\Line}{\operatorname{Line}}
\newcommand{\C}{\mathbb C}
\newcommand{\R}{\mathbb R}
\newcommand{\Z}{\mathbb Z}
\newcommand{\N}{\mathbb N}
\newcommand{\NN}{\mathcal N}
\newcommand{\Q}{\mathbb Q}

\newcommand{\E}{\mathcal E}
\newcommand{\M}{\mathcal M}
\newcommand{\bm}{\overline{\mathcal M}}

\newcommand{\DD}{\mathcal D}
\newcommand{\XX}{\mathcal X}
\newcommand{\HH}{\mathcal H}
\newcommand{\LL}{\mathcal Log}
\newcommand{\LLL}{\mathcal L}

\newcommand{\trop}{\operatorname{trop}}
\newcommand{\Val}{\operatorname{Val}}
\newcommand{\val}{\operatorname{val}}
\newcommand{\Arg}{\operatorname{Arg}}

\newcommand{\ntrop}{N^{\operatorname{irr}}_{\operatorname{trop}}}
\newcommand{\ntropr}{N_{\operatorname{trop}}}
\newcommand{\ntropR}{N^{\operatorname{irr}}_{\operatorname{trop},{\mathbb R}}}
\newcommand{\ntropRr}{N_{\operatorname{trop},{\mathbb R}}}
\newcommand{\ntropRW}{N^{\operatorname{irr}}_{\operatorname{trop},{\mathbb R},W}}
\newcommand{\ntropRWr}{N_{\operatorname{trop},{\mathbb R},W}}
\newcommand{\ncl}{N^{\operatorname{irr}}}
\newcommand{\nclr}{N}

\newtheorem{thm}{Theorem}

\newtheorem{lem}{Lemma}[section]

\newtheorem{cor}[thm]{Corollary}
\newtheorem{thma}[lem]{Theorem}
\newtheorem{coro}[lem]{Corollary}
\newtheorem{prop}[lem]{Proposition}

\theoremstyle{definition}
\newtheorem{defn}[lem]{Definition}
\newtheorem{exa}[lem]{Example}
\theoremstyle{remark}
\newtheorem{rmk}[lem]{Remark}
\newtheorem{rem}[lem]{Remark}

\newtheorem{exam}[lem]{Example}

\newcommand{\Rtr}{{\mathbb R}_{\operatorname{trop}}}
\newcommand{\tor}{(\C^*)^{2}}
\newcommand{\rtor}{(\R^*)^{2}}
\newcommand{\torn}{(\C^*)^{n}}

\newcommand{\sign}{\operatorname{Sign}}
\newcommand{\subdiv}{\operatorname{Subdiv}}
\newcommand{\irr}{\operatorname{irr}}
\newcommand{\mult}{\operatorname{mult}}

\newcommand{\conj}{\operatorname{conj}}
\newcommand{\dd}{\partial}
\newcommand{\am}{\mathcal{A}}
\newcommand{\ce}{\mathcal{C}}
\newcommand{\z}{\mathcal{Z}}
\newcommand{\D}{\mathcal{D}}
\newcommand{\cp}{{\mathbb C}{\mathbb P}}
\newcommand{\rp}{{\mathbb R}{\mathbb P}}
\newcommand{\pp}{{\mathbb P}}
\newcommand{\ppp}{{\mathcal P}}
\newcommand{\qq}{{\mathcal Q}}
\newcommand{\qqq}{{\mathcal R}}

\newcommand{\Log}{\operatorname{Log}}

\newcommand{\Int}{\operatorname{Int}}
\newcommand{\ov}{\operatorname{ov}}
\renewcommand{\setminus}{\smallsetminus}

\newcommand{\Area}{\operatorname{Area}}
\newcommand{\area}{\operatorname{Area}}

\begin{document}

\title
{Enumerative tropical algebraic geometry in $\R^2$}
\author{Grigory Mikhalkin}
\address{Department of Mathematics\\
University of Toronto\\
100 St. George St.\\
Toronto, Ontario, M5S 3G3 Canada}
\address{St. Petersburg Branch of Steklov Mathematical Institute,
Fontanka 27, St. Petersburg, 191011 Russia}
\email{mikha@math.toronto.edu}
\thanks{The author would like to acknowledge
a partial support of the NSF and NSERC}
\begin{abstract}
The paper establishes a formula for enumeration of curves
of arbitrary genus in toric surfaces.
It turns out that such curves can be counted by means
of certain lattice paths in the Newton polygon.
The formula was announced earlier in \cite{M}.

The result is established with the help of the so-called
tropical algebraic geometry.
This geometry allows one to replace complex toric varieties with the real
space $\R^n$ and holomorphic curves with certain piecewise-linear
graphs there.
\end{abstract}

\maketitle

\section{Introduction}
Recall the basic enumerative problem in the plane.
Let $g\ge 0$ and $d\ge 1$ be two numbers
and let $\z=(z_1,\dots,z_{3d-1+g})$ be a collection of points $z_j\in\cp^2$
in general position.
A holomorphic curve $C\subset\cp^2$ is parameterized by a Riemann surface
$\tilde{C}$ under a holomorphic map $\phi:\tilde{C}\to\cp^2$ so that
$C=\phi(\tilde{C})$.
Here we choose the minimal parametrization, i.e. such that no component
of $\tilde{C}$ is mapped to a point by $\phi$.
The curve $C$ is irreducible if and only if
$\tilde{C}$ is connected.
The number $N^{\irr}_{g,d}$ of irreducible
curves of degree $d$ and genus $g$
passing through $\z$ is finite and does not depend on the choice of $z_j$ as
long as this choice is generic.

Similarly we can set up the problem of counting all (not necessarily irreducible)
curves. Define the genus of $C\subset\cp^2$ to be $\frac{1}{2}(2-\chi(\tilde{C}))$.
Note that the genus can take negative values for reducible curves.
The number $N^{\mult}_{g,d}$ of curves of degree $d$ and genus $g$
passing through $\z$ is again finite and does not depend on the choice of $z_j$ as
long as this choice is generic.
Figure \ref{chisla}
lists some (well-known) first few numbers $N^{\irr}_{g,d}$ and $N^{\mult}_{g,d}$.
\begin{figure}[h]
\centerline{
\begin{tabular}{|r|r|r|r|r|}
\hline
$g\backslash d$ & 1 & 2 & 3 & 4\\
\hline
0 & 1 & 1 & 12 & 620\\
\hline
1 & 0 & 0 & 1 & 225\\
\hline
2 & 0 & 0 & 0 & 27\\
\hline
3 & 0 & 0 & 0 & 1\\
\hline
\end{tabular}
\ \hspace{1in}
\begin{tabular}{|r|r|r|r|r|}
\hline
$g\backslash d$ & 1 & 2 & 3 & 4\\
\hline
-1 & 0 & 3 & 21 & 666\\
\hline
0 & 1 & 1 & 12 & 675\\
\hline
1 & 0 & 0 & 1 & 225\\
\hline
2 & 0 & 0 & 0 & 27\\
\hline
\end{tabular}
}
\caption{\label{chisla} Numbers \protect{$N^{\irr}_{g,d}$ and $N^{\mult}_{g,d}$.}}
\end{figure}

The numbers $N^{\irr}_{g,d}$ are known as the Gromov-Witten invariants
of $\cp^2$ (see \cite{KM}) while the numbers $N^{\mult}_{g,d}$ are sometimes
called the multicomponent Gromov-Witten invariant.
One series of numbers determines another by a simple combinatorial
relation (see e.g. \cite{CH}).
A recursive relation which allows one to compute the numbers $N^{\irr}_{0,d}$
(and thus the numbers $N^{\mult}_{0,d}$) was given by Kontsevich.
This relation came from
the associativity of the quantum cohomology (see \cite{KM}).
In the arbitrary genus case
Caporaso and Harris \cite{CH} gave an algorithm (bases on a degeneration
of $\cp^2$) which allows one to compute
the numbers $N^{\mult}_{g,d}$ (and thus the numbers $N^{\irr}_{g,d}$).

The main result of this paper gives a new way of computation
of these numbers as well as the $\R$-counterparts of these numbers
(that appear in real algebraic geometry).
The number $N^{\mult}_{g,d}$ turns out to be the number of certain lattice paths
of length $3d-1+g$ in the triangle $\Delta_d\subset\R^2$ with vertices
$(0,0)$, $(d,0)$ and $(0,d)$. The paths have to be counted with certain
non-negative multiplicities. Furthermore, this formula works not only
for $\cp^2$ but for other toric surfaces as well.
For other toric surfaces we just have to replace the triangle $\Delta_d$
by other convex lattice polygons.
The polygon should be chosen so that it determines the corresponding (polarized)
toric surface.

The formula comes as an application of the so-called {\em tropical geometry}
whose objects are certain piecewise-linear polyhedral complexes in $\R^n$.
These objects are the limits of the amoebas of holomorphic varieties
after a certain degeneration of the complex structure.
The idea to use these objects for enumeration of holomorphic curves
is due to Kontsevich.

In \cite{KS} Kontsevich and Soibelman proposed a program
linking Homological Mirror Symmetry and torus fibrations
from the Strominger-Yau-Zaslow conjecture \cite{SYZ}.
The relation is provided by passing to the so-called
{\em ``large complex limit"} which deforms a complex
structure on a manifold to its worst possible degeneration.
Similar deformations appeared in other areas
of mathematics under different names.
The {\em patchworking} in Real Algebraic Geometry was
discovered by Viro \cite{V}.
Maslov and his school studied the so-called
dequantization of the semiring of positive real numbers
(cf. \cite{Ma}).
The limiting semiring is isomorphic to
the $(\max,+)$-semiring $\Rtr$, the semiring
of real numbers equipped with taking the maximum
for addition and addition for multiplication.

The semiring $\Rtr$ is known to computer scientists
as one of {\em tropical} semirings, see e.g.
\cite{Pin}. In Mathematics this semiring appears
from non-Archimedean fields $K$ under a certain
pushing forward to $\R$ of the arithmetic operations in $K$.

In this paper we develop some basic algebraic geometry
over $\Rtr$ with a view towards counting curves.
In particular, we rigorously set up some enumerative
problems over $\Rtr$ and prove their equivalence to the relevant
problems of complex and real algebraic geometry. The reader can refer to Chapter 9
of Sturmfels' recent book
\cite{S} for some first steps in tropical algebraic geometry.
See also \cite{SS}, \cite{SST}, \cite{SpSt}
for some of more recent development.

We solve the corresponding tropical enumerative problem in $\R^2$.
As an application we get a formula counting the number of curves
of given degree and genus in terms of certain lattice paths of a
given length in the relevant Newton polygon. In particular this
gives an interpretation of the Gromov-Witten invariants in $\pp^2$
and $\pp^1\times\pp^1$ via lattice paths in a triangle and a
rectangle respectively. This formula was announced in \cite{M}.
For the proof we use the patchworking side of the story which is
possible to use since the ambient space is 2-dimensional and the curves
there are hypersurfaces. An alternative approach (applicable to
higher dimensions as well) is to use the symplectic field theory
of Eliashberg, Givental and Hofer \cite{EGH}.
Generalization of this formula to higher
dimensions is a work in progress.
In this paper we only define the enumerative multiplicity
for the 2-dimensional case.
There is a similar definition
(though no longer localized at the vertices) for multiplicities
of isolated curves in higher-dimensional tropical enumerative problems.
However, in higher dimensions there might be families of tropical
curves (of positive genus) for enumerative problems
with finite expected numbers of solutions (this phenomenon
already appears for curves in $\R^3$ passing through a finite
collection of points in general position) which seem to pose
a serious problem (that perhaps asks for development
of tropical virtual classes).


The main theorems are stated in section 7 and proved in section 8.
In section 2 we define tropical curves geometrically
(in a way similar to webs of Aharony, Hanany and Kol \cite{AH}, \cite{AHK}).
In section 3 we exhibit them as algebraic objects
over the tropical semifield.
In section 4 we define the tropical enumerative problems in $\R^2$,
in section 5 recall those in $\tor$.
Section 6 is auxiliary for section 8 and deals with certain
piecewise-holomorphic piecewise-Lagrangian objects in $\tor$
called {\em complex tropical curves}.
An outline of the approach taken in this paper
can also be found in \cite{It}. A somewhat
different approach can be found in \cite{Shust}.

The author is grateful to Y. Eliashberg, K. Hori, I. Itenberg,
M. Kapranov, M. Kontsevich,
A. Okounkov, B. Sturmfels,
R. Vakil, O. Viro and J.-Y. Welschinger  for useful discussions.
The author is also indebted
to the referees for thorough and helpful remarks, corrections
and suggestions as well as to B. Bertrand, F. Bihan and H. Markwig
who found many inconsistencies and misprints in earlier
versions of the paper.

\section{Tropical curves as graphs in $\R^n$}
In this section we geometrically define tropical curves
in $\R^n$ and set up the corresponding enumerative problem.
We postpone the algebraic treatment of the tropical curves
(which explains the term ``tropical" among other things)
until the next section.

\subsection{Definitions and the first examples}
Let $\bar\Gamma$ be a {\em weighted} finite graph. The weights
are natural numbers prescribed to the edges. Clearly, $\bar\Gamma$
is a compact topological space. We make it non compact by removing
the set of all 1-valent vertices $\mathcal V_1$,
$$\Gamma=\bar\Gamma\setminus {\mathcal V_1}.$$

\begin{rmk}
Removal of the 1-valent vertices is due to a choice we made in the
algebraic side of the treatment.
In this paper we chose the semifield $\Rtr=(\R,\max,+)$
as our ``ground semifield" for tropical variety, see the
next section. The operation $\max$ plays the r\^ole of addition
and thus we do not have an additive zero.
Non-compactness of $\Gamma$ is caused by this choice.
Should we have chosen $\Rtr\cup\{-\infty\}$ instead
for our ground semifield
we wouldn't need to remove the 1-valent vertices but then we would
have to consider tropical toric compactification of the ambient
space $\R^n$ as well. The approach of this paper is chosen for
the sake of simplicity. The other approach has its own advantages
and will be realized in a forthcoming paper.
\end{rmk}

\begin{defn}
\label{tropcur} A proper map $h:\Gamma\to\R^n$ is called {\em a
parameterized tropical curve} if it satisfies to the following two
conditions.
\begin{itemize}
\item For every edge $E\subset\Gamma$ the restriction $h|_E$ is either
an embedding or a constant map.
The image $h(E)$ is contained in a line $l\subset\R^2$
such that the slope of $l$ is rational.
\item For every vertex
$V\in\Gamma$ we have the following property. Let
$E_1,\dots,E_m\subset\Gamma$ be the edges adjacent to $V$, let
$w_1,\dots,w_m\in\N$ be their weights and let
$v_1,\dots,v_m\in\Z^n$ be the primitive integer vectors
at the point $h(V)$
in the direction of the edges $h(E_j)$ (we take $v_j=0$
if $h(E_j)$ is a point). We have
\begin{equation}
\label{balance}
\sum\limits_{j=1}^m w_jv_j=0.
\end{equation}
\end{itemize}
We say that two parameterized tropical curves $h:\Gamma\to\R^n$
and $h':\Gamma'\to\R^n$ are {\em equivalent} if there exists a
homeomorphism $\Phi:\Gamma\to\Gamma'$ which respects the weights
of the edges and such that $h=h'\circ\Phi$.
We do not distinguish equivalent parameterized tropical curves.
The image $$C=h(\Gamma)\subset\R^n$$ is called the unparameterized
tropical curve or just a {\em tropical 1-cycle} if no connected
component of $\Gamma$ gets contracted to a point.
The 1-cycle $C$
is a piecewise-linear graph in $\R^n$ with natural weights on its edges
induced from the weights on $\Gamma$.
If $E$ is an edge of $C$ then $h^{-1}(E)$ is a union of subintervals
of the edges of $\Gamma$. The weight of $E$ is the sum of the weights
of these edges.
\end{defn}

\begin{rmk}
In dimension 2 the notion of tropical curve coincides with
the notion of $(p,q)$-webs introduced by Aharony, Hanany and Kol
in \cite{AHK} (see also \cite{AH}).
\end{rmk}

\begin{rmk}
The map $h$ can be used to induce
a certain structure on $\Gamma$ from the affine space $\R^n$.
It is an instance of the so-called {\em $\Z$-affine structure}.
For a graph $\Gamma$ such a structure is equivalent to a metric
for every edge of $\Gamma$. Here is a way to obtain such a metric
for the edges that are not contracted to a point.

Let $E\subset\Gamma$ be a compact edge of weight $w$
that is not contracted to a point by $h$.
Such edge is mapped to a finite
straight interval with a rational slope in $\R^n$.
Let $l$ be the length of a primitive rational vector in the
direction of $h(E)$.
We set the length of $E$ to be $\frac{|h(E)|}{lw}$.

Note that $\Gamma$ also has non-compact edges (they result
from removing one-valent vertices from $\bar\Gamma$).
Such edges are mapped to unbounded straight intervals by $h$.

It is possible to consider abstract tropical curves as graphs
equipped with such $\Z$-affine structures. Then tropical maps
(e.g. to $\R^n$) will be maps that respect such structure.
Abstract tropical curves have genus (equal to $b_1(\Gamma)$)
and the number of punctures (equal to the number of ends of $\Gamma$)
and form moduli space in a manner similar to that of the classical
Riemann surfaces.
This point of view will be developed in a forthcoming paper.
\end{rmk}

\begin{exa}
\label{linexa} Consider the union of three simple rays
$$Y=\{(x,0)\ |\ x\le 0\}\cup\{(0,y)\ |\ y\le 0\}\cup
\{(x,x)\ |\ x\ge 0\}\subset\R^2.$$ This graph (considered as a tautological
embedding in $\R^2$) is a tropical curve since
$(-1,0)+(0,-1)+(1,1)=0$. A parallel translation of $Y$ in any
direction in $\R^2$ is clearly also a tropical curve. This gives
us a 2-dimensional family of curves in $\R^2$. Such curves are
called {\em tropical lines}.
\end{exa}

\begin{rmk}
The term {\em tropical line} is justified
in the next section dealing with the underlying algebra.
So far we would like to note the following properties of this family,
see Figure \ref{lines}.
\begin{itemize}
\item For any two points in $\R^2$ there is a tropical line
passing through them. \item Such a line is unique if the choice of
these two points is generic. \item Two generic tropical lines
intersect in a single point.
\end{itemize}
\begin{figure}[h]
\centerline{\psfig{figure=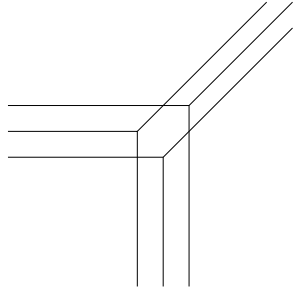,height=1.5in,width=1.5in}}
\caption{\label{lines} Three distinct tropical lines.}
\end{figure}
\end{rmk}

Somewhat more complicated tropical curves (corresponding to projective
curves of degree 3) are pictured on Figure \ref{cubics}

\subsection{The degree of a tropical curve in $\R^n$}
Let ${\mathcal T}=\{\tau_1,\dots,\tau_q\}\subset\Z^n$
be a set of non-zero integer vectors such that $\sum\limits_{j=1}^q\tau_j=0$.
Suppose that in this set we do not have positive multiples
of each other, i.e. if $\tau_j=m\tau_k$ for $m\in\N$
then $\tau_j=\tau_k$.
The degree of a tropical 1-cycle $C\subset\R^n$ takes
values in such sets according to the following construction.

By our definition a tropical curve $h:\Gamma\to\R^n$ has a finite number
of ends, i.e. unbounded edges (rays).
Let $\tau\in\Z^n$ be a primitive vector.
A positive multiple of $\tau$ is included in ${\mathcal T}$ if and only if
there exists an end of $\Gamma$ which is mapped in the direction of $\tau$.
In such case we include $m\tau$ into ${\mathcal T}$, where $m$
is the sum of multiplicities of all such rays.

\begin{defn}\label{multdeg}
The resulting set ${\mathcal T}$ is called {\em the toric
degree} of $C\subset\R^n$. Accordingly, the degree of
a parameterized tropical curve $h:\Gamma\to\R^n$ is the degree
of its image $h(\Gamma)$.
\end{defn}
Note that the sum of all vectors in ${\mathcal T}$ is zero.
This follows from adding the conditions \eqref{balance} from
Definition \ref{tropcur} in all
vertices of $C$.

For example the degree of both curves from Figure \ref{cuba}
is $\{(-1,-1)$, $(2,-1)$, $(-1,2)\}$ while the degree of both curves
from Figure \ref{cubics}
is
$\{(-3,0),(0,-3),(3,3)\}$.

\begin{defn}\label{prdeg}
If the toric degree of a tropical 1-cycle $C\subset\R^n$ is
$\{(-d,0,\dots,0),\dots,(0,\dots,0,-d),(d,\dots,d)\}$
then $C$ is called a tropical {\em projective curve of degree $d$}.
\end{defn}
The curves from Figure \ref{cubics} are examples of planar projective cubics.



\subsection{Genus of tropical curves and tropical 1-cycles}
We say that a tropical curve $h:\Gamma\to\R^n$ is {\em reducible} if
$\Gamma$ is disconnected. We say that a tropical 1-cycle $C\subset\R^n$
is reducible if it can be presented
as a union of two distinct tropical 1-cycles.
Clearly, every reducible 1-cycle can be presented as an image of
a reducible parameterized curve.
\begin{defn}
{\em The genus of a parameterized tropical curve} $\Gamma\to\R^n$
is $\dim H_1(\Gamma)-\dim H_0(\Gamma)+1$. In particular,
for irreducible parameterized curves the genus is the first Betti number of $\Gamma$.
{\em The genus of a tropical 1-cycle} $C\subset\R^n$ is the minimum
genus among all parameterizations of $C$.
\end{defn}
Note that according to this definition the genus can be negative.
E.g. the union of the three lines from
Figure \ref{lines} has genus $-2$.

If $C\subset\R^n$ is an embedded 3-valent graph
then the parameterization is unique.
However, in general, there might be several parameterizations of different
genus and taking the minimal value is essential.

\begin{exa}
The tropical 1-cycle on the right-hand side of Figure \ref{cubics} can be
parameterized by a tree once we ``resolve" its 4-valent vertex to make the
parameterization domain into a tree.
Therefore, its genus is 0.
\begin{figure}[h]
\centerline{\psfig{figure=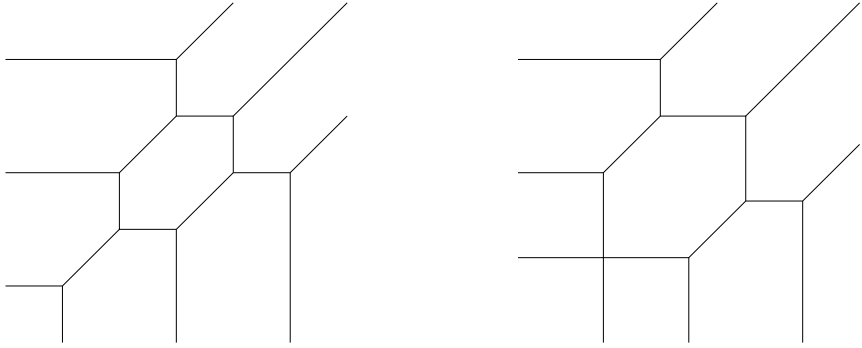,height=1.8in,width=4.5in}}
\caption{\label{cubics} A smooth projective tropical cubic and a
rational (genus 0) projective tropical cubic.}
\end{figure}
\end{exa}

\subsection{Deformations of tropical curves within
their combinatorial type}
As in the classical complex geometry case the deformation space
of a tropical curve $h:\Gamma\to\R^n$
is subject to the constraint coming from the Riemann-Roch formula.
Let $x$ be the number of ends of $\Gamma$.
\begin{rmk}
The number $-x$ is a tropical counterpart of the value
of the canonical class of the ambient complex variety on the curve $h(\Gamma)$.
The ambient space $\R^n$ corresponds to the torus $\torn$ classically.
Let $V\subset\torn$ be a holomorphic curve with a finite number of ends.
The space $\torn$ is not compact, but one can always choose
a toric compactification $\C T\supset\torn$ such that every point
of the closure $\bar{V}\supset V$
in $\C T$ intersects not more than one boundary divisor
(i.e. a component of $\C T\setminus\torn$).
Then every end of $V$ can be prescribed a multiplicity equal to
the intersection number of the point of $\bar{V}$ and
the corresponding boundary divisor. The value of the canonical
class of $\C T$ on $\bar{V}$ equals to the sum of these multiplicities.
\end{rmk}
\begin{defn}
The curves $h:\Gamma\to\R^n$  and $h':\Gamma\to\R^n$
(parameterized by the same graph $\Gamma$)
are said to be {\em of the same combinatorial type}
if for any edge $E\subset\Gamma$
the segments $h(E)$ and $h'(E)$ are parallel.
\end{defn}
Note that if two tropical curves $\Gamma\to\R^n$ are isotopic
in the class of tropical curves (with the same domain $\Gamma$)
then they are of the same combinatorial type.

The {\em valence} of a vertex
of $\Gamma$ is the number of adjacent edges regardless of their
weights. The graph $\Gamma$ is called 3-valent if every its
vertex is 3-valent. The parameterized tropical curve $h:\Gamma\to\R^n$
is called 3-valent if $\Gamma$ is 3-valent.

\begin{prop}
\label{RR}
Let $\Gamma$ be a 3-valent graph.
The space of all tropical curves $\Gamma\to\R^n$
in the same combinatorial type (up to their equivalence
from Definition \ref{tropcur}) is an open convex polyhedral domain
in a real affine $k$-dimensional space, where
$$k\ge x+(n-3)(1-g).$$
\end{prop}
\begin{proof}
It suffices to prove this for a connected graph $\Gamma$
since different components of $\Gamma$ vary independently.
and, furthermore, both sides of the inequality are
additive with respect to taking the union of components
(note that $1-g=b_0(\Gamma)-b_1(\Gamma)=\chi(\Gamma)$).
Let $\Gamma'\subset\Gamma$ be a finite tree containing all the
vertices of $\Gamma$. Note that the number of finite edges in
$\Gamma\setminus\Gamma'$ is $g$. By an Euler characteristic
computation we get that the number of finite edges of $\Gamma'$ is equal
to $x-3+2g$.

Maps $\Gamma'\to\R^n$ vary in a linear $(x-3+2g+n)$-dimensional
family if we do not change the slopes of the edges.
The $(x-3+2g)$-dimensional part comes from varying of the lengths of
the edges while the $n$-dimensional part comes from translations
in $\R^n$.
Such a map is extendable to a tropical map
$\Gamma\to\R^n$ if the pairs of vertices corresponding to the $g$
remaining edges define the lines with the correct slope. Each
of the $g$ edges imposes a linear condition of codimension at most
$n-1$. Thus tropical perturbations of $\Gamma\to\R^n$ are
contained in a linear family of dimension at least
$x-3+2g+n-(n-1)g=x+(n-3)(1-g)$. They form an open convex polyhedral
domain there defined by the conditions that the lengths of all
the edges are positive.
\end{proof}

\ignore{
\begin{coro}
A simple tropical curve $C\subset\R^n$ locally varies
in a (real) linear $k$-dimensional space, where
$$k\ge x+(n-3)(1-g).$$
\end{coro}
}

Consider the general case now and
suppose that $\Gamma$ has vertices of valence higher than 3.
How much $\Gamma$ differs from a 3-valent graph is measured
by the following characteristic.
Let the {\em overvalence} $\ov(\Gamma)$ be the sum of
the valences of all vertices of valence higher than 3 minus
the number of such vertices. Thus $\ov(\Gamma)=0$ if and only
if no vertex of $\Gamma$ has valence higher than 3.
\begin{prop}
\label{RRov}
The space of all tropical curves $\Gamma\to\R^n$
in the same combinatorial type (up to their equivalence
from Definition \ref{tropcur}) is an open convex polyhedral domain
in a real affine $k$-dimensional space, where
$$k\ge x+(n-3)(1-g)-\ov(\Gamma)-c,$$
where $c$ is the number of edges of $\Gamma$ that are mapped to a point.
\end{prop}
\begin{proof}
The proof is similar to that of Proposition \ref{RR}.
If the image of an edge is a point in $\R^n$ then we cannot
vary its length. Similarly we are lacking some degrees of freedom
(with respect to the set-up of Proposition \ref{RR}) if $\ov>0$.
\end{proof}
Note that $\ov+c$ can be interpreted as the overvalence
of the image $h(\Gamma)$.
\subsection{Changing the combinatorial type of $\Gamma$}
\begin{figure}[h]
\centerline{\psfig{figure=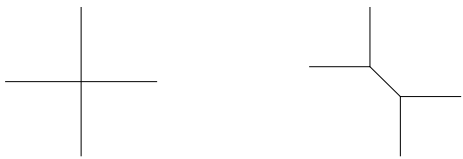,height=0.8in,width=2.4in}}
\caption{\label{3val} Smoothing a vertex of higher valence.}
\end{figure}

Sometimes we can deform $\Gamma$ {\em and} $h:\Gamma\to\R^n$
by the following procedure reducing $\ov$. If we have $n>3$  edges
adjacent to the same vertex then we can separate them into two
groups so that each group contains at least 2 edges.
Let us
insert a new edge $E'$ separating these groups as shown in Figure
\ref{3val}. This replaces the initial $n$-valent vertex with 2
vertices (the endpoints of $E'$) of smaller valence. There is a
``virtual slope" of $E'$ determined by the slopes of the edges in
each group. This is the slope to appear in local perturbation of the
tropical map $h:\Gamma\to\R^n$ (if such a perturbation exists).
Note that the weight of the new edge does not have to be equal to 1.

There is another modification of a tropical curve
near its vertex by changing the combinatorial type of $\Gamma$
which works even for some 3-valent vertices.
\begin{exa}
Let $\Gamma$ be the union of 3 rays in $\R^2$ in the direction
$(-2,1)$, $(1,-2)$ and $(1,1)$ emanating from the origin (pictured
on the left-hand side of Figure \ref{cuba}).
This curve is a simple tropical curve of genus 0.

It can be obtained as a $t\to 0$ limit of the family
of genus 1 curves given by the union of 3 rays in $\R^2$
in the direction $(-2,1)$, $(1,-2)$ and $(1,1)$ emanating from
and $(-2t,t)$, $(t,-2t)$ and $(t,t)$ respectively and the three
intervals $[(-2t,t),(t,-2t)]$, $[(-2t,t),(t,t)]$ and $[(t,t),(t,-2t)]$
as pictured in Figure \ref{cuba}.
\begin{figure}[h]
\centerline{\psfig{figure=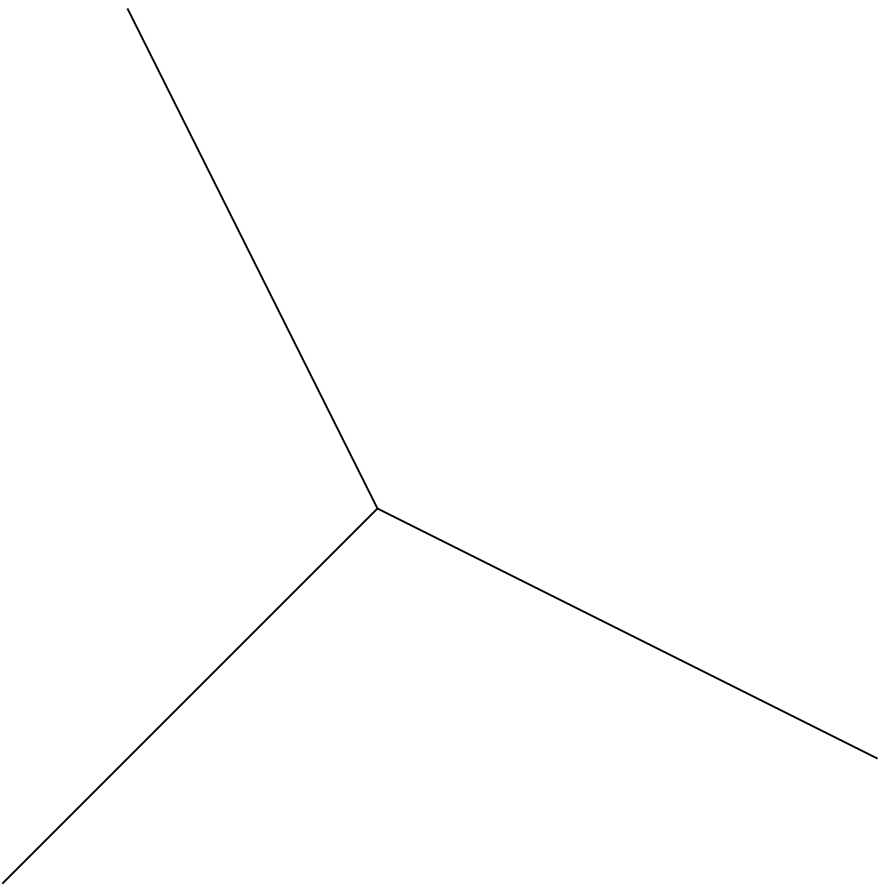,height=1.5in,width=1.5in}\hspace{1in}
\psfig{figure=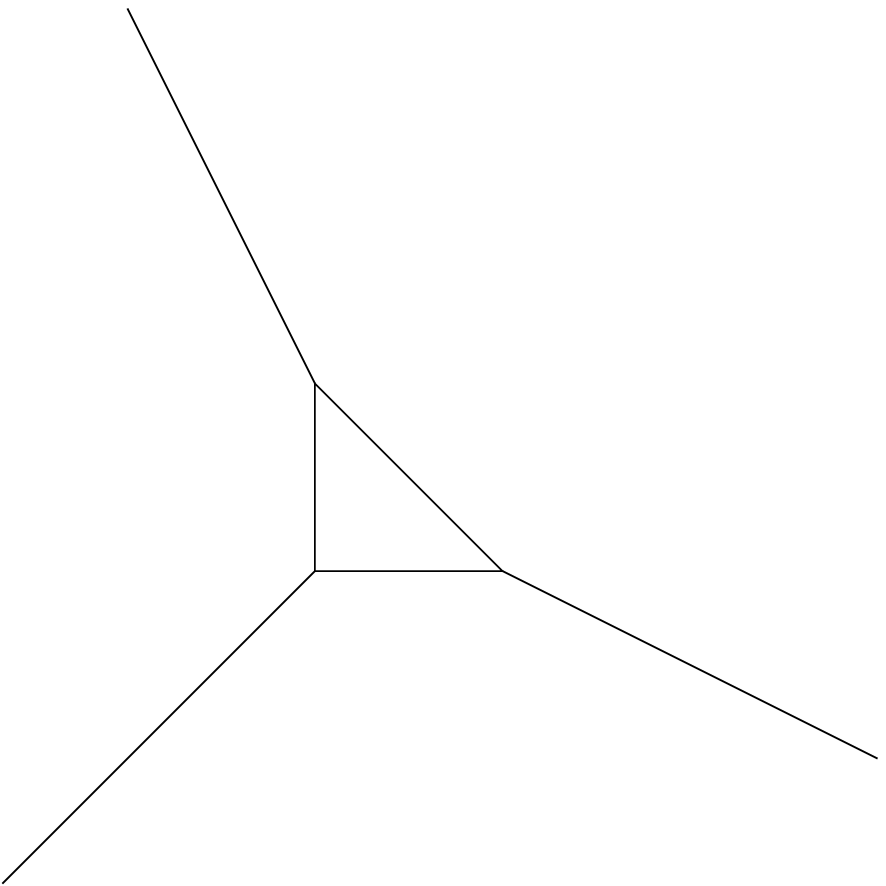,height=1.5in,width=1.5in}}
\caption{\label{cuba} Perturbation at a non-smooth 3-valent vertex}
\end{figure}
\end{exa}

Let $V$ be a 3-valent vertex of $\Gamma$.
As in Definition \ref{tropcur} let
$w_1,w_2,w_3$ be the weights of the edges adjacent to $V$ and let
$v_1,v_2,v_3$ be the primitive integer vectors in the direction of the edges.
\begin{defn}\label{multvert}
The {\em multiplicity} of $C$ at its 3-valent vertex $V$
is $w_1w_2|v_1\wedge v_2|$.
Here $|v_1\wedge v_2|$ is
the area of the parallelogram spanned
by $v_1$ and $v_2$.
Note that
$$w_1w_2|v_1\wedge v_2|=w_2w_3|v_2\wedge v_3|=w_3w_1|v_3\wedge v_1|$$
since $v_1w_1+v_2w_2+v_3w_3=0$ by Definition \ref{tropcur}.
\end{defn}
Note that the multiplicity of a vertex is always divisible
by the product of
the weights of any two out of the three adjacent edges.

\begin{defn}
We say that $h':\Gamma'\to\R^n$ is a perturbation of
$h:\Gamma\to\R^n$
if there exists a family $h'_t:\Gamma'\to\R^n$, $t>0$, in the same combinatorial
type as $h'$ and the pointwise limit $h'_0=\lim\limits_{t\to 0}h'_t$
such that $h'_0(\Gamma')$ coincides with $h(\Gamma)$ (as tropical 1-cycles).
\end{defn}

\begin{defn}
A tropical curve $h:\Gamma\to\R^n$ is called {\em smooth}
if $\Gamma$ is 3-valent, $h$ is an embedding and
the multiplicity of every vertex of $C$ is 1.
\end{defn}
\begin{prop}
A smooth curve does not admit perturbations of different
combinatorial types.
\end{prop}
\begin{proof}
Suppose that $h'_t:\Gamma'\to\R^n$ is a perturbation
of a smooth curve $h:\Gamma\to\R^n$.
Since $h$ is an embedding and $h(\Gamma)=h'_0(\Gamma')$
we have a map $$\psi:\Gamma'\to\Gamma.$$
Note that the weight of every edge from $\Gamma$ is 1
since otherwise the endpoints of multiple edges would
have multiplicity greater than 1. Thus the inverse image
of every open edge of $\Gamma$ under $\psi$ is a single
edge of $\Gamma'$.

Thus $\psi$ must be a homeomorphism near
the inner points of the edges of $\Gamma'$. Let $a\in\Gamma$
be a vertex and $U\ni a$ be its small neighborhood in $\Gamma$.
Note that $\psi^{-1}(U)$ is connected since $a$ is 3-valent
(otherwise we can divide the adjacent edges to $a$ into two
groups with zero sums of the primitive integer vectors).

Suppose that $\psi^{-1}(a)$ is not a point. Then $\psi^{-1}(a)$
is a graph which has 3 distinguished vertices that are adjacent
to the edges of $\psi^{-1}(U)\setminus\psi^{-1}(a)$.
The graph $A=\psi'_t(\psi^{-1}(U))$ must be contained in
the affine 2-plane in $\R^n$ containing the ends $A$.
This follows from the balancing condition for $A$.

The 3-valent vertices of $A$ have multiplicities from
Definition \ref{multvert}. Since $A$ is planar we can extend
the definition of the multiplicity to higher-valent vertices
as follows. Let $V\in A$ be a $k$-valent vertex,
$v_1,\dots,v_k$ be the primitive integer vectors in
the directions of the adjacent edges to $A$ numbered consistently
with the cyclic order in the ambient 2-plane and
$w_1,\dots,w_k$ be the corresponding weights.
We set the multiplicity of $V$ to be
$$\prod\limits_{l=2}^{k-1}|v_l\wedge \sum\limits_{j=1}^l v_j|.$$
It is easy to see that the multiplicity of $a$ in $\Gamma$
is equal to the sum of multiplicities of all the vertices of $A$.
The multiplicities of all vertices are positive integers.
Therefore, the multiplicity of $a$ is greater than $1$ unless
$\psi^{-1}(a)$ is a point.
\end{proof}

Proposition \ref{RR} can be generalized in the following way to incorporate
possible perturbations.
\begin{prop}
\label{RRg} The space of deformations of a parameterized tropical curve
$h:\Gamma\to\R^n$ is locally a cone
$${\mathcal C}=\bigcup\limits_j {\mathcal C_j},$$
where ${\mathcal C_j}$ is a polyhedral convex cone corresponding to
tropical curves parameterized by perturbations $\Gamma_j\to\R^n$
of $h:\Gamma\to\R^n$ of a given combinatorial type.
The union is taken over all possible combinatorial
types of perturbations.
We have $$\dim{\mathcal C_j}\ge x+(n-3)(1-g)-\ov(\Gamma_j)-c,$$
where $c$ is the number of the edges of $\Gamma_j$ that are mapped
to a point.
\end{prop}
\begin{proof}
This proposition follows from Proposition \ref{RRov} applied
to all possible perturbation of $\Gamma$.
\end{proof}

\begin{rmk}
Not all conceivable perturbations of $h:\Gamma\to\R^n$
are realized as the following example shows.
Let $C_1\subset\R^2\times\{0\}\subset\R^3$ be a
tropical 1-cycle of genus 1.
Let $C_2=\{y\}\times\R\subset\R^3$ be a vertical line such that
$y$ is a point on $C_1$ such that $C_1\setminus\{y\}$ is contractible.
The curve $$C=C_1\cup C_2$$
has a 4-valent vertex that cannot be perturbed (as any such perturbation
would force $C_1$ out of the plane $\R^2\times\{0\}$).
Thus any 3-valent perturbation of the tautological embedding
$C_1\cup C_2\subset\R^3$ has to have an edge mapping to a point.

This phenomenon is related to the so-called {\em superabundancy}
phenomenon.
\end{rmk}

\subsection{Superabundancy and regularity}
Some curves vary in a family strictly larger than
``the prescribed dimension" $x+(n-3)(1-g)-\ov-c$.
\begin{defn}
A parameterized tropical curve $h:\Gamma\to\R^n$ is called {\em regular}
if the space of the curves of this combinatorial type
(which is a polyhedral domain in an affine space by Proposition \ref{RRov})
has dimensions $x+(n-3)(1-g)-\ov-c$. Otherwise it is called {\em superabundant}.
\end{defn}
In contrast to the classical case
tropical superabundancy can be easily seen geometrically.
By the proof of Proposition \ref{RR} the superabundancy appears
if the cycles of the graph $\Gamma$ do not provide transversal conditions
for the length of the bounded edges of the subtree $\Gamma'$.
This is the case if some of the cycles of $C\subset\R^n$ are contained in
smaller-dimensional affine-linear subspaces of $\R^n$,
e.g. if a non-trivial cycle of $\Gamma$ gets contracted
or it a spatial curve develops a planar cycle.
More generally, this is the case if several non-degenerate
``spatial" cycles combine to a degenerate ``flat" cycle.

Clearly, no irreducible tropical curve of genus 0 can be superabundant
as it has no cycles. Furthermore, tropical immersions of 3-valent
graphs to the plane $\R^2$ are never superabundant as the following proposition shows.

\begin{prop}
\label{3vimm}
Every tropical immersion $h:\Gamma\to\R^2$ is regular if $\Gamma$
is 3-valent.
If $\Gamma$ has vertices of valence higher than 3 them
$h:\Gamma\to\R^2$ varies in at most $(x+g-2)$-dimensional family.
\end{prop}
\begin{proof}
Recall the proof of Propositions \ref{RR} and \ref{RRov}. Once again
we may assume that $\Gamma$ is connected. Let $V\in\Gamma$ be any vertex.

We may choose an order on the vertices of $\Gamma$ so that it is consistent
with the distance from $V$, i.e. so that the order of a vertex $V'$ is greater
than the order of a vertex $V''$ whenever $V'$ is strictly further from $V$
than $V''$.
The balancing condition for $h(\Gamma)$ implies the following
maximum principle for $\Gamma$:
{\em any vertex of $\Gamma$ is either adjacent to an unbounded edge of $\Gamma$
or is connected with a bounded edge to a higher order vertex}.
Inductively one may choose a maximal tree $\Gamma'\subset\Gamma$
so that this maximum principle also holds for $\Gamma'$.
Note that the the set of vertices of $\Gamma'$ coincides with
the set of vertices of $\Gamma$. Note also that our choice of
order on this set gives the orientation on the edges of $\Gamma$:
every edge is directed from a smaller to a larger vertex.

The space of deformation of $h$ within the same
combinatorial type is open in a $k$-dimensional real affine space
that is cut by $g$ hyperplanes in $\R^{l+n}$ where $l$ is the number
of bounded edges of $\Gamma'$.
Each of these $g$ hyperplanes is non-trivial if $h$ is an immersion
and $\Gamma$ is 3-valent, since then there can be no parallel
edges adjacent to the same vertex.

We have regularity if these hyperplanes intersect transversely.
The hyperplanes are given by a $g\times (l+2)$-matrix with real values.
The rows of this matrix correspond to the edges of $\Gamma\setminus\Gamma'$
while the first $l$ columns correspond to the edges of $\Gamma'$ (the
remaining $2$ columns correspond to translations in $\R^2$).
To show that the rank of this matrix is $g$
in the 3-valent case we exhibit an upper-triangular
$g\times g$-minor with non-zero elements on its diagonal.

For each edge $E$ of $\Gamma\setminus\Gamma'$ we include the column corresponding
to the (bounded) edge of $\Gamma'$ directed {\em toward} the highest endpoint
of $E$. If $\Gamma$ is 3-valent different edges of $\Gamma-\Gamma'$
correspond under this construction to different edges of $\Gamma'$.
This produces the required $g\times g$-minor.

If $\Gamma$ is not 3-valent then the number of the bounded edges
of $\Gamma'$ is $x-3+2g-\ov$. This number is the same as the
number of vertices of $\Gamma'$ other than $V$.
We can do the construction of the minor as above but only using
one edge of $\Gamma\setminus\Gamma'$ at every vertex of $\Gamma$ other
than $V$.
In such a way we can get a non-degenerate
$(g-\ov)\times (g-\ov)$-minor and thus the dimension is at most
$x-3+2g-\ov-(g-\ov)+2=x-1+g$ (2 comes from translations in $\R^2$).
If there exists a vertex of valence higher than 3
then we may choose such a vertex for the root $V$ of the tree $\Gamma'$
and this gives a non-degenerate minor of size strictly larger
than $(g-\ov)$.
\end{proof}
\begin{coro}
\label{RR2}
An immersed 3-valent tropical curve $h:\Gamma\to\R^n$ locally varies
in a (real) linear $k$-dimensional space, where
$$k= x+(n-3)(1-g)$$
if either $n=2$ or $g=0$.
\end{coro}
\begin{rmk}
There exist superabundant tropical immersions $\Gamma\to\R^2$
if $\Gamma$ is not 3-valent. A nice example is given by the Pappus
theorem configuration that is a union of 9 lines, see Figure \ref{pappus}.
\begin{figure}[h]
\centerline{\psfig{figure=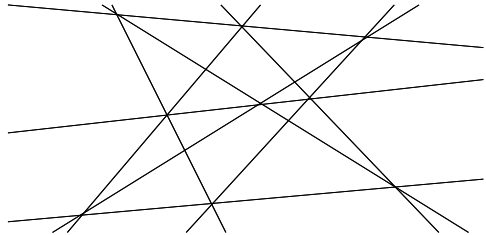,height=1.2in,width=2.5in}}
\caption{\label{pappus} Pappus configuration is superabundant.}
\end{figure}
Assume that the nine Pappus lines have rational slopes and
take $\Gamma$ to be their union in $\R^2$
so that our tropical curve $\Gamma\to\R^2$ is a tautological embedding.
We have $x=18$, $g=22$, $\ov=39$, $c=0$.
Therefore $x+g-1-\ov-c=0$, yet our configuration varies at least
in a 3-dimensional family (as we can
apply any translation and homothety in $\R^2$ without changing the slopes
of our lines.

Clearly, there also exist superabundant immersed 3-valent tropical
curves in $\R^n$, $n>2$. E.g. if $h:\Gamma\to\R^2$ is a (regular)
tropical immersion of a 3-valent graph $\Gamma$ then its composition
with the embedding $\R^2\subset\R^n$, $n>2$, is superabundant.
\end{rmk}

\section{Underlying tropical algebra}
In this section we exhibit the tropical curves as algebraic
varieties with respect to a certain algebra and also define
some higher-dimensional tropical algebraic varieties in $\R^n$.
\subsection{The tropical semifield $\Rtr$}
Consider the semiring $\Rtr$ of real numbers equipped with the
following arithmetic operations called {\em tropical} in Computer Science:
$$``x+y"=\max\{x,y\}\hspace{20pt}``xy"=x+y,$$
$x,y\in\Rtr$. We use the quotation marks to distinguish the
tropical operations from the classical ones. Note that addition is
idempotent, $``x+x=x"$. This makes $\Rtr$ to a semiring without
the additive zero (the role of such zero would be played by
$-\infty$).

\begin{rmk}
According to \cite{Pin} the term ``tropical" appeared in
Computer Science in honor of Brazil and, more specifically,
after Imre Simon (who is a Brazilian computer scientist)
by Dominique Perrin.
In Computer Science the term is usually applied to $(\min,+)$
semirings. Our semiring $\Rtr$ is $(\max,+)$ by our definition but
isomorphic to the $(\min,+)$-semiring, the isomorphism is
given by $x\mapsto -x$.
\end{rmk}

As usual, a (Laurent) polynomial in $n$ variables over $\Rtr$ is
defined by
$$f(x)=``\sum\limits_{j\in A} a_j x^j"=\max\limits_{j\in A} (<j,x>+a_j),$$
where $x=(x_1,\dots,x_n)\in\R^n$, $j=(j_1,\dots,j_n)$,
$x^j=x_1^{j_1}\dots x_n^{j_n}$,
${<j,x>}=j_1x_1+\dots+j_nx_n$ and $A\subset\Z^n$ is a finite set.
Note that $f:\R^n\to\R$ is a convex piecewise-linear function. It
coincides with the Legendre transform of a function $j\mapsto
-a_j$ defined on the finite set $A$.

\begin{defn}
The polyhedron $\Delta=\operatorname{Convex Hull}(A)$ is called
the {\em Newton polyhedron} of $f$. It can be treated
as a refined
version of the degree of the polynomial $f$ in toric geometry.
\end{defn}

\subsection{Tropical hypersurfaces: the variety of a tropical polynomial}
For a tropical polynomial $f$ in $n$ variables we define its {\em
variety} $V_f\subset\R^n$ as the set of points where the
piecewise-linear function $f$ is not smooth, cf. \cite{Ka},
\cite{M} and \cite{S}. In
other words, $V_f$ is the corner locus of $f$.
\begin{prop}
$V_f$ is the set of points in $\R^n$ where more than one
monomial of $f$ reaches its maximal value.
\end{prop}
\begin{proof}
If exactly one monomial of
$f(x)=``\sum\limits_{j\in A} a_j x^j"=\max\limits_{j\in A} (<j,x>+a_j)$
is maximal at $x\in\R^n$ then $f$ locally coincides with this monomial
and, therefore, linear and smooth. Otherwise $f$ has a corner at $x$.
\end{proof}

\begin{rmk}
At a first glance this definition might appear to be unrelated to
the classical definition of the variety as the zero locus of a
polynomial. To see a connection recall that there is no additive
zero in $\Rtr$, but its r\^ole is played by $-\infty$.

Consider the graph $\Gamma_f\subset\R^n\times\R$ of a tropical
polynomial $f:\R^n\to\R$. The graph $\Gamma_f$ itself is not a
tropical variety in $\R^{n+1}$ but it can be completed to the
tropical variety
$$\bar\Gamma_f=\Gamma_f\cup\{(x,y)\ |\ x\in V_f,\ y\le f(x) \},$$
see Figure \ref{parabola}.

\begin{prop}
$\bar\Gamma_f$ coincides with the variety of the polynomial in $(n+1)$
variables $``y+f(x)"$ (where $y\in\R$, $x\in\R^n$).
\end{prop}
\begin{proof}
If $(x,y)\in\Gamma_f$ then we have $y$ and one of the monomials of $f$
both reaching the maximal values in $``y+f(x)"=\max\{y,f(x)\}$.
If $x\in V_f$ and $y<f(x)$ then two monomials of $f(x)$ are reaching
the maximal value at the expression $``y+f(x)"$.
\end{proof}
Note that we have $V_f=\bar\Gamma_f \cap \{y=t\}$ for $t$
sufficiently close to $-\infty$. This is the sense in which $V_f$
can be thought as a zero locus.

One may argue that $\Gamma_f$ itself is a {\em subtropical}
variety (as in subanalytic vs. analytic sets) while $\bar\Gamma_f$
is its tropical closure. Figure \ref{parabola} sketches the graph
$y=``ax^2+bx+c"$ and its tropical closure.
\begin{figure}[h]
\centerline{\psfig{figure=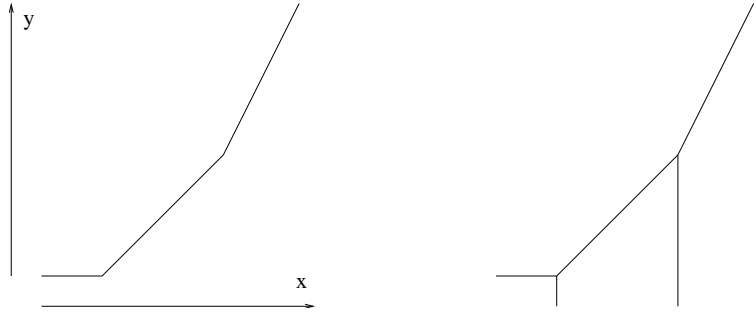,height=1.6in,width=3.9in}}
\caption{\label{parabola} The graph $y=``ax^2+bx+c"$ and its
closure, the tropical parabola}
\end{figure}
\end{rmk}

\begin{defn}
Varieties $V_f\subset\R^n$ are called {\em tropical hypersurfaces}
associated to $f$.
\end{defn}

\begin{rmk}\label{ambig}
Different tropical polynomials may define the same varieties.
To see this let us first extend the notion of concavity
to those $\R$-valued functions which are only defined on
a finite set $A\subset\R^n$.
We call a function $\phi:A\to\R$ {\em concave} if
for any (possibly non-distinct) $b_0,\dots,b_{n}\in A\subset\R^n$
and any $t_0,\dots,t_n\ge 0$ with $\sum\limits_{k=0}^n t_k=1$
and $\sum\limits_{k=0}^n t_kb_k\in A$ we have
$$\phi(\sum\limits_{k=0}^n t_kb_k)\ge
\sum\limits_{k=0}^n t_k\phi(b_k).$$ We have three types of
ambiguities when $f\neq g$ but $V_f=V_g$.
\begin{itemize}
\item $g=``x_j f"$, where $x_j$ is a coordinate in $\Rtr^n$. Note
that in this case the Newton polyhedron of $g$ is a translate of
the Newton polyhedron of $f$. \item $g=``c f"$, where $c\in\Rtr$
is a constant. \item The function $\Delta\cap\Z^n\ni j\mapsto a_j$
is not concave, where $f=``\sum\limits_{j\in A}a_jx^j"$ and we set
$a_j\mapsto -\infty$ if $j\notin A$. Then the variety of $f$
coincides with the variety of $g$ where $g$ is the smallest
concave function such that $g\ge f$ (in other words $g$ is a
concave hull of $f$).
\end{itemize}
Thus to define tropical hypersurfaces it suffices to consider only
tropical polynomials whose coefficients satisfy the concavity
condition above.
\end{rmk}

\begin{prop}[\cite{M-pp}]\label{Mdelta}
The space of all tropical hypersurfaces with a given Newton
polyhedron $\Delta$ is a closed convex polyhedral cone
$\M_\Delta\subset\R^m$, $m=\#(\Delta\cap\Z^n)-1$. The cone
$\M_\Delta\subset\R^m$ is well-defined up to the natural action of
$SL_m(\Z)$.
\end{prop}
\begin{proof}
All concave functions $\Delta\cap\Z^n\to\R, j\mapsto a_j$ form a
closed convex polyhedral cone $\tilde{\M}_\Delta\subset\R^{m+1}$.
But the function $j\mapsto a_j+c$ defines the same curve as the
function $j\mapsto a_j$. To get rid of this ambiguity we choose
$j'\in\Delta\cap\Z^n$ and define $\M_\Delta$ as the image of
$\tilde{\M}_\Delta$ under the linear projection $\R^{m+1}\to\R^m,
a_j\mapsto a_j-a_{j'}$.
\end{proof}

\subsection{Compactness of the space of tropical hypersurfaces}
Clearly, the cone $\M_\Delta$ is not compact. Nevertheless it gets
compactified by the cones $\M_{\Delta'}$ for all non-empty lattice
subpolyhedra $\Delta'\subset\Delta$ (including polygons with the
empty interior). Indeed, we have the following proposition.
\begin{prop}\label{compact}
Let $C_k\subset\R^n$, $k\in\N$ be a sequence of tropical
hypersurfaces whose Newton polyhedron is $\Delta$. There exists a
subsequence which converges to a tropical hypersurface $C$ whose Newton
polyhedron $\Delta_C$ is contained in $\Delta$ (note that $C$ is
empty if $\Delta_C$ is a point). The convergence is in the
Hausdorff metric when restricted to any compact subset in $\R^n$.
Furthermore, if the Newton polyhedron of $C$ coincides with
$\Delta$ then the convergence is in the Hausdorff metric in the
whole $\R^n$.
\end{prop}
\begin{proof}
Each $C_k$ is defined by a tropical polynomial
$f^{C_k}(x)=``\sum\limits_{j}a^{C_k}_jx^j"$. We may assume that
the coefficients $a^{C_k}_{j}$ are chosen so that they
satisfy the concavity condition and so that $\max\limits_{j}
a^{C_k}_{j}=0$. This takes care of the ambiguity in the choice of
$f^{C_k}$ (since the Newton polyhedron is already fixed).

Passing to a subsequence we may assume that $a^{C_k}_{j}$ converge
(to a finite number or $-\infty$) when $k\to\infty$ for all
$j\in\Delta\cap\Z^n$. By our assumption one of these limits is 0.
Define $C$ to be the variety of $``\sum
a^{\infty}_{j}x^j",$ where we take only finite coefficients
$a^{\infty}_{j}=\lim\limits_{k\to\infty} a^{C_k}_{j} > -\infty$.
\end{proof}

\subsection{Lattice subdivision of $\Delta$ associated to a tropical hypersurface}
\label{Dsubdiv} A tropical polynomial $f$ defines a lattice
subdivision of its Newton polyhedron $\Delta$ in the following way
(cf. \cite{GKZ}). Define the (unbounded) {\em extended polyhedral
domain}
$$\tilde{\Delta}=\operatorname{Convex Hull}
\{(j,t)\ |\ j\in A, t\le a_j  \}\subset\R^n\times\R.$$ The
projection $\R^{n}\times\R\to\R^n$ induces a homeomorphism from
the union of all closed bounded faces of $\tilde\Delta$ to
$\Delta$.
\begin{defn}
The resulting lattice subdivision $\operatorname{Subdiv}_f$ of $\Delta$
is called the {\em subdivision associated to $f$}.
\end{defn}

\begin{prop}\label{dualsubdiv}
The lattice subdivision $\operatorname{Subdiv}_f$ is dual to the
tropical hypersurface $V_f$.
Namely, for every $k$-dimensional polyhedron $\Delta'\in\operatorname{Subdiv}_f$
there is a convex closed (perhaps unbounded) polyhedron
$V_f^{\Delta'}\subset V_f\subset\R^n$.
This correspondence has the following properties.
\begin{itemize}
\item $V_f^{\Delta'}$ is contained
in an $(n-k)$-dimensional affine-linear subspace $L^{\Delta'}$ of $\R^n$ orthogonal to $\Delta'$.
\item The relative interior $U_f^{\Delta'}$
of $V_f^{\Delta'}$ in $L^{\Delta'}$ is not empty.
\item $V_f=\bigcup
U_f^{\Delta'}$.
\item $U_f^{\Delta'}\cap U_f^{\Delta''}=\emptyset$ if $\Delta'\neq\Delta''$.
\item $V_f^{\Delta'}$ is compact if and only if $\Delta'\subset\Delta$.
\end{itemize}
\end{prop}
\begin{proof}
For every $\Delta'\in\operatorname{Subdiv}_f$,
consider the
{\em truncated} polynomial $$f^{\Delta'}(x)=\sum\limits_{j\in\Delta'} a_jx^j$$
(recall that $f(x)=\sum\limits_{j\in\Delta} a_jx^j$).
Define
\begin{equation}\label{mnog}
V_f^{\Delta'}=V_f\cap \bigcap\limits_{\Delta''\subset\Delta'} V_{f^{\Delta''}}.
\end{equation}
Note that for any face $\Delta''\subset\Delta'$ we have the variety
$V_{f^{\Delta''}}$ orthogonal to $\Delta''$ (as moving in the direction
orthogonal to $\Delta''$ does not change the value of $\Delta''$-monomials)
and therefore to $\Delta'$.
To verify the last item of the proposition we restate the defining equation
\eqref{mnog} algebraically: $V_f^{\Delta'}$ is the set of points where
all monomials of $f$ indexed by $\Delta'$ have equal values while the value of
any other monomial of $f$ could only be smaller.
\end{proof}

\begin{exa}
Figure \ref{razb} shows the subdivisions dual to the curves from
Figure \ref{cubics} and \ref{cuba}.
\begin{figure}[h]
\centerline{\psfig{figure=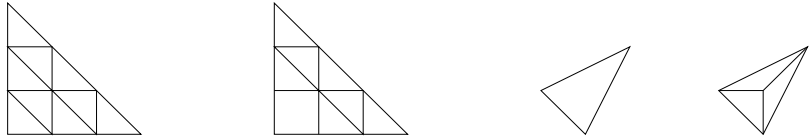,height=0.7in,width=4.2in}}
\caption{\label{razb} Lattice subdivisions associated to
the curves from Figure \ref{cubics}
and Figure \ref{cuba}.}
\end{figure}
\end{exa}

It was observed in \cite{Ka}, \cite{M} and \cite{S} that $V_f$ is
an $(n-1)$-dimensional polyhedral complex dual to the subdivision
$\operatorname{Subdiv}_f$. The complex $V_f\subset\R^n$ is a union
of convex (not necessarily bounded) polyhedra or {\em cells} of
$V_f$. Each $k$-cell (even if it is unbounded) of $V_f$ is dual to
a bounded $(n-k-1)$-face of $\tilde\Delta$, i.e. to an
$(n-k-1)$-cell of $\operatorname{Subdiv}_f$. In particular, the
slope of each cell of $V_f$ is rational.

In particular, an $(n-1)$ dimensional cell is dual to an interval
$I\subset\R^n$ whose both ends are lattice points. We define the
{\em lattice length} of $I$ as $\#(I\cap\Z^n)-1$. (Such length is
invariant with respect to $SL_n(\Z)$.) We can treat $V_f$ as a
weighted piecewise-linear polyhedral complex in $\R^n$, the
weights are natural numbers associated to the $(n-1)$-cells. They
are the lattice lengths of the dual intervals.

\begin{defn}\label{tip}
The {\em combinatorial type} of a tropical hypersurface
$V_f\subset\R^n$ is the equivalence class of all $V_g$ such that
$\operatorname{Subdiv}_g=\operatorname{Subdiv}_f$.
\end{defn}

Let $\mathcal S$ be such a combinatorial type.
\begin{lem}\label{linearitysubdiv}
All tropical hypersurfaces of the same combinatorial
type $\mathcal S$ form a convex polyhedral domain
$\M_{\mathcal S}\subset\M_{\Delta}$ that is open in
its affine-linear span.
\end{lem}
\begin{proof}
The condition $\operatorname{Subdiv}_f={\mathcal S}$
can be written in the following way in terms of the
coefficients of $f(x)=``\sum\limits_j a_jx^j"$.
For every $\Delta'\in{\mathcal S}$
the function $j\mapsto -a_j$ for $j\in\Delta'$
should coincide with some linear function $\alpha:\Z^n\to\R$
such that $-a_j>\alpha(j)$ for every $j\in\Delta\setminus\Delta'$.
\end{proof}

It turns out that the weighted piecewise-linear complex $V_f$
satisfies the balancing property at each $(n-2)$-cell,
see Definition 3 of \cite{M}.
Namely, let $F_1,\dots,F_k$ be the $(n-1)$-cells adjacent
to a $(n-2)$-cell $G$ of $V_f$. Each $F_j$ has a rational slope
and is assigned a weight $w_j$. Choose a direction of rotation
around $G$ and let $c_{F_j}:\Z^n\to\Z$ be linear maps whose
kernels are planes parallel to $F_j$ and such that they are primitive (non-divisible)
and agree with the chosen direction of rotation.
The balancing condition states that
\begin{equation}\label{ball}
\sum\limits_{j=1}^k w_jc_{F_j}=0.
\end{equation}
As it was shown in \cite{M} this balancing condition at every
$(n-2)$-cell of a rational piecewise-linear $(n-1)$-dimensional polyhedral
complex in $\R^n$ suffices for such
a polyhedral complex to be the variety of some tropical
polynomial.
\begin{thma}[\cite{M-pp}]
\label{trDelta} A weighted $(n-1)$-dimensional polyhedral complex
$\Pi\subset\R^n$ is the variety of a tropical polynomial if and
only if each $k$-cell of $\Pi$ is a convex polyhedron sitting in a
$k$-dimensional affine subspace of $\R^n$ with a rational slope
and $\Pi$ satisfies
the balancing condition \eqref{ball} at each $(n-2)$-cell.
\end{thma}
This theorem implies that the definitions of tropical curves
and tropical hypersurfaces agree if $n=2$.
\begin{coro}\label{curr2}
Any tropical curve $C\subset\R^2$ is a tropical hypersurface
for some polynomial $f$. Conversely, any tropical hypersurface
in $\R^2$ can be parameterized by a tropical curve.
\end{coro}
\begin{rmk}\label{degr2}
Furthermore, the degree of $C$ is determined by the Newton
polygon $\Delta$ of $f$ according to the following recipe.
For each side $\Delta'\subset\dd\Delta$ we take the primitive
integer outward normal vector and multiply it by the lattice
length of $\Delta'$ to get the degree of $C$.
\end{rmk}

\ignore{
\subsection{Smooth and singular tropical hypersurfaces}
\begin{defn}
A tropical hypersurface is called {\em smooth with a regular
$\infty$-behavior}
if the associated subdivision of its Newton polygon $\Delta$ (see
Definition \ref{subdiv}) is unimodular, i.e. each $n$-dimensional
polygon in the subdivision is a simplex of $n$-volume
$\frac{1}{n!}$.
A tropical hypersurface is called {\em smooth} if the associated
subdivision of $\Delta$ is a triangulation and every point of
$\Int\Delta\cap\Z^2$ is a vertex of this triangulation.
\end{defn}
}

\subsection{Tropical varieties and non-Archimedean amoebas}\label{Kap}
Polyhedral complexes resulting from
tropical varieties appeared in \cite{Ka} in the following context.
Let $K$ be a {\em complete algebraically closed non-Archimedean
field}. This means that $K$ is an algebraically closed field and
there is a {\em valuation} $\val:K^*\to\R$ defined on
$K^*=K\setminus\{0\}$ such that $e^{\val}$ defines a complete
metric on $K$. Recall that a valuation $\val$ is a map such that
$\val(xy)=\val(x)+\val(y)$ and
$\val(x+y)\le\max\{\val(x),\val(y)\}$.

Our principal example of such $K$ is a field of {\em Puiseux series}
with real powers. To construct $K$ we take the algebraic closure
$\overline{\C((t))}$ of the field of Laurent series $\C((t))$.
The elements of $\overline{\C((t))}$ are
formal power series in $t$ $a(t)=\sum\limits_{k\in A} a_kt^k$,
where $a_k\in\C$ and $A\subset\Q$ is a subset bounded from below and
contained in an arithmetic progression.
We set $\val(a(t))=-\min{A}$.
We define $K$ to be the completion of $\overline{\C((t))}$ as the
metric space with respect to the norm $e^{\val}$.

Let $V\subset (K^*)^n$ be an algebraic variety over $K$. The image
of $V$ under the map $\Val:(K^*)^n\to\R^n$,
$(z_1,\dots,z_n)\mapsto(\val(z_1),\dots,\val(z_n))$ is called {\em
the amoeba of $V$} (cf. \cite{GKZ}). Kapranov \cite{Ka} has shown
that the amoeba of a non-Archimedean hypersurface is the variety
of a tropical polynomial. Namely, if $\sum\limits_{j\in A}
a_jz^j=0$, $0\neq a_j\in K$ is a hypersurface in $(K^*)^n$ then
its amoeba is the variety of the tropical polynomial
$\sum\limits_{j\in A} \val(a_j)x^j$.

More generally, if $F$ is a field with a real-valued norm then the
amoeba of an algebraic variety $V\subset (F^*)$ is
$\Log(V)\subset\R^n$, where
$\Log(z_1,\dots,z_n)=(\log||z_1||,\dots,\log||z_n||)$.
Note that $\Val$ is such map with respect to the non-Archimedean
norm $e^{\val}$ in $K$.

Another particularly interesting case is
if $F=\C$ with the standard norm
$||z||=\sqrt{z\bar{z}}$ (see \cite{GKZ},\cite{Mi},\cite{PR}, etc.). The
non-Archimedean hypersurface amoebas appear as limits in the
Hausdorff metric of $\R^n$ from the complex hypersurfaces amoebas
(see e.g. \cite{M-pp}).

It was noted in \cite{SST} that the non-Archimedean approach can
be used to define tropical varieties of arbitrary
codimension in $\R^n$.
Namely, one can define
the {tropical varieties} in $\R^n$ to be the images $\Val(V)$
of arbitrary algebraic varieties $V\in(K^*)^n$.
This definition allows one to avoid dealing with the intersections
of tropical hypersurfaces in non-general position.
We refer to \cite{SST} for relevant discussions.

\section{Enumeration of tropical curves in $\R^2$}
\subsection{Simple curves and their lattice subdivisions}
Corollary \ref{curr2} states that any tropical 1-cycle in $\R^2$ is a
tropical hypersurface, i.e. it is the variety of a tropical
polynomial $f:\R^2\to\R$. By Remark \ref{ambig} the Newton
polygon $\Delta$ of such $f$ is well-defined up to a translation.
\begin{defn} We call $\Delta$ {\em the degree of a tropical
curve in $\R^2$.}
\end{defn}
By Remark \ref{degr2} this degree supplies the same amount
of information as the toric degree from Definition \ref{multdeg}.
We extract two numerical characteristics from the polygon $\Delta\subset\R^2$:
\begin{equation}\label{ml}
s=\#(\dd\Delta\cap\Z^2),\ l=\#(\Int\Delta\cap\Z^2).
\end{equation}
\ignore{
It is convenient also to consider $m=\#(\Delta\cap\Z^2)-1=s+l-1$.
The number $m$ is the dimension of the space of tropical curves
of degree $\Delta$ (since
$m+1$ is the number of coefficients in the corresponding polynomials).
}
The number $s$ is the number of unbounded edges of the curve
if counted with multiplicities (recall that we denoted the
number of unbounded edges ``counted simply" with $x\le s$).
The number $l$ is the genus of a smooth tropical curve of degree $\Delta$.
To see this let us note that every lattice point of $\Delta$
is a vertex of the associated subdivision for a smooth curve $C$.
Therefore, the homotopy type of $C$ coincides with $\Int\Delta\setminus\Z^2$.
Note also that smooth curves are dense in $\M_\Delta$.

There is a larger class of tropical curves in $\R^2$ whose behavior is as
simple as that of smooth curves.
\begin{defn}\label{simplecurve}
A parameterized tropical curve $h:\Gamma\to\R^2$ is called {\em simple}
if it satisfies to all of the following conditions.
\begin{itemize}
\item The graph $\Gamma$ is 3-valent.
\item The map $h$ is an immersion.
\item For any $y\in\R^n$ the inverse image $h^{-1}(y)$ consists
of at most two points.
\item If $a,b\in\Gamma$, $a\neq b$, are such that $h(a)=h(b)$ then neither $a$
nor $b$ can be a vertex of $\Gamma$.
\end{itemize}
A tropical 1-cycle $C\subset\R^2$ is called simple
if it admits a simple parameterization.
\end{defn}

\begin{prop}
\label{switch}
A simple tropical 1-cycle $C\subset\R^2$
admits a unique simple tropical
parameterization. The genus of a simple 1-cycle coincides
with the genus of its simple parameterization.
Furthermore, any of its non simple parameterization has
a strictly larger genus.
\end{prop}
\begin{proof}
By Definition \ref{simplecurve}
$C$ has only 3- and 4-valent vertices, where 4-valent vertices
are the double points of a simple immersion. Any other parametrization
would have to have a 4-valent vertex in the parameterizing graph.
\end{proof}
Proposition \ref{switch} allows us to switch back and forth
between parameterized tropical curves and tropical 1-cycles in
the case of simple curves in $\R^2$. Thus we refer to them just
as {\em simple tropical curves}. In a sense they are a tropical
counterpart of nodal planar curves in classical complex geometry.
\begin{rmk}
\label{immpar}
More generally, every tropical 1-cycle $C\subset\R^2$
admits a parameterization by an immersion of genus
not greater than $g(C)$. Start from
an arbitrary parameterization $h:\Gamma\to\R^2$.
To get rid from an edge $E\subset\Gamma$ such that
$h(E)$ is contracted to a point
we take the quotient of $\Gamma$ by $E$ for a
new domain of parameterization. This procedure does not
change the genus of $\Gamma$.

Therefore, we may assume that $h:\Gamma\to\R^2$ does
not have contracting edges. This is an immersion away
from such vertices of $\Gamma$ that there exist two
distinct adjacent edges $E_1,E_2$ that $h(E_1)\cap h(E_2)\neq\emptyset$.
Changing the graph $\Gamma$
by identifying the points on $E_1$ and $E_2$ with the same image
can only decreases the genus of $\Gamma$ (if $E_1$ and $E_2$ were
distinct edges connecting the same pair of vertices).
Inductively we get an immersion.
\end{rmk}

\begin{lem}
A tropical curve $C\subset\R^2$ is simple (see Definition \ref{simplecurve})
if and only if it is the variety of a tropical polynomial
such that $\operatorname{Subdiv}_f$ is a subdivision into triangles
and parallelograms.
\end{lem}
\begin{proof}
The lemma follows from Proposition \ref{dualsubdiv}. The 3-valent
vertices of $C$ are dual to the triangles of $\operatorname{Subdiv}_f$
while the intersection of edges are dual to the parallelograms
(see e.g. the right-hand side of Figure \ref{cubics} and the corresponding
lattice subdivision in Figure \ref{razb}).
\end{proof}

We have the following formula which expresses the genus
of a simple tropical curve
$V_f$ in terms of
the number $r$ of triangles in $\operatorname{Subdiv}_f$.
\begin{lem}[cf.  \cite{It}]\label{genusvertex}
If a curve $V_f\subset\R^2$ is simple then $g(V_f)=\frac{r-x}{2}+1$.
\end{lem}
\begin{proof}
Let $\Delta_0$ be the number of vertices of $\operatorname{Subdiv}_f$
while $\Delta_1$ and $\Delta_2$ be the number of its edges and
(2-dimensional) polygons.
Out of the $\Delta_2$ 2-dimensional polygons
$r$ are triangles and $(\Delta_2-r)$ are parallelograms.

We have
$$\chi(V_f)=-2\Delta_2+\Delta_1.$$
Note that $3r+4(\Delta_2-r)=2\Delta_1-x$.
Thus, $\Delta_1=\frac32 r + 2(\Delta_2-r)+\frac{x}{2}$ and
$$g(V_f)=1-\chi(V_f)=1+\frac{r-x}{2}.$$
\end{proof}

\subsection{Tropical general positions of points in $\R^2$}
\begin{defn}\label{genpos}
Points $p_1,\dots,p_k\in\R^2$ are said to be {\em in general
position tropically} if for any tropical curve
$h:\Gamma\to\R^2$ of
genus $g$ and with $x$ ends
such that $k\ge g+x-1$ and
$p_1,\dots,p_k\in h(\Gamma)$ we have the following conditions.
\begin{itemize}
\item The curve $h:\Gamma\to\R^2$ is simple (see Definition \ref{simplecurve}).
\item Inverse images $h^{-1}(p_1),\dots,h^{-1}(p_k)$ are
disjoint from the vertices of $C$.
\item $k=g+x-1$.
\end{itemize}
\end{defn}
\begin{exa} Two
distinct points $p_1,p_2\in\R^2$ are in general position tropically if
and only if the slope of the line in $\R^2$ passing through $p_1$
and $p_2$ is irrational.
\end{exa}

\begin{rem}
\label{redkriv} Note that we can always find a curve
with $g+x-1=k$ passing through $p_1,\dots,p_k$.
For such a curve we can
take a reducible curve consisting of $k$ affine (i.e. classical)
lines in $\R^2$ with rational
slope each passing through its own point $p_j$. This curve has
$2k$ ends while its genus is $1-k$.
\end{rem}

\begin{prop}\label{trsubset}
Any subset of a set of points in tropically general position is
itself in tropically general position.
\end{prop}
\begin{proof}
Suppose the points $p_1,\dots,p_j$ are not in general position.
Then there is a curve $C$ with $x$ ends of genus $j+2-x$ passing
through $p_1,\dots,p_j$ or of genus $j+1-x$ but with a non-generic
behavior with respect to $p_1,\dots,p_j$. By Remark \ref{redkriv}
there is a curve $C'$ passing through $p_{j+1},\dots,p_k$ of genus
$k-j+1-x'$. The curve $C\cup C'$ supplies a contradiction.
\end{proof}

\begin{prop}
\label{cl} For each $\Delta\subset\R^2$
the set of configurations $\ppp=\{p_1,\dots,p_k\}\subset\R^2$ such that
there exists a curve $C$ of degree $\Delta$ such that
the conditions of Definition \ref{genpos} are violated by $C$
is closed and nowhere dense.
\end{prop}
\begin{proof}
By Remark \ref{immpar} it suffices to consider
only immersed tropical curves $h:\Gamma\to\R^2$.
We have only finitely many combinatorial types of
tropical curves of genus $g$ with the Newton polygon
$\Delta$ as there are only finitely many lattice subdivisions
of $\Delta$. By Proposition \ref{3vimm}
for each such combinatorial type we have an
$(x+g-1)$-dimensional family of simple curves
or a smaller-dimensional family of non-simple curves.
For a fixed $C$ each of the $k$ points $p_j$ can vary
in a 1-dimensional family on $C$
or in a 0-dimensional family if $p_j$ is a vertex of $C$.
Thus the dimension of the
space of ``bad'' configurations $\ppp\in\operatorname{Sym}^k(\R^2)$
is at most $2k-1$.
\end{proof}

\begin{coro}\label{genposi}
The configurations $\ppp=\{p_1,\dots,p_k\}$ in general position tropically
form a dense set which can be obtained as an intersection of countably
many open dense sets in $\operatorname{Sym}^k(\R^2)$.
\end{coro}

\subsection{Tropical enumerative problem in $\R^2$}\label{trenpr}
To set up an enumerative problem we fix
the {\em degree}, i.e. a polygon $\Delta\subset\R^2$
with $s=\#(\dd\Delta\cap\Z^2)$, and
the {\em genus}, i.e. an integer number $g$.
Consider a configuration $\ppp=\{p_1,\dots,p_{s+g-1}\}\subset\R^2$
of $s+g-1$ points in tropical general position.
Our goal is to count tropical curves $h:\Gamma\to\R^2$ of genus $g$
such that $h(\Gamma)\supset\ppp$ and has degree $\Delta$.
\begin{prop}\label{tropfinite}
There exist only finitely many such curves $h:\Gamma\to\R^2$.
Furthermore, each end of $C=h(\Gamma)$ is of weight 1 in this case,
so $\Gamma$ has $s$ ends.
\end{prop}
Finiteness follows from Lemma \ref{1comb} proved
in the next subsection.
If $C$ has ends whose weight is greater than 1 then
the number of ends is smaller than $s$ and the existence
of $C$ contradicts to the general position of $\ppp$.
Recall that since $\ppp$ is in general position any
such $C$ is also simple and the vertices of $C$
are disjoint from $\ppp$.

\begin{figure}[h]
\centerline{\psfig{figure=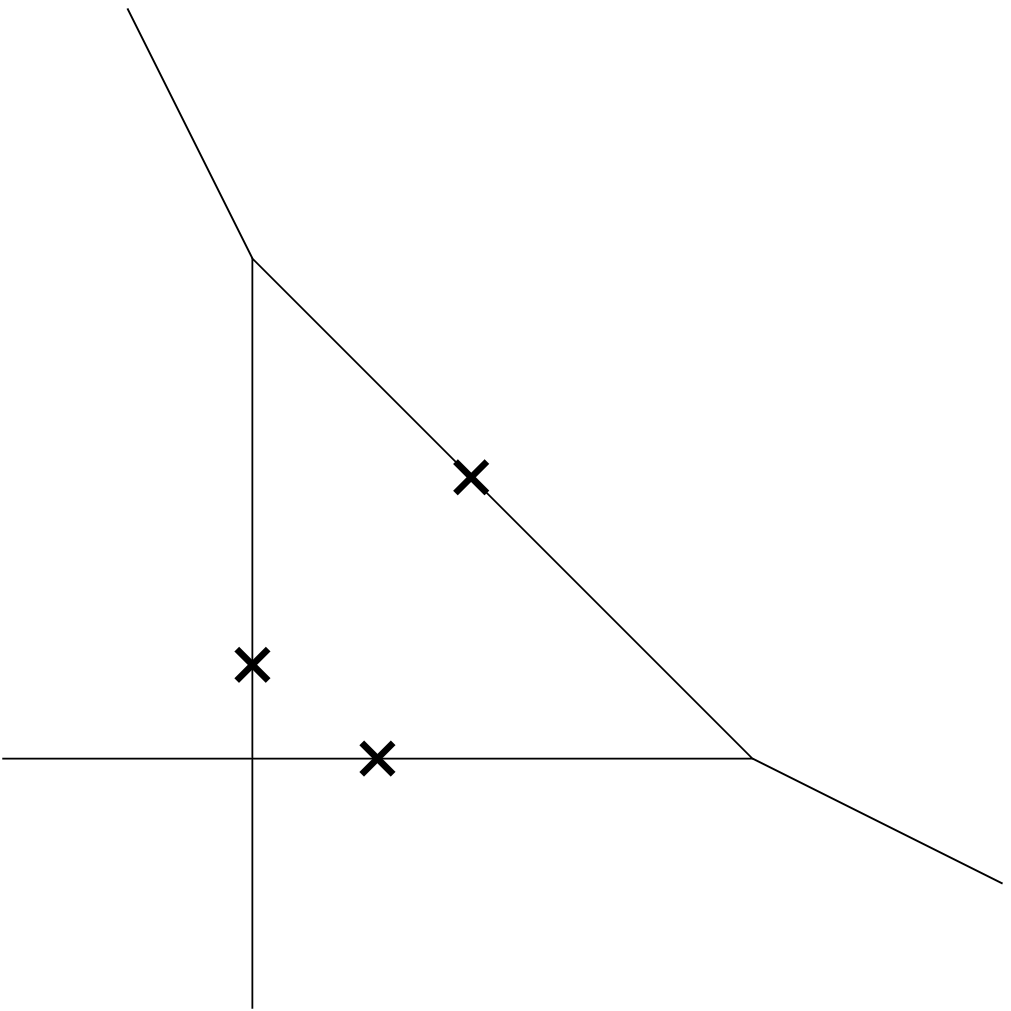,height=1.5in,width=1.5in}
\hspace{.1in}
\psfig{figure=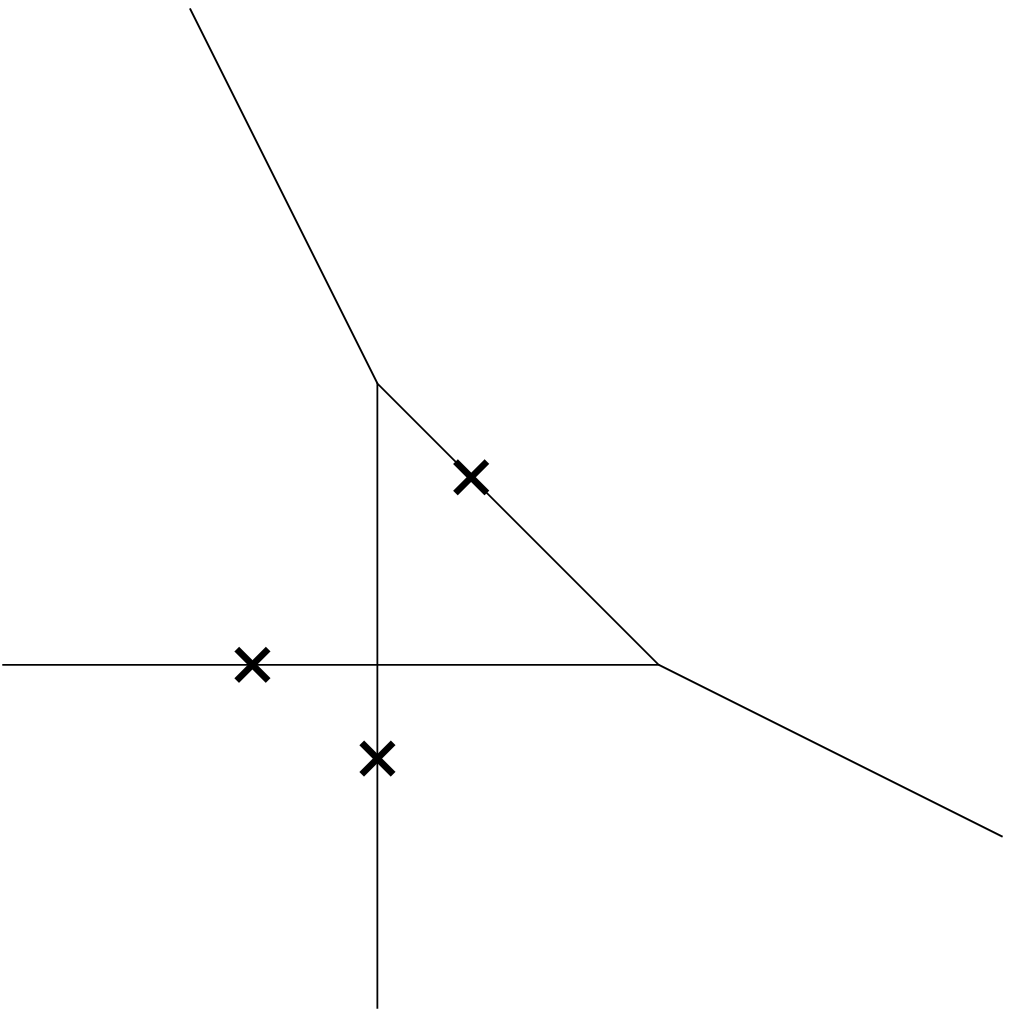,height=1.5in,width=1.5in}
\hspace{.1in}
\psfig{figure=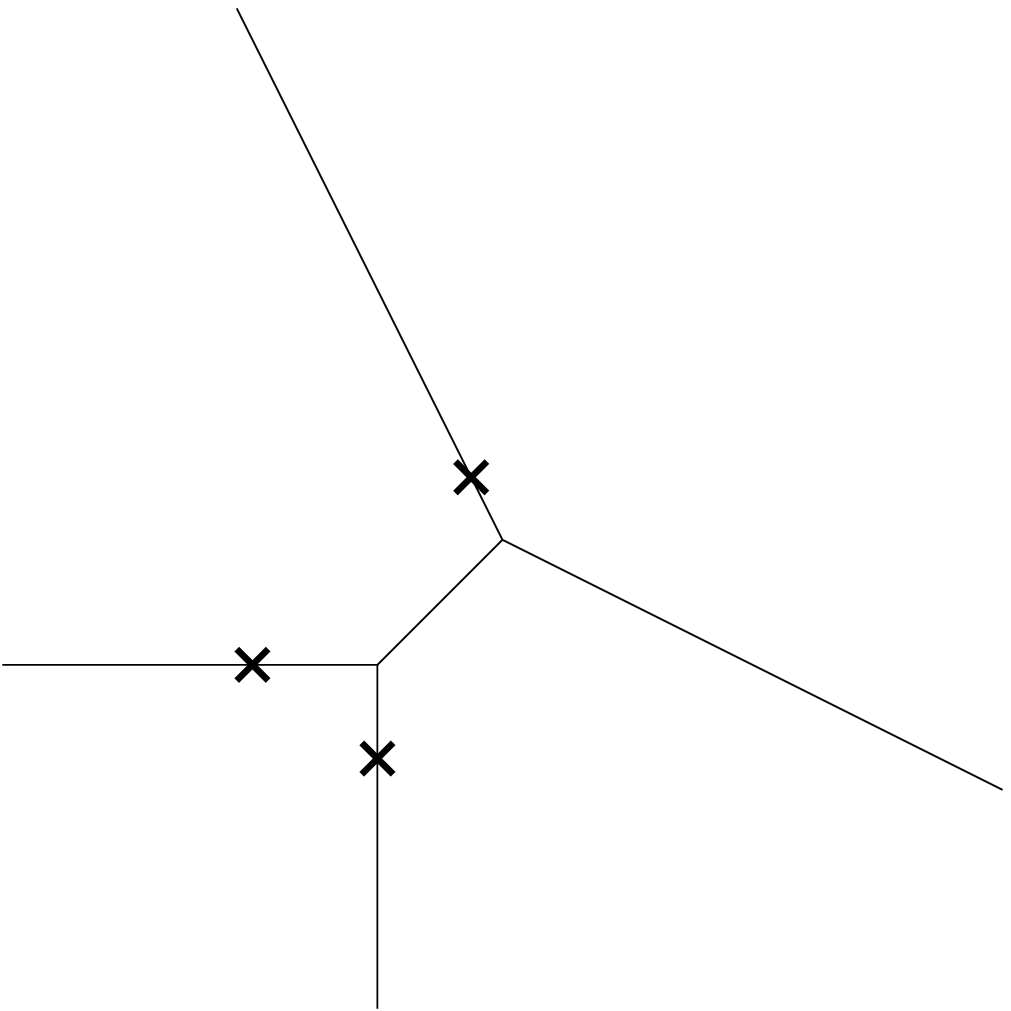,height=1.5in,width=1.5in}
}
\caption{\label{tren1} Tropical curves through a configuration of 3 points,
$\ntrop(g,\Delta)=5$.}
\end{figure}

\begin{figure}[h]
\centerline{\psfig{figure=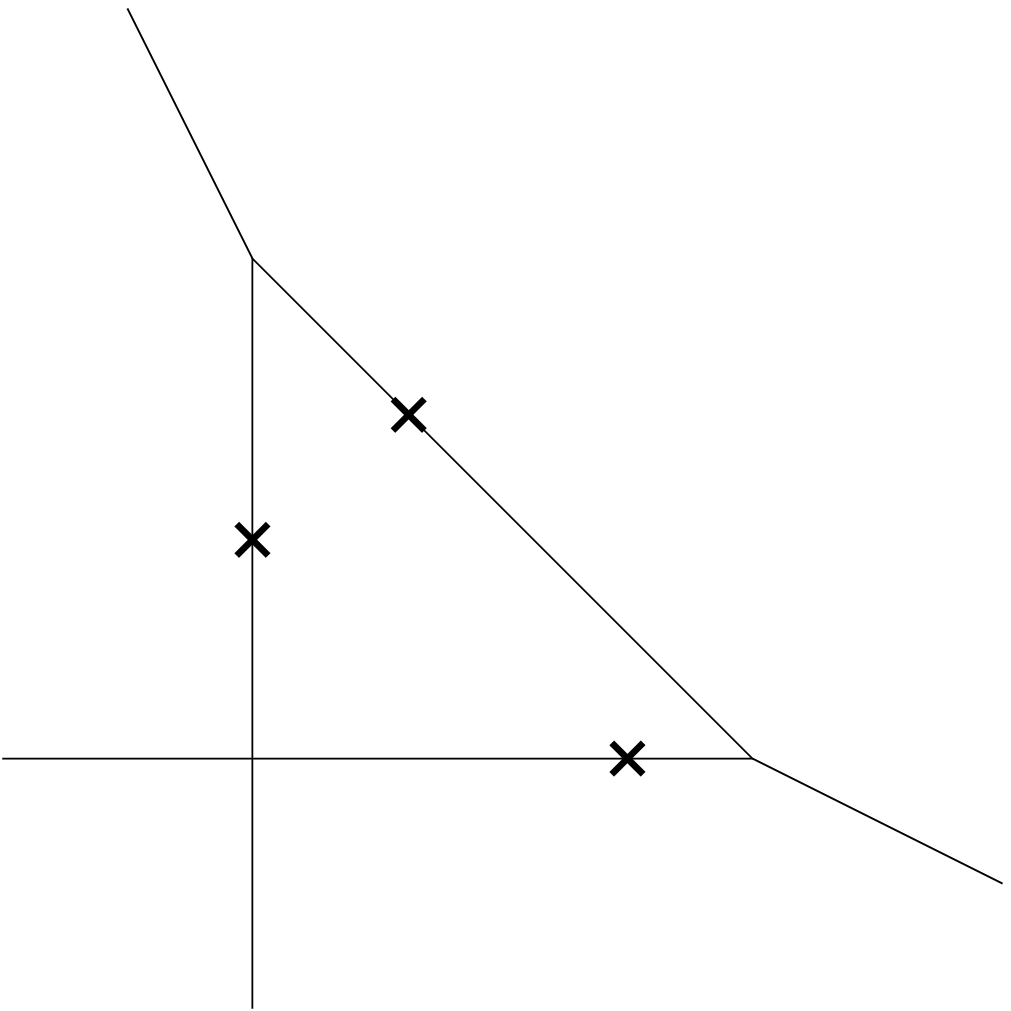,height=1.5in,width=1.5in}
\hspace{.1in}
\psfig{figure=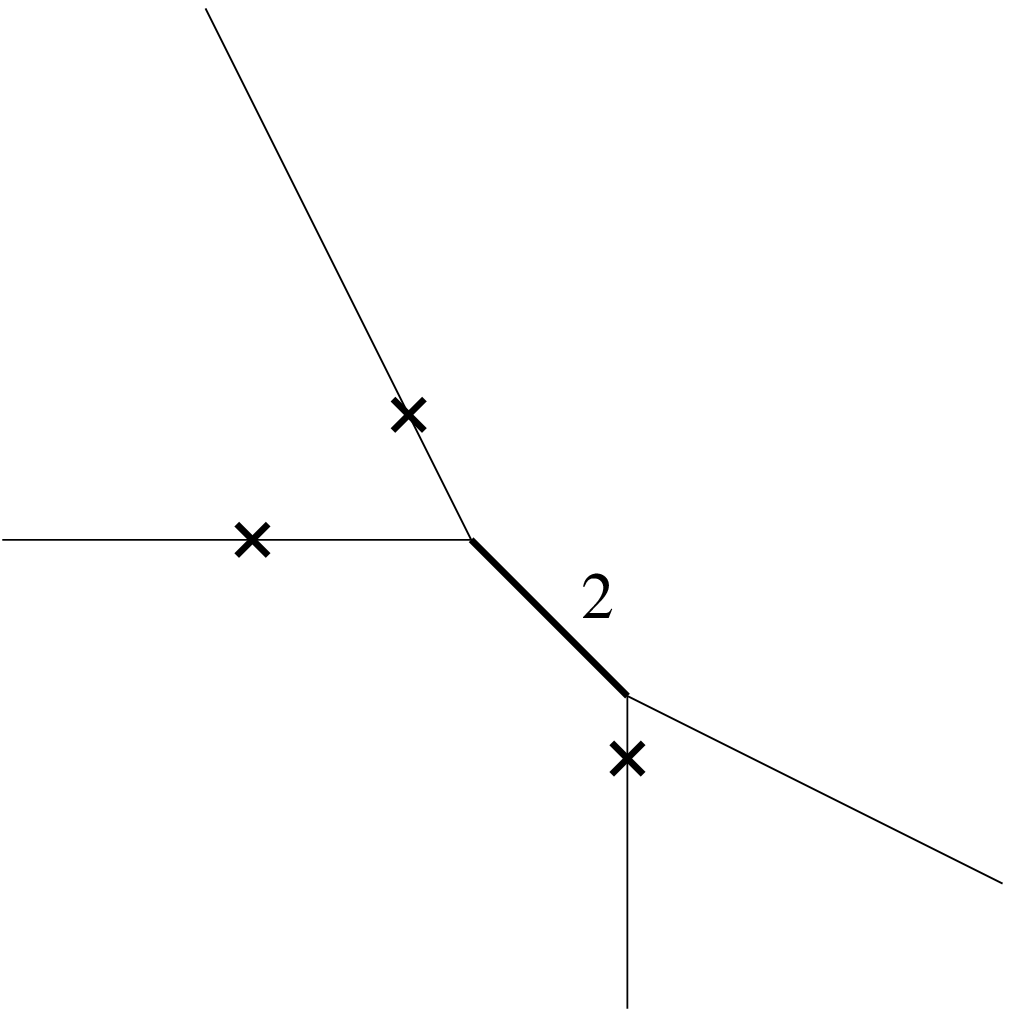,height=1.5in,width=1.5in}
}
\caption{\label{tren2} Tropical curves through another configuration of 3 points.
Note that the bounded edge in the right-hand curve has weight 2,
$\ntrop(g,\Delta)=5$.}
\end{figure}

\begin{exa}\label{exa-cusp}
Let $g=0$ and $\Delta$ be the quadrilateral whose
vertices are $(0,0)$, $(1,0)$, $(0,1)$ and $(2,2)$
(so that the number $s$ of the lattice points
on the perimeter $\dd\Delta$ is 4).
For a configuration $\ppp$ of 3 points in $\R^2$ pictured
in Figure \ref{tren1} we have 3 tropical curves passing.
In Figure \ref{tren2} the corresponding number is 2.
\end{exa}

\begin{defn}\label{multdim2}
The multiplicity $\mult(C)$ of a tropical curve $C\subset\R^2$ of degree $\Delta$ and
genus $g$ passing via $\ppp$ equals to the product of
the multiplicities of all the 3-valent vertices of $C$
(see Definition \ref{multvert}).
\end{defn}

\begin{defn}\label{tropnumbers}
We define the number $\ntrop(g,\Delta)$ to be the number of
irreducible tropical curves of genus $g$ and degree $\Delta$
passing via $\ppp$ where each such curve is counted with the multiplicity
$\mult$ from Definition \ref{multdim2}.
Similarly we define the number $\ntropr(g,\Delta)$ to be the number of
all tropical curves of genus $g$ and degree $\Delta$
passing via $\ppp$. Again each curve is counted with the multiplicity
$\mult$ from Definition \ref{multdim2}.
\end{defn}

The following proposition is a corollary of Theorem \ref{main}
formulated below in section \ref{sectionmain}.
\begin{prop}
The numbers $\ntropr(g,\Delta)$ and $\ntrop(g,\Delta)$ are finite
and do not
depend on the choice of $\ppp$.
\end{prop}
E.g. the 3-point configurations from
Figures \ref{tren1} and \ref{tren2} have the same number
$\ntrop(g,\Delta)$.

\subsection{Forests in the polygon $\Delta$}
\ignore{
Recall that we assume that our configuration $\ppp$ consists
of points $p_1,\dots,p_k\in\R^2$ in general position (cf. Corollary \ref{genposi}).
In addition, we assume that $k=m-l+g$, where $g$ is the genus of $C$
while $m$ and $l$ are the numerical characteristics of its Newton polygon
$\Delta$ from \eqref{ml}.
}

Recall that every vertex of a tropical curve $C\subset\R^2$ corresponds to a
polygon in the dual lattice subdivision of the Newton polygon
$\Delta$ while every edge of $C$ corresponds to an edge of the
dual subdivision $\subdiv_C$ (see Proposition \ref{dualsubdiv}).
If $C\subset\R^2$ is a tropical curve
passing through $\ppp$ then we can mark
the $k$ edges of $\subdiv_C$ dual to $p_1,\dots,p_k$.
Let $\Xi\subset\Delta$ be the union of the marked $k$ edges.
\begin{defn}\label{defcombtype}
The {\em combinatorial type} of the pair $(C,\ppp)$ passing
via the configuration $\ppp$ is the lattice subdivision $\subdiv_C$
together with the graph $\Xi\subset\Delta$ formed by
the $k$ marked edges of this subdivision.

The {marked combinatorial type} of a parameterized tropical
curve $h:\Gamma\to\R^2$ passing via $\ppp$ is the combinatorial
type of $h$ together with the marking of the edges containing
$h^{-1}(p_j)$.
\end{defn}
Note that these two notions are equivalent in the
case of simple curves.

\begin{prop}\label{forest}
The graph $\Xi\subset\Delta$ is a forest (i.e. a disjoint union of trees)
if $\ppp$ is in general position.
\end{prop}
\begin{proof}
From the contrary, suppose that $\Xi$ contains a cycle $Z\subset\Xi$
formed by $q$ edges. Since $\ppp$ is in general position our curve $C$
is simple.
Suppose that $p_1,\dots,p_q$ are the
marked points on the edges of $C$ dual to $Z$.
We claim that $p_1,\dots,p_q$ are {\em not} in tropical general position
which leads to a contradiction with Proposition \ref{trsubset}.
To show this we exhibit a curve of non-positive genus with $q$ ends
at infinity passing through $p_1,\dots,p_q$.


Consider the union $D$ of the (closed) polyhedra from $\operatorname{Subdiv}_C$
that are enclosed by $Z$. Let
$$C^D=\bigcup\limits_{\Delta'\subset D} U^{\Delta'},$$
where $U^{\Delta'}$ is the stratum of $C$ dual to $\Delta'$
(see Proposition \ref{dualsubdiv}). The set $C^D$ can be extended
to a tropical curve $\tilde{C}^D$ by extending all non-closed bounded
edges of $C^D$ to infinity. These extensions can intersect each other
if $D$ is not convex so the Newton polygon of $\tilde{C}^D$
is a convex polyhedron $\tilde\Delta$.
In other words,
$\tilde{C}^D$ is given by a tropical polynomial
$\tilde{f}^D(x)=\sum\limits_{j\in\tilde\Delta}b_jx^j$
with some choice of $b_j\in\Rtr$ (note that $b_j=a_j$ if $j\in D$).
The corresponding subdivision $\subdiv_{\tilde{f}^D}$
of $\tilde\Delta$ is $\subdiv_f|_D$ union with some parallelograms.

Consider the polynomial
$$\tilde{f}^Z(x)=
\sum\limits_{j\in\tilde\Delta\setminus\Int(\operatorname{Convex Hull}(D))}b_jx^j.$$
We have $V_{\tilde{f}^Z}\ni p_1,\dots,p_q$ since
$V_{\tilde{f}^D}\ni p_1,\dots,p_q$ and the edges corresponding
to $p_1,\dots,p_q$ are on the boundary of $D$.
On the other hand the genus of $V_{\tilde{f}^Z}$ is non-positive
since no vertex of $\operatorname{Subdiv}_{\tilde{f}^Z}$
is in the interior of $D$ (cf. Lemma \ref{genusvertex}).
\end{proof}

\begin{figure}[h]
\centerline{\psfig{figure=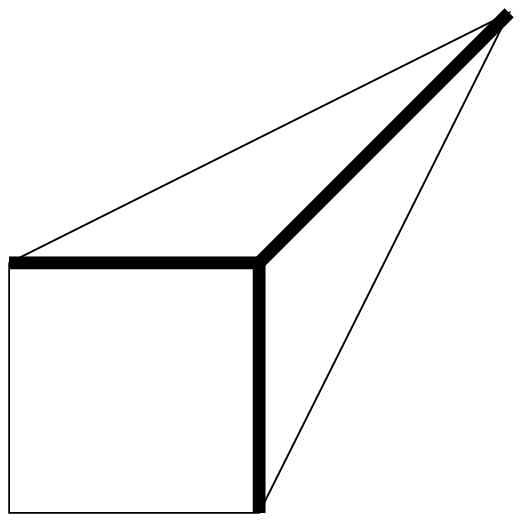,height=.7in,width=.7in}
\hspace{.1in}
\psfig{figure=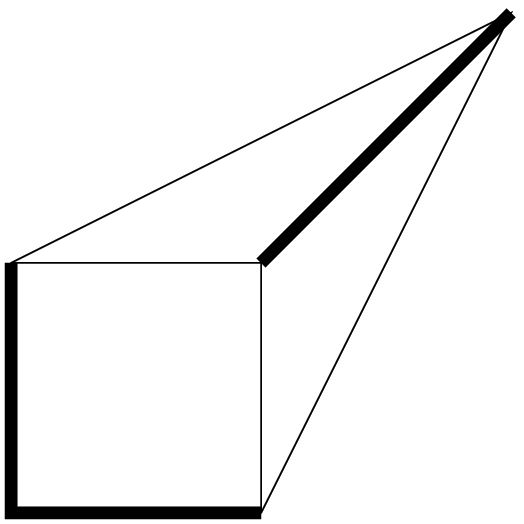,height=.7in,width=.7in}
\hspace{.1in}
\psfig{figure=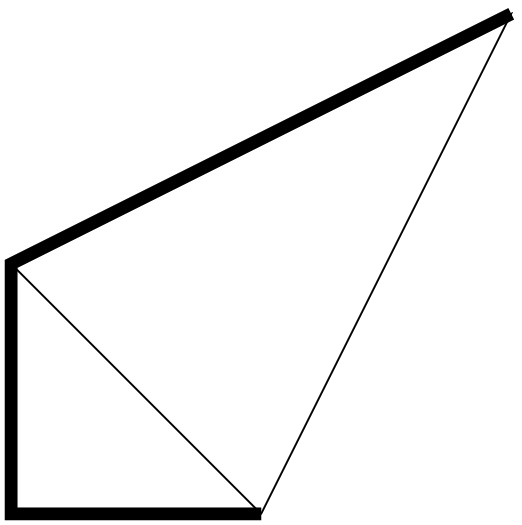,height=.7in,width=.7in}
\hspace{.1in}
\psfig{figure=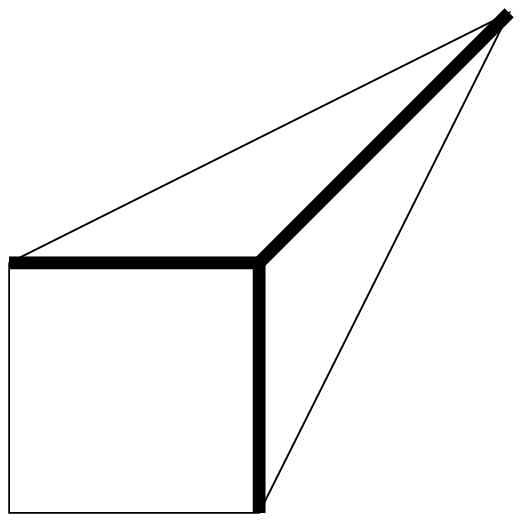,height=.7in,width=.7in}
\hspace{.1in}
\psfig{figure=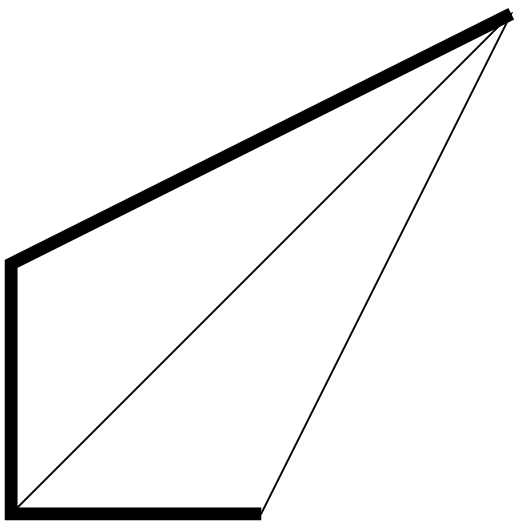,height=.7in,width=.7in}}
\caption{\label{forests} Forests corresponding to the curves
passing through the marked points from Figures \ref{tren1} and \ref{tren2}.}
\end{figure}

\subsection{The tropical curve minus the marked points}
The following lemma strengthens the Proposition \ref{forest} but
is stated in a dual language, in terms of the graph parameterizing
the tropical curve.
\begin{lem}\label{descompl}
Let $C$ be a simple curve of genus $g$ and degree $\Delta$.
Suppose that $C$ is parameterized by
$h:\Gamma\to\R^2$.
\begin{itemize}
\item Suppose that $C$ passes through a configuration $\ppp$ of $s+g-1$ points
in general position. Then each component $K$ of $\Gamma\setminus h^{-1}(\ppp)$
is a tree and the closure of $h(K)\subset\R^2$ has exactly one end of weight
one at infinity.
\item Conversely, suppose that $\ppp\subset C$ is a finite set
disjoint from the vertices of $C$ and
such that each component $K$ of $\Gamma\setminus h^{-1}(\ppp)$
is a tree while the closure of $h(K)\subset\R^2$ has exactly one end at infinity.
Then the combinatorial type of $(C,\ppp)$ is realized by $(C',\ppp')$,
where $C'$ is a curve of genus $g$ and $\ppp'$
is a configuration of points in tropical general position
which is a result of a small perturbation of $\ppp$.
Furthermore, the number of points in $\ppp$ is
$x+g-1$, where $x$ is the number of ends of $C$.
\end{itemize}
\end{lem}
\begin{proof}
Each component $K$ of $C\setminus\ppp$ in the first part of the statement
has to be a tree. Otherwise we can reduce the genus
of $C$ by the same trick as in Proposition \ref{forest} keeping the
number of ends of $C$ the same. This leads to a contradiction
with the assumption that $\ppp$ is in general position.
Also similarly to Proposition \ref{forest} we get a contradiction
if $h(K)$ is bounded.
If $h(K)$ has more than one end then by Theorem \ref{RR} it can
be deformed keeping the marked points from $\overline{h(K)}\setminus h(K)\subset\ppp$
fixed. This supplies a contradiction with Proposition \ref{tropfinite}.

For the second part let us slightly deform $\ppp$ to bring it to
a tropical general position. We can deform $h:\Gamma\to\R^2$
at each component $K$ individually to ensure $C'\supset\ppp'$.
\end{proof}

\begin{figure}[h]
\centerline{\psfig{figure=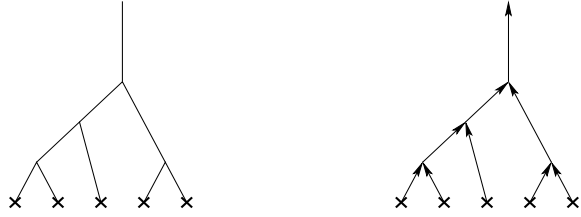,height=1.1in,width=3in}}
\caption{\label{typkomp} A component of
$\Gamma\setminus h^{-1}(\ppp)$ and its orientation.}
\end{figure}

Lemma \ref{descompl} allows one to extend the forest $\Xi$
from Proposition \ref{forest} to a tree $\XX\subset\Delta$
that spans all the vertices of $\subdiv_C$ in the case when
the number of points in the configuration $\ppp$ is $s+g-1$.

Each parallelogram $\Delta'$ corresponds to an intersection
of two edges $E$ and $E'$ of $\Gamma$.
Lemma \ref{descompl} allows one to orient the edges of each component
of $\Gamma\setminus h^{-1}(\ppp)$ consistently toward the end at infinity
(see Figure \ref{typkomp}).
Let $U'\subset\R^2$ be a small disk with a center at $E\cap E'$.
Each component of $U'\setminus (E\cup E')$ corresponds to a vertex
of $\Delta'$. Two of these components are distinguished
by the orientations of $E$ and $E'$.
One is adjacent to the sources of the edges while the other is adjacent
to the sinks. We connect the corresponding vertices of $\Delta'$
with an edge.

\begin{figure}[h]
\centerline{\psfig{figure=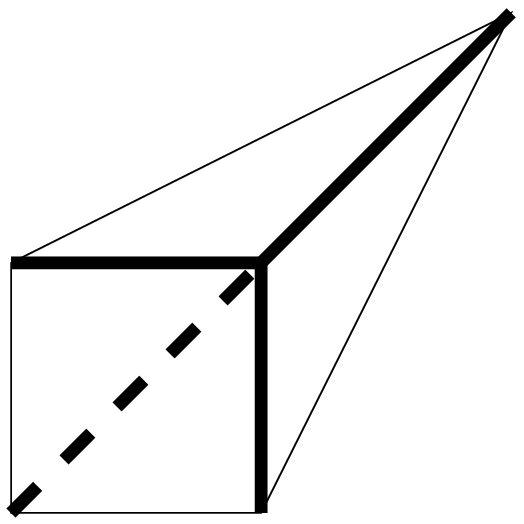,height=.7in,width=.7in}
\hspace{.1in}
\psfig{figure=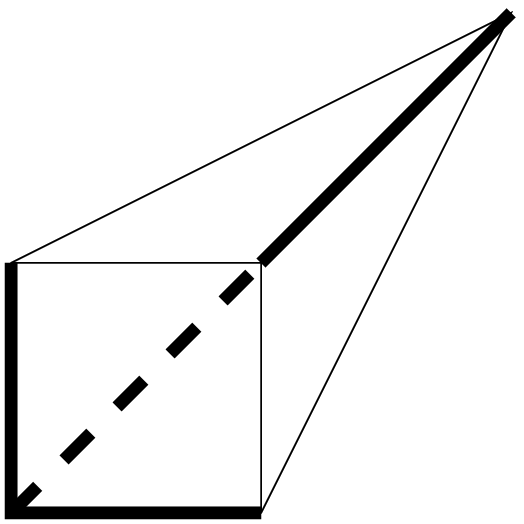,height=.7in,width=.7in}
\hspace{.1in}
\psfig{figure=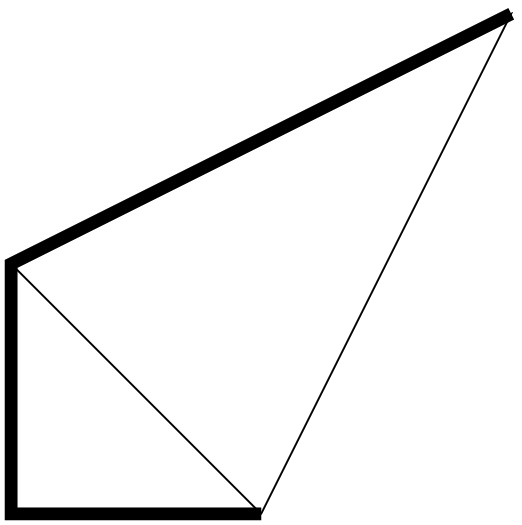,height=.7in,width=.7in}
\hspace{.1in}
\psfig{figure=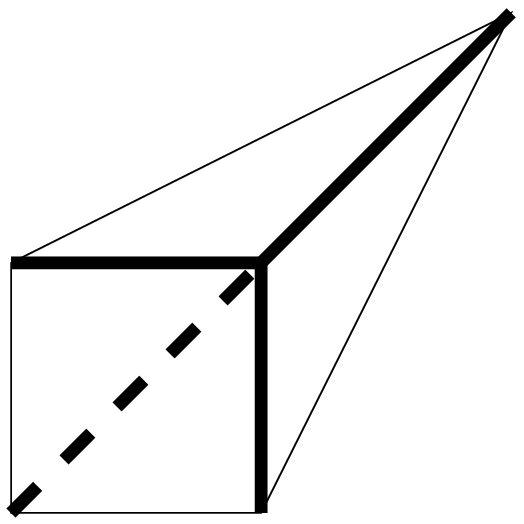,height=.7in,width=.7in}
\hspace{.1in}
\psfig{figure=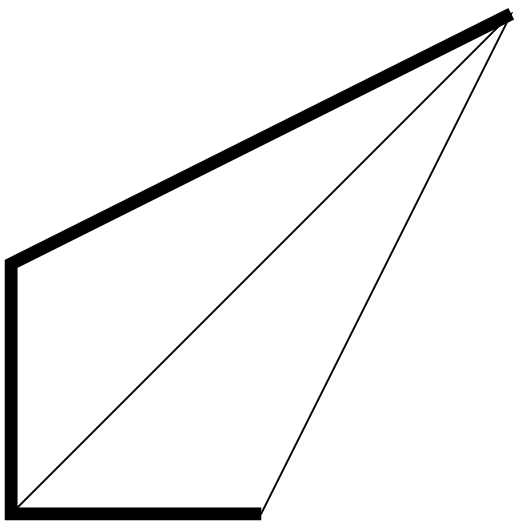,height=.7in,width=.7in}}
\caption{\label{trees} Trees obtained from the forests
in Figure \ref{forests}.}
\end{figure}

We form the graph $\XX\subset\Delta$ by taking the union of $\Xi$
with such edges for all the parallelograms from $\subdiv_C$,
see Figure \ref{trees}.

\begin{prop}\label{tree}
The graph $\XX$ is a tree that contains all the vertices of $\subdiv_C$.
\end{prop}
\begin{proof}
Suppose that $\subdiv_C$ does not contain parallelograms.
Then $\XX=\Xi$. Let $K\subset\XX$ be a component of $\XX$.
If there exists a vertex $v\in\subdiv_C$ not contained in $K$
then we can form a 1-parametric family of curves of genus $g$ and
degree $\Delta$ passing via $\ppp$. Indeed,
let $f(x)=``\sum \beta_j x^j$ be the tropical polynomial that
defines $C$. To get rid of the ambiguity in the choice of $f$
we assume that $j$ runs only over the vertices of $\subdiv_C$ and
that $\beta_{j_0}=0$ for a choice of the base index
$j_0\in K$.
Let us deform the coefficient $\beta_v$.
If $v$ belongs to component $K'$ of $\Xi$ different from $K$
then we also inductively deform coefficients at the other
vertices of $\subdiv_C$ that belong to $K'$ to make sure
that the curve $C'$ corresponding to the result of deformation
still contains $\ppp$. Clearly the genus of $C'$ is still $g$.
Thus we get a contradiction to the assumption that $\ppp$ is in
tropically general position.


We can reduce the general case to this special case
by the following procedure.
For each parallelogram $\Delta'\in\subdiv_C$ consider a point $p'$
and its small neighborhood $U'\subset\R^2$ that
is obtained by a small shift of the intersection $E\cap E'$ to the
target component of $U'\setminus (E\cup E')$
where $U'\subset\R^2$ is a small neighborhood of $p'$.
Let $\ppp'\supset\ppp$ be the resulting configuration.
The curve $C$ can be deformed to a curve $C'\supset\ppp'$
with the corresponding forest $\Xi'$ equal to $\XX$.
\end{proof}

\ignore{
\begin{coro}\label{coro-forest}
If $\ppp$ consists of $s+g-1$ points then the forest $\Xi$
from Proposition \ref{forest} contains all the vertices of $\subdiv_C$.
\end{coro}
}

\subsection{Uniqueness in a combinatorial type}
Enumeration of tropical curves is easier than that of complex
or real curves thanks to the following lemma.
\begin{lem}\label{1comb}
In each combinatorial type of marked tropical curves of genus $g$
with $x$ ends
there is either one or no curves passing through
$\ppp\in\R^2$ as long as $\ppp$ is a configuration of $x+g-1$
points in general position.
\end{lem}

\begin{proof}
The tropical immersions $h:\Gamma\to\R^2$ of a given
combinatorial type form a convex polyhedral domain $P\subset\R^{x+g-1}$
by Propositions \ref{RR} and \ref{3vimm}.
Since $\ppp$ is in general position $h$ is a simple curve.
The condition that the image of a particular edge of $\Gamma$ contains
$p_j$ is a hyperplane in $\R^{x+g-1}$.

Suppose that these conditions
cut a positive-dimensional polyhedral domain $Q\subset P$.
A point from $Q\cap\dd P$ (the boundary is taken in $\R^{x+g-1}$)
is a non-simple curve with perhaps even smaller $g$ and $x$ and
thus cannot pass through $\ppp$.
If $Q\cap \dd P=\emptyset$ but $Q\subset P$ is positive-dimensional
then it contains a line. Since the edge lengths cannot be negative
this means that this line corresponds to a family obtained from
a single curve by translations. This supplies a contradiction.

\ignore{
If $Q\cap \dd P=\emptyset$
then we can go to infinity along $Q$ and then by the compactness
theorem \ref{compact} the corresponding sequence of curves
converges to a tropical curve $C'$ parameterized by a part $\Gamma'$ of
the graph $\Gamma$.
The complement $\Gamma-\Gamma'$ parameterizes a tropical curve
in the complement of a large radius disk $D\subset\R^2$.
Each component of $\Gamma-\Gamma'$ has ends of two types:
going to infinity in $\R^2$ or ending on $\dd D$.

Thus existence of $C'$ contradicts to the general position of $\ppp$.
}

Note that even if the combinatorial type is generic we still may
have no curves passing through $\ppp$ since this linear
system of equations defined not in the whole $\R^n$ but in an open
polyhedral domain there.
\end{proof}


\section{Algebraic curves in $\tor$ and a classical enumerative problem}
\subsection{Enumerative problem in $\tor$}
As in subsection \ref{trenpr} we fix a number $g\in\Z$
and a convex lattice polygon $\Delta\subset\R^2$.
As before let $s=\#(\dd\Delta\cap\Z^2)$.
Let $\qq=\{q_1,\dots,q_{s+g-1}\}\in\tor$ be a configuration
of points in general position.
A complex algebraic curve $C\subset\tor$ is defined
by a polynomial $f:\tor\to\C$ with complex coefficients.
As in the tropical set-up we refer to the Newton polygon $\Delta$
of $f$ as the {\em degree of $C$}.
\begin{defn}
We define the number $\ncl(g,\Delta)$ to be the number of
irreducible complex curves of genus $g$ and degree $\Delta$
passing via $\qq$.
Similarly we define the number $\nclr(g,\Delta)$ to be the number of
all complex curves of genus $g$ and degree $\Delta$
passing via $\qq$.
\end{defn}
Note that here we count every relevant complex curve simply,
i.e. with multiplicity 1.
\begin{prop}\label{cgeneric}
For a generic choice of $\qq$ the numbers $\ncl(g,\Delta)$
and $\nclr(g,\Delta)$ are finite and do not depend on $\qq$.
\end{prop}
This proposition is well-known, cf. e.g. \cite{CH} (also later
in this section the invariance of $\ncl(g,\Delta)$
and $\nclr(g,\Delta)$ is reduced to invariance
of certain Gromov-Witten numbers).

In modern Mathematics there are two ways to interpret the numbers
$\ncl(g,\Delta)$ and $\nclr(g,\Delta)$.
A historically older interpretation is via the degree of Severi varieties.
A more recent interpretation (introduced in \cite{KM}) is via the
Gromov-Witten invariants. In both interpretations
it is convenient to consider the compactification of the problem
with the help of the toric surface associated to the polygon $\Delta$.

\subsection{Toric surfaces and Severi varieties}
\begin{figure}[h]
\centerline{\psfig{figure=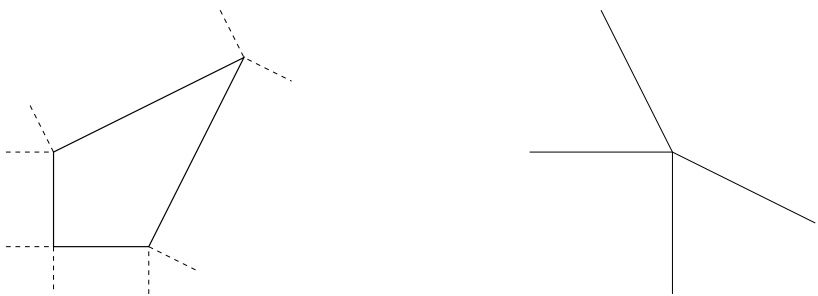,height=1.5in,width=4.25in}}
\caption{\label{fan} A polygon and its normal fan.}
\end{figure}

Recall that a convex polygon $\Delta$ defines a compact toric
surface $\C T_\Delta\supset\tor$, see e.g. \cite{GKZ}. (Some
readers may be more familiar with the definition of toric surfaces
by fans, in our case the fan is formed by the dual cones at the
vertices of $\Delta$, see Figure \ref{fan}.)
The sides of the polygon $\Delta$
correspond to the divisors in $\C T_\Delta\setminus\tor$. These
divisors intersect at the points corresponding to the vertices of
$\Delta$. This surface is non-singular if every vertex of $\Delta$
is {\em simple}, i.e. its neighborhood in $\Delta$ is mapped to a
neighborhood of the origin in the positive quadrant angle under a
composition of an element of $SL_2(\Z)$ and a translation in
$\R^2$. Non-simple vertices of $\Delta$ correspond to
singularities of $\C T_\Delta$.

\begin{exam}\label{mn-ki}
Let $\Delta_d$ be the convex hull of $(d,0)$, $(0,d)$ and $(0,0)$.
We have $\C T_{\Delta_d}=\cp^2$ no matter what is $d$.
If $\Delta=[0,r]\times[0,s]$, $r,s\in\N$ then $\C
T_{\Delta}=\cp^1\times\cp^1$ no matter what are $r$ and $s$. All
the vertices of such polygons are simple.

If $\Delta^{\frac32}=\operatorname{Convex
Hull}\{(0,0),(2,1),(1,2)\}$ then $\C T_\Delta=\cp^2/\Z^3$, where
the generator of $\Z^3$ acts on $\cp^2$ by $[x:y:z]\mapsto
[x:e^{\frac{2\pi i}{3}}y:e^{\frac{4\pi i}{3}}z]$. This action has
3 fixed points which give the singularities of $\C
T_{\Delta^{\frac32}}$ corresponding to the three (non-simple)
vertices of $\Delta$.
\end{exam}

In addition to a complex structure (which depends only on the dual
fan) the polygon $\Delta$ defines a holomorphic linear bundle
$\HH$ over $\C T_\Delta$. Let $\LLL_\Delta=\Gamma(\HH)$ be the
vector space of sections of $\HH$. The projective space
$\pp\LLL_\Delta$ is our system of the curves. Note that it can be
also considered as the space of all holomorphic curves in $\C
T_\Delta$ such that their homology class is Poincar\'e dual to
$c_1(\HH)$.

Returning to Example \ref{mn-ki} we note that $\Delta_d$ gives us
the projective curves of degree $d$. The polygon
$[0,r]\times[0,s]$ gives us the curves of bidegree $(r,s)$ in the
hyperboloid $\cp^1\times\cp^1$. The polygon $\Delta^{\frac32}$
gives us the images in $\C T_{\Delta^{\frac32}}$ of the cubic
curves in $\cp^2$ that are invariant with respect to the
$\Z_3$-action.

Any curve in $\LLL_\Delta$ is the closure in $\C T_\Delta$ of the
zero set of a polynomial whose Newton polygon is contained in
$\Delta$. Thus $\dim\LLL_\Delta=\#(\Delta\cap\Z^2)=s+l$
(see \eqref{ml}) and
$\dim\pp\LLL_\Delta=s+l-1$. A general curve $C\subset\C T_\Delta$ from
$\pp\LLL_\Delta$ is a smooth curve that is transverse to
$\C T_\Delta\setminus\tor$
(this means that it is transverse to all divisors
corresponding to the sides of $\dd\Delta$ and does not pass
through their intersection points).

By the genus formula we have $g(C)=l=\#(\Int\Delta\cap\Z^2)$ for a
smooth curve $C$ in $\pp\LLL_\Delta$. However singular curves have
smaller geometric genus. More precisely let $C\subset\C T_\Delta$
be the curve from $\pp\LLL_\Delta$ and let $\tilde{C}\to C$ be its
normalization. We define {\em the geometric genus} as
$g(C)=\frac12 (2-\chi(\tilde{C}))$. Note that if $C$ is not
irreducible then $\tilde{C}$ is disconnected and then $g(C)$ may
take a negative value.

Fix a number $g\in\Z$. The curves of genus not greater than $g$
form in the projective space $\pp\LLL_\Delta$
an algebraic variety known as {\em the Severi variety} of $\C T_\Delta$.
This variety may have several components.
E.g. if $\C T_\Delta$ has an exceptional divisor $E$
corresponding to a side $\Delta_E\subset\Delta$ then reducible
curves $E\subset C'$ where $C'$ corresponds to the polygon
$\Delta'=\operatorname{Convex Hull}((\Delta\setminus\Delta_E)\cap\Z^2)$
form a component (or a union of components)
of the Severi variety. Such components correspond to
smaller polygons $\Delta'\subset\Delta$.
We are interested only in those components that correspond to
the polygon $\Delta$ itself.
\begin{defn}
{\em The irreducible Severi variety
$\Sigma^{\operatorname{irr}}_g\subset\pp\LLL_\Delta$ corresponding to $\Delta$
of genus $g$} is the closure of the set formed by all irreducible
curves whose Newton polygon is $\Delta$ and genus is not more then $g$.
{\em The Severi variety $\Sigma_g\subset\pp\LLL_\Delta$
corresponding to $\Delta$ of genus $g$} is the closure of the set formed by all
curves whose Newton polygon is $\Delta$ and genus is not more then $g$.
Clearly, $\Sigma_g\supset\Sigma^{\operatorname{irr}}_g$.
\end{defn}

Note that $\Sigma_g$ is empty unless $g\le l$. If $g=l$ we
have $\Sigma_g=\pp\LLL$. If $g=l-1$ then $\Sigma_g$
is the (generalized) $\Delta$-discriminant variety.
It is sometimes convenient to set $\delta=l-g$.
Similarly to \cite{CH} it can be shown that $\Sigma_g$ is
the closure in $\pp\LLL$ of immersed nodal curves with $\delta$
ordinary nodes.
In the same way, $\Sigma^{\operatorname{irr}}_g$ is
the closure in $\pp\LLL$ of irreducible immersed nodal
curves with $\delta$ ordinary nodes.

It follows from the Riemann-Roch formula that $\Sigma_g$ and
$\Sigma^{\operatorname{irr}}_g$ have pure dimension $s+g-1$.
The Severi numbers $\ncl(g,\Delta)$ and $\nclr(g,\Delta)$
can be interpreted as the degrees of
$\Sigma^{\operatorname{irr}}_g$ and $\Sigma_g$ in $\pp\LLL$.

\begin{exam}
Suppose $\Delta=\Delta_d$ so that $\C T_\Delta=\cp^2$.
We have $s=3d$ and $l=\frac{(d-1)(d-2)}{2}$. The
number $\ncl(g,\Delta)$ is the number of genus $g$,
degree $d$ (not necessarily irreducible) curves passing through
$3d+g-1$ generic points in $\cp^2$.

%
The formula $\nclr(l-1,\Delta_d)=3(d-1)^2$
is well-known as the degree of the discriminant (cf. \cite{GKZ}).
(More generally, if $\C T_\Delta$ is smooth then
$\nclr(l-1,\Delta)=6\Area(\Delta)-2\operatorname{Length}(\dd\Delta)+
\#\operatorname{Vert}\Delta$, where $\operatorname{Length}(\dd\Delta)=s$
is the lattice length of $\dd\Delta$ and $\#\operatorname{Vert}\Delta$
is the number of vertices, see \cite{GKZ}).

An elegant recursive formula
for $\ncl(0,\Delta_d)$ was found by Kontsevich \cite{KM}.
Caporaso and Harris \cite{CH} discovered an algorithm
for computing $\nclr(g,\Delta_d)$ for arbitrary $g$.
See \cite{V} for computations for some other
rational surfaces, in particular, the Hirzebruch surfaces (this
corresponds to the case when $\Delta$ is a trapezoid).
\end{exam}

\subsection{Gromov-Witten invariants}
If $\Delta$ is a polygon with simple vertices then $\C
T_\Delta$ is a smooth 4-manifold.
This manifold is equipped with
a symplectic form $\omega_\Delta$ defined by $\Delta$. The linear
system $\LLL_\Delta$ gives an embedding $\C
T_\Delta\subset\pp\LLL_\Delta\approx\cp^m$. This embedding induces
$\omega_\Delta$. As we have already seen $\Delta$ also defines a
homology class $\beta_\Delta\in H_2(\C T_\Delta)$, it is the
homology class of the curves from $\pp\LLL_\Delta$.

To define the Gromov-Witten invariants
of genus $g$ one takes a generic almost-complex structure
on $\C T_\Delta$ that is compatible with $\omega_\Delta$ and counts
the number of pseudo-holomorphic curves of genus $g$ via generic
$s+g-1$ points in $\C T_\Delta$ in the following sense
(see \cite{KM} for a precise definition).
\footnote{These are the Gromov-Witten invariants
evaluated on the cohomology classes dual to a point, this is the
only non-trivial case for surfaces. In this discussion we
completely ignore the gravitational descendants.}
Consider the space $\M_{g,s+g-1}(\C T_\Delta)$
of all stable (i.e. those with finite automorphism group)
parameterized pseudo-holomorphic curves with $s+g-1$ marked points.
Evaluation at each marked point produces a map
$\M_{g,s+g-1}(\C T_\Delta)\to\C T_\Delta$.
With the help of this map
we can pull back to $\M_{g,s+g-1}(\C T_\Delta)$ any cohomology
class in $\C T_\Delta$, in particular the cohomology class of a point.
Doing so for each of the $s+g-1$ point and taking the cup-product
of the resulting classes we get the Gromov-Witten invariant
$I^{\C T_\Delta}_{g,s+g-1,\beta_\Delta}<pt^{\otimes^{s+g-1}}>$.

The result is invariant with respect to deformations of the almost-complex
structure.
In many cases it is useful to pass to a generic almost-complex structure to make
sure that for any stable curve $C\to\C T_\Delta$ passing through
our points we have $H^1(C,N_{C/\C T_\Delta})=0$. But we have this
condition automatically if $\C T_\Delta$ is a smooth Fano surface
(or, equivalently, all exceptional divisors have self-intersection
$-1$), cf. e.g. \cite{V} for details. (However, if $\C T_\Delta$ has
exceptional divisors of self-intersection $-2$ and less we need
either to perturb the almost-complex structure or to consider a
virtual fundamental class.)

The Gromov-Witten invariant
$I^{\C T_\Delta}_{g,s+g-1,\beta_\Delta}<pt^{\otimes^{s+g-1}}>$
coincides with the number $\ncl(g,\Delta)$ if $\C T_\Delta$ is
a toric Fano surface (cf. e.g. \cite{V}).
In particular this is the case for $\cp^2$ or $\cp^1\times\cp^1$
(these are the only smooth toric surfaces without exceptional
divisors, in other words minimal Fano).

The Gromov-Witten invariants corresponding to disconnected curves
are sometimes called
multicomponent Gromov-Witten invariants.

\section{Complex tropical curves in $\tor$ and
connection between classical and tropical geometries}
\subsection{Degeneration of complex structure on $\tor$}\label{defstr}

Let $t>1$ be a real number. We have the
following self-diffeomorphism $\tor\to\tor$
\begin{equation}\label{Ht}
H_t:(z,w)\mapsto
(|z|^{\frac{1}{\log t}}\frac{z}{|z|},|w|^{\frac{1}{\log t}}\frac{w}{|w|}).
\end{equation}
For each $t$ this map induces a new complex structure on $\tor$.

Here is a description of the complex structure induced by $H_t$ in
logarithmic polar coordinates $\tor\approx \R^2\times iT^2$. (This
identification is induced by the holomorphic logarithm $\LL$ from
the identification $\C^2\approx\R^2\times i\R^2$.) If $v$ is a
vector tangent to $iT^2$ we set $J_tv=\frac{1}{\log(t)}iv$. Note that
$J_tv$ is tangent to $\R^2$.

Clearly, a curve $V_t$ is holomorphic with respect to $J_t$ if and
only if $V_t=H_t(V)$, where $V$ is a holomorphic curve with
respect to the standard complex structure, i.e. $J_e$-holomorphic.
Let
$\Log:\tor\to\R^2$
be the map defined by $\Log(z,w)=(\log|z|,\log|w|)$.
We have $$\Log\circ H_t=\Log_t.$$
Note that
$H_t$ corresponds to a $\log(t)$-contraction $(x,y)\mapsto
(\frac{x}{\log(t)},\frac{y}{\log(t)})$ under $\Log$.

\subsection{Complex tropical curves in $\tor$}
There is no limit (at least in the usual sense) for the complex
structures $J_t$, $t\to\infty$. Nevertheless, as in section 6.4 of
\cite{M-pp} we can define the $J_\infty$-holomorphic curves
which happen to be the limits of
$J_t$-holomorphic curves, $t\to\infty$.

There are several way to define them. An algebraic definition is
the shortest and involves varieties over a non-Archimedean field.
Let $K$ be the field of the (real-power) Puiseux series
$$a=\sum\limits_{j\in I_a}\xi_j t^j,\ \xi_j\in\C,$$
where $I_a\subset\R$ is a well-ordered set (cf. \cite{Ka}).
The field $K$ is algebraically closed and of
characteristic 0. The field $K$ has a non-Archimedean valuation
$\val(a)=-\min I_a$, $$\val(a+b)\le\max\{\val(a),\val(b)\}.$$

As usual, we set $K^*=K\setminus\{0\}$. The multiplicative
homomorphism $\val:K^*\to\R$ can be ``complexified" to
\begin{equation}\label{mapw}
w:K^*\to\C^*\approx\R\times S^1
\end{equation}
by setting
$w(a)=e^{\val(a)+i\arg(a_{\val(a)})}$. Applying this map
coordinatewise we get the map $$W:(K^*)^2\to\R^2\times (S^1)^2\approx\tor.$$
Applying the map $\val$ coordinatewise we get the map
$$\Val:(K^*)^2\to\R^2,$$
$\Val=\Log\circ W$.
The image of an algebraic curve $V_K$ under $W$ turns out to be a
$J_\infty$-holomorphic curve (cf. Proposition \ref{jinfty}).

Note the following special case. Let
$$V_K=\{z\in \tor\ |\ \sum\limits_{j\in \operatorname{Vert}\Delta}
a_jz^j=0\},$$ where $z\in(K^*)^2$, $j$ runs over all vertices of
$\Delta$ and $a_j$ is such that $\val(a_j)=0$. By Kapranov's
theorem \cite{Ka} $\Val(V_K)=\Log(W(V_K))$ is the tropical curve
$C_f$ defined by $f(x)=``\sum\limits_{j\in
\operatorname{Vert}\Delta}x^j"$, $x\in\R^2$.
Thus $C_f$ is a union of rays starting from the origin and
orthogonal to the sides of $\Delta$, in other words it is the
1-skeleton of the {normal fan} to $\Delta$.
However, $W(V_K)\subset\tor$ depends on the argument of the
leading term of $a_j\in K^*$. Thus different choices of $W(V_K)$
give different {\em phases} for the lifts of $C_f$. We may
translate $C_f$ in $\R^2$ so that it has a vertex in a point
$x\in\R^2$ instead of the origin. Corresponding translations of
$W(V_K)$ give a set of possible lifts.

This allows one to give a more geometric description of
$J_\infty$-curves. They are certain 2-dimensional objects in
$\tor$ which project to tropical curves under $\Log$. Namely, let
$C\subset\R^2$ be a tropical curve, $x\in C$ be a point
and $U\ni x\in\R^2$ be a convex neighborhood such that
$U\cap C$ is a cone over $x$ (i.e. for every $y\in C\cap
U$ we have $[x,y]\subset C$). Note that if $x$ is a point on
an open edge then it is dual (in the sense of \ref{Dsubdiv}) to a
segment in $\Delta$. If $x$ is a vertex of $ C$ then it is
dual to a 2-dimensional polygon in $\Delta$. In both cases we
denote the dual polygon with $\Delta'$. We say that a
2-dimensional polyhedron $V_\infty\subset\tor$ is {\em
$( C\cap U)$-compatible} if $\Log^{-1}(U)\cap V_\infty=
\Log^{-1}(U)\cap W$, where $W$ is a translation of $W(V_K)$ while
$V_K=\{z\in \tor\ |\ \sum\limits_{j}
a_jz^j\}$, $j$ runs over some lattice points of $\Delta'$
and $\val(a_j)=0$.

\begin{prop}\label{jinfty}
Let $V_\infty\subset\tor$. The
following conditions are equivalent.
\begin{enumerate}
\item $V_\infty=W(V_K)$, where $V_K\subset (K^*)^2$ is an
algebraic curve.
\item $C=\Log(V_\infty)\subset\R^2$
is a graph such that there exists a choice of
natural weights on its edges turning it to a tropical curve
such that for every $x\in C$ there exists a small open convex
neighborhood $x\in U\subset\R^2$ such that $\Log^{-1}(U)\cap
V_\infty$ is $( C\cap U)$-compatible.
\item $V_\infty$ is the
limit when $k\to\infty$ in the Hausdorff metric of a sequence of
$J_{t_k}$-holomorphic curves $V_{t_k}$ with
$\lim\limits_{k\to\infty}t_k=\infty$.
\end{enumerate}
\end{prop}
We precede the proof of Proposition \ref{jinfty}
with some definitions and remarks.
\begin{defn}
Curves satisfying any of the equivalent conditions of Proposition \ref{jinfty}
are called (unparameterized) {\em complex tropical curves}
or {\em complex tropical 1-cycles in $\tor$}.
\end{defn}
Proposition \ref{jinfty} allows us to think of complex tropical curves
both as tropical curves equipped with a {\em phase}, i.e.
a lifting to $\tor$, and as $J_\infty$-holomorphic curves,
i.e. as limits of $J_t$-holomorphic curves when $t\to\infty$.

\begin{prop}\label{ctrop-trop}
Let $V_\infty=W(V_K)$, where $V_K\subset(K^*)^2$ is an algebraic curve
with the Newton polygon $\Delta$.
Then
$\Log(V_\infty)\subset\R^2$ is a graph. Furthermore,
it is possible to equip the edges of $\Log(V_\infty)$
with natural weight so that the result is a tropical
curve of degree $\Delta$ in $\R^2$.
\end{prop}
\begin{proof}
The proposition follows from Kapranov's theorem \cite{Ka},
cf. \ref{Kap} since $\log\circ W=\val:K^*\to\R$.
The edge weights come from the lattice lengths of the edges
of the corresponding lattice subdivision of $\Delta$.
\end{proof}
\begin{defn}
We say that a complex tropical curve $V_\infty$
with a choice of natural weights for the edges
of $C=\Log(V_\infty)$ has degree $\Delta$ if
these weights turn $C$ to a tropical curve of degree $\Delta$.
\end{defn}

\begin{proof}[Proof of Proposition \ref{jinfty}]

{$(1)\implies (2).$}
Let $f:(K^*)^2\to K$ be the polynomial defining $V_K$. The image
$C=\Log(V_\infty)$ is a tropical curve by Proposition \ref{ctrop-trop}.
Let $x\in C\subset\R^2$. The lowest $t$-powers
of $f(x)$ come only from the polygon $\Delta'\subset\Delta$ dual
to the stratum containing $x$. The compatibility curve is given by
the sum of the $\Delta'$-monomials of $f$.

$(2)\implies (1).$ Consider the subdivision of $\Delta$ dual to the
tropical curve $ C$. The compatibility condition gives us a
choice of monomials for each polygon $\Delta'\subset\Delta$ in the
subdivision dual to $ C$. However, the choice is not unique,
due to the higher $t$-power contributions. On the other hand, a
monomial corresponds to a lattice point of $\Delta$ which may
belong to several subpolygons in the subdivision.

We have to choose the coefficients for the monomials so that they
would work for all subpolygons of the subdivision. Let
$j\in\Delta\cap\Z^2$. The coefficient $a_j$ is a Puiseux series in
$t$. The lowest power $\val(a_j)$ is determined from the
tropical curve $C$ as the coefficient of
the corresponding tropical monomial.
We set $a_j=\alpha_jt^{\val(a_j)}$, where
$\alpha_j\in\C^*$, i.e. our coefficient Puiseux series are
actually monomials.

Namely, let $\Delta'$ be a polygon in the subdivision of $\Delta$
dual to $ C$. A point $x$ on the corresponding stratum of
$ C$ is compatible with $W(\{f^{\Delta'}=0\})$ for a
polynomial $f^{\Delta'}$ over $K$ with the Newton polygon
$\Delta'$. \ignore{ Let $\Delta'=[j,k]$ be an edge of the
subdivision of $\Delta$ dual to $ C$, $j,k\in\Delta\cap\Z^2$.
A point $x$ on the corresponding edge of $ C$ is compatible
with $W(\{f^{\Delta'}=0\})$ for a polynomial $f^{\Delta'}$ over
$K$ with the Newton polygon $\Delta'$. The polynomial
$f^{\Delta'}$ is a polynomial in one variable after a suitable
change of variables since $\Delta'$ is an edge. Therefore the
curve $W(\{f^{\Delta'}=0\})$ is a union of no more $n$ holomorphic
cylinders, where $n$ is the integer length of $[j,k]$, i.e. the
maximal divisor of the integer vector $k-j\in\Z^2$. Note that some
of those holomorphic cylinders may coincide if the corresponding
one-variable polynomial has multiple roots. } The curve
$W(\{f^{\Delta'}=0\})$ coincides with the curve
$W(\{f_{\min}^{\Delta'}=0\})$, where the polynomial
$f_{\min}^{\Delta'}$ is obtained by replacing each coefficient
series $a_j=a_j(t)$ of $f^{\Delta'}$ with its lower $t$-power
monomial $\alpha_jt^{\val(a_j)}$. The polynomial
$f_{\min}^{\Delta'}$ is well-defined up to a multiplication by a
complex number $c_{\Delta_j}$.
To finish the proof we need to make the equations for different
$\Delta'\in\subdiv_C(\Delta)$ agree on common monomials.
For that we order the polygons $\Delta_j\in\subdiv_C(\Delta)$ so
that $\Delta_{j+1}\cap\bigcup\limits_{l=1}^j\Delta_l$ is connected
and choose $c_{\Delta_{j+1}}$ inductively so that the monomials
of the same multidegree have the same coefficient.

$(1)\implies (3).$ This implication is a version of the so-called {\em
Viro patchworking} \cite{Vi} in real algebraic geometry. If
$f:(K^*)^2\to K$ is the polynomial defining $V_K$ then we can
construct the sequence $V_{t_k}\subset\tor$ in the following way,
cf. \cite{M-pp}.

First we truncate the polynomial $f$ by replacing each
coefficient $a=\sum\limits_{j\in I_a}\xi_j t^j,\ \xi_j\in\C,$
at a monomial of $f$ with $a_{\min}=\xi_{-\val(a)}t^{-\val(a)}$.
Denote the result with $f_{\min}$ and its zero set with
$V_K^{\min}$. As in the proof of the implication $(2)\implies (1)$
we have $W(V_K^{\min})=W(V_K)$.

Let $f^{t_k}:\tor\to\C$ be the complex polynomial
obtained from $f_{\min}$ by plugging $t=t_k$ to the coefficients of $f$.
We set $V_k$ to be the image of the zero set
of $f^{t_k}$ by the self-diffeomorphism $H_{t_k}$ defined in
\eqref{Ht}.

$(3)\implies (1).$ Propositions \ref{compact} and \ref{trvnutri}
imply that
$\Log(V_\infty)$ is a tropical curve.
Let $``\sum\limits_j \alpha_jx^j"$, $j\in\Delta\cap\Z^2$,
$\alpha_j\in\Rtr$, be a tropical polynomial defining $\Log(V_\infty)$.
To find a presentation $V_\infty=W(V_K)$ we take the polynomial
with coefficients $\beta_j t^{\alpha_j}\in K$ for some $\beta_j\in S^1$.

To find $\beta_j$ we note that if $V_k$ is a $J_{t_k}$-holomorphic
curve then $H^{-1}_{t_k}(V_k)$ is (honestly) holomorphic and is
given by a complex polynomial $f_k=\sum\limits_j a^{V_k}_jz^j$.
To get rid of the ambiguity resulting from multiplication by a
constant we may assume that $a^{V_k}_{j'}=1$ for a given
$j'\in\Delta\cap\Z^2$ and for all sufficiently large $k$.

Let $\Arg:\tor\to S^1\times S^1$ be defined by
$\Arg(z_1,z_2)=(\arg(z_1),\arg(z_2))$.
We take for $\beta_j$ an accumulation point of $\Arg(a^{V_k}_j)$,
$k\to\infty$.
Note that this accumulation point is unique and thus equal
to $\lim\limits_{k\to\infty}\Arg(a^{V_k}_j)$ since $V_{\infty}$
is the limit of $V_{t_k}$.
\end{proof}

\begin{prop}\label{jinfty-arg}
Suppose that $V_\infty\subset\tor$ is a complex tropical
curve, $ C=\Log(V_\infty)\subset\R^2$ is the corresponding
``absolute value" tropical curve
and $x\in C$ is either a vertex dual to a polygon
$\Delta'\subset\Delta$ or a point on an open edge dual to an edge
$\Delta'\subset\Delta$. We have $\Arg(\Log^{-1}(x)\cap
V_\infty)=\Arg(V')$ for some holomorphic curve $V'\subset\tor$
with the Newton polygon $\Delta'$.
\end{prop}
This proposition follows from the second characterization of
$J_\infty$-holomorphic curves in Proposition \ref{jinfty}.

The third equivalent description of Proposition \ref{jinfty}
allows one to define the genus and the number of ends
for a complex tropical curve.
\begin{defn}\label{genusinfty}
A complex tropical curve $V_{\infty}\subset\tor$
is said to have genus $g$
if $V_{\infty}$ is the limit (in the sense of the Hausdorff metric in $\tor$)
of a sequence of $J_{t_k}$-holomorphic curves in $\tor$
with $t_k\to +\infty$ of genus $g$ and cannot be presented as a limit
of a sequence of $J_{t_k}$-holomorphic curves of smaller genus.

Similarly, $V_{\infty}\subset\tor$
is said to have $x$ ends at infinity
if $V_{\infty}$ is the limit
of a sequence of $J_{t_k}$-holomorphic curves in $\tor$
with $t_k\to +\infty$ with $x$ ends at infinity and cannot be presented as a limit
of a sequence of $J_{t_k}$-holomorphic curves with smaller number of ends.
\end{defn}

Since both the genus and the number of ends are upper-semicontinuous,
this definition makes sense.

\begin{prop}\label{genusupdown}
Let $V_\infty\subset\tor$ be a complex tropical curve and
$C=\Log(V_\infty)\subset\R^2$ be the corresponding tropical curve.
We have $g(C)\le g(V_\infty)$ and
$x(C)\le x(V_\infty)$.
\end{prop}
\begin{proof}
The inequality on the number of ends follows from properness
of the map $\Log$. To get the genus inequality it suffices
to exhibit a parameterization of $C$ of genus
not more than $g=g(V_\infty)$.

We use a sequence of $J_{t_k}$-holomorphic curves
$V_k\subset\tor$, $t_k\to\infty$, $k\in\N$,
of genus $g$ approximating $V_\infty$ from
Definition \ref{genusinfty} to find such a parameterization.
Let $\nu_k:\tilde{V}_k\to V_k$ be
the normalization of the $J_{t_k}$-holomorphic curve $V_k$
(induced by the normalization of the (honestly) holomorphic
curve $H^{-1}_{t_k}(V_k)$).

Let $\{U_\alpha\}$, $U_\alpha\subset\R^2$, be a collection of
small disks in $\R^2$ centered at the vertices of $C\subset\R^2$.
Consider a component of the inverse image
$(\Log\circ\nu_k)^{-1}(U_\alpha)$.
For a sufficiently large $k$ such a component
has at least two ends (this follows from the maximum principle since
the image $\Log(\Log^{-1}(U_\alpha)\cap V_k)$ approximates $C\cap U_\alpha$)
while each end corresponds to an end of $C\cap U_\alpha$.
Thus each such component gives a subgraph
of $C\cap U_\alpha$ that is a locally tropical 1-cycle (i.e.
can be presented as an intersection of a tropical 1-cycle in $\R^2$
with $U_\alpha$). Together these components give (locally)
a parameterization of $C$ that can be extended to the edges of $C$
to get a parameterized tropical curve $h:\Gamma\to C\subset\R^2$.
Clearly, the number of components $\kappa$ of $\Gamma$ coincides
with the number of components of $\tilde{V}_k$. The same holds for
the number of ends, $x(\Gamma)=x(\tilde{V}_k)=x(V_k)$.
Denote this number with $x$.

Note that $$\chi(\tilde{V}_k)\le -v(\Gamma),$$
where $v(\Gamma)$ is the number of vertices
of $\Gamma$ since each component of $(\Log\circ\nu_k)^{-1}(U_\alpha)$
that makes up a vertex of $\Gamma$ has at least two ends.
Computing the Euler characteristic of $\tilde{V}_k$
we get $\chi(\tilde{V}_k)=2\kappa-2g-x$ and therefore
\begin{equation}\label{chicomp}
2\kappa-2g-x\le -v(\Gamma).
\end{equation}

Also we have
$2e(\Gamma)\ge 3v(\Gamma)+x$, where $e(\Gamma)$
is the number (of both bounded and unbounded) edges
of $\Gamma$, since all vertices of $\Gamma$ are at least 3-valent.
On the other hand $\chi(\Gamma)=x+v(\Gamma)-e(\Gamma)$ and
thus $\chi(\Gamma)\le\frac{x-v(\Gamma)}{2}$.
Therefore $-v(\Gamma)\le 2\chi(\Gamma)-x$.

Recall that $\chi(\Gamma)=\kappa-g(\Gamma)$ by definition
of the genus of the graph $\Gamma$.
Thus $$-v(\Gamma)\le 2\kappa-g(\Gamma)-x.$$
Combining this with \eqref{chicomp} we get the genus inequality.
\end{proof}

\begin{rmk}
It can happen that $g(C)< g(V_\infty)$.
For example take
$V^a_\infty$ to be the limiting curve for the family
$1+z_1+z_2+at^{-1}z_1z_2$, where $a\in\C$, $|a|=1$, and $t>0$,
$t\to\infty$. The image $C=\Log(V^a_\infty)$ does not depend
on the choice of $a$.
We have $g(C)=-1$.
The curve $V^a_\infty$ is a union of the lines
given by equations $z_1=-1$ and $z_2=-1$ if
$a=1$.
However for all other values of $a$ the curve
$V^a_\infty$ is irreducible and $g(V^a_\infty)=0$.
\end{rmk}

\subsection{Simple complex tropical curves
and their parameterizations}

A basic example of a complex tropical curve is a
{\em complex tropical line} $\Lambda\subset\tor$
that is a complex tropical curve
of degree $\Delta_1=\operatorname{Convex Hull}((0,0),(1,0),(0,1))$.
It is easy to see
that
any two complex tropical lines differ by
a (multiplicative) translation in $\tor$.

One can generalize this example.
Let $\Delta$ be a triangle and $C\subset\R^2$ be a tropical
curve of degree $\Delta$ with no bounded edges.
Such curve has genus 0 and 3 unbounded edges.
(Note that if $\#(\dd\Delta\cap\Z^2)>3$ then
some of the unbounded edges have weight greater than 1.)
\ignore{
Let $V_\infty\subset\tor$ be
a simple complex tropical curve $V_\infty$
projecting to $C$.
\begin{lem}\label{cttri}
The curve $V_\infty\subset\tor$ lifts to a complex tropical projective line
under some linear map $M_\Delta:\tor\to\tor$ that is self-covering
of degree $2\Area(\Delta)$.
\end{lem}
Clearly, we have $\deg(M_\Delta)=2\Area(\Delta)$ of such lifts.
\begin{proof}[Proof of Lemma \ref{cttri}]
}

Since $\Delta\subset\R^2$ is a lattice triangle
there exists an affine-linear surjection $\Delta_1\to\Delta$.
Let
\begin{equation}\label{LDelta}
L_\Delta:\R^2\to\R^2
\end{equation}
be the linear part of this map.
Note that $\det(L_\Delta)=2\Area(\Delta)$.
The matrix of $L_\Delta$ written multiplicatively defines
a map
\begin{equation}\label{MDelta}
M_\Delta:\tor\to\tor.
\end{equation}
Alternatively we can define $M_\Delta$ as the map covered
by $L_\Delta\otimes\C:\C^2\to\C^2$ under $\exp:\C^2\to\tor$.
We have $\deg(M_\Delta)=\det(L_\Delta)=2\Area(\Delta)$.

Note that $M_\Delta$ extends to a holomorphic map
\begin{equation}\label{barMDelta}
\bar{M}_\Delta: \cp^2\to\C T_\Delta.
\end{equation}
If $\Delta'\subset\dd\Delta$ is a side of $\Delta$ then
$(\bar{M}_\Delta)^{-1}(\Delta')$ consists of $2\Area(\Delta)/|\Delta'|$
components, where $|\Delta'|=\#(\Delta'\cap\Z^2)-1$ is the lattice length
of $\Delta'$.

\begin{prop}\label{obrazM}
The image $M_\Delta(\Lambda)$ is a complex tropical curve of degree $\Delta$.
\end{prop}
\begin{proof}
The lemma follows from the first description of complex tropical curves
in Proposition \ref{jinfty} since the map $L_\Delta$ also defines
a linear endomorphism $(K^*)^2\to (K^*)^2$.
\end{proof}

\begin{defn}
A curve $C\subset\tor$ is called a {\em holomorphic cylinder in $\tor$}
if $C=M_{\Delta}(\tilde{C})$, where $\tilde{C}=\{(z,w)\in\tor\ |\ w=c\}$,
for some lattice triangle $\Delta$ and $c\in\C$.
The map $M_{\Delta}|_{\tilde{C}}:\tilde{C}\to{C}$ is called
{\em a $d$-fold covering of
a holomorphic cylinder $C$} if the upper left element of
the matrix $L_{\Delta}$ is $d$.
A subset $A$ of a holomorphic cylinder $C$
is called a {\em holomorphic annulus in $\tor$}
if $A=\Log^{-1}(U)\cap C$ for a convex set $U\subset\R^2$.
Then $$M_{\Delta}|_{\tilde{A}}:\tilde{A}\to A$$ is called {\em a $d$-fold covering of
a holomorphic annulus $A$}, $\tilde{A}=M_{\Delta}^{-1}(A)\cap\tilde{C}$.
\end{defn}

Note that holomorphic cylinders in $\tor$ are complex tropical
curves and in the same time complex tropical subvarieties of $\tor$.

\begin{prop}
\label{cline-desc}
A complex tropical line $\Lambda\subset\tor$ is homeomorphic
to a sphere punctured in 3 points.
There exists a point $p\in\R^2$ such that $\Lambda\setminus\Log^{-1}(p)$
is a union of three holomorphic annuli while $\Lambda\cap\Log^{-1}(p)$
is homeomorphic to a union of two triangles whose vertices pairwise
identified.
\end{prop}
\begin{proof}
Since all complex tropical lines are multiplicative translates of
each other in $\tor$ it suffices to check the proposition for
the $\Lambda=W(\{(z,w)\in (K^*)^2\ |\ z+w+1=0\}$.
For this case the holomorphic annuli are subsets of the
cylinders $\{(z,w)\ |\ z=-1\}$, $\{(z,w)\ |\ w=-1\}$ and
$\{(z,w)\ |\ z=-w\}$ while $\Lambda\cap\Log^{-1}(p)=
\Arg(\{(z,w)\in\tor\ |\ z+w+1=0\}$ consists of all
points $(e^{i\alpha},e^{i\beta})\in (S^1)^2$, $\alpha,\beta\in\R$,
with $|\alpha|,|\beta|\le\pi$,
$\alpha\beta\ge 0$ and $|\alpha|+|\beta|\ge\pi$.
Clearly, the union of this figure with 3 holomorphic annuli along
$\{(e^{i\alpha},e^{i\beta})\in (S^1)^2\ |\ \alpha=\pi\}$,
$\{(e^{i\alpha},e^{i\beta})\in (S^1)^2\ |\ \beta=\pi\}$
and $\{(e^{i\alpha},e^{i\beta})\in (S^1)^2\ |\ \alpha+\beta=\pi\}$
is a sphere punctured 3 times.
\end{proof}

\begin{defn}
\label{psctc}
A proper map $h:\tilde{V}_\infty\to\tor$ is called a {\em
parameterized simple complex tropical curve} if
the following conditions hold.
\begin{itemize}
\item $\Log(h(\tilde{V}_\infty))\subset\R^2$
is a simple tropical curve.
\item $\tilde{V}_\infty$ is homeomorphic to
a (smooth) orientable surface.
\item If $p\in C=\Log(h(\tilde{V}_\infty))$
is a point different from a double point or a 3-valent
vertex of the simple
tropical curve $C\subset\R^2$ then there exists
a neighborhood $U\ni p$ in $\R^2$ such that
$\kappa=(\Log\circ h)^{-1}(U)$ is homeomorphic to an annulus
and $h|_{\kappa}$ is a $d$-fold covering of a holomorphic
annulus in $\tor$, where $d$ is the weight of the edge of $C$
containing $p$.
\item If $p\in C=\Log(h(\tilde{V}_\infty))$
is a double point of $C$ then there exists
a neighborhood $U\ni p$ in $\R^2$ such that
$(\Log\circ h)^{-1}(U)$ is homeomorphic to a disjoint
union of two annuli $\kappa_1$ and $\kappa_2$ while
$h|_{\kappa_j}$ is a $d_j$-fold
coverings a holomorphic annulus in $\tor$, $j=1,2$,
where $d_1$ and $d_2$ are the weights of the two edges
of $C$ containing $p$.
\item If $p\in C=\Log(h(\tilde{V}_\infty))$
is a 3-valent vertex of $C$ dual to a triangle $\Delta'\in\subdiv_C$
then there exists a neighborhood $U\ni p$ in $\R^2$ such that
$\kappa=(\Log\circ h)^{-1}(U)$ is connected and
the map
$$h|_{\kappa}:\kappa\to\tor$$
coincides with the map
$$M_{\Delta'}|_{(\Log\circ M_{\Delta'})^{-1}(U)}:
\Lambda\cap(\Log\circ M_{\Delta'})^{-1}(U)\to\tor$$
for some complex tropical line $\Lambda\subset\tor$.
\end{itemize}
As usual we consider simple parameterized tropical curves
up to the reparameterization equivalence.
Two simple parameterized tropical curves $h:\tilde{V}_\infty\to\tor$
and $h':\tilde{V}'_\infty\to\tor$ are equivalent if there
exists a homeomorphism $\Phi:\tilde{V}_\infty\to\tilde{V}'_\infty$
such that $h=h'\circ\Phi$.
\end{defn}

\begin{rmk}
Note that if $a\in\R^2$ is a 3-valent
vertex of the simple tropical curve $\Log(h(\tilde{V}_\infty))$
then $(\Log\circ h)^{-1}(a)$ is topologically
a union of two triangles whose vertices are pairwise identified.
Indeed, by Definition \ref{psctc} the inverse image $(\Log\circ h)^{-1}(a)$
is homeomorphic to the closure of the image of the argument map
of a complex tropical line in $\tor$.
\end{rmk}

By their definition the simple complex tropical curves are those
maps which locally coincide with the maps of complex tropical lines
$\Lambda\subset\tor$ by $M_\Delta$.
Enumerative geometry of lines is straightforward: there
is a single line via a pair of generic points.
Interestingly, this allows one to
make enumerative geometry of simple complex tropical curves
straightforward as well.

\begin{defn}
The genus of a simple parameterized complex tropical curve
$h:\tilde{V}_\infty\to\tor$ is the genus of the surface $\tilde{V}$.
The number of ends of $h$ is the number of ends of
the surface $\tilde{V}$.
\end{defn}
In Section \ref{prma} we shall see that this definition
agrees with Definition \ref{genusinfty}.
So far we note that these definitions agree
with the definition of genus and number of ends for
the tropical 1-cycle $\Log(h(\tilde{V}_\infty))\subset\R^2$.

\begin{prop}\label{simplegen}
The genus and the number of ends of a simple parameterized
complex tropical curve $h:\tilde{V}_\infty\to\tor$
coincide with the genus and the number of ends of
the tropical 1-cycle $C=\Log(h(\tilde{V}_\infty))$.
\end{prop}
\begin{proof}
Since the map $\Log\circ h$ is proper the number of ends of $\tilde{V}_\infty$
is not less than the number of ends of $C$.
On the other hand
each end of $C$ is a ray of weight $w$ going to infinity.
Over this ray the surface $\tilde{V}$ has a holomorphic annulus
wrapped $w$ times and, therefore, it has exactly one end.

Let $\Gamma\to C\subset\R^2$ be a simple parameterization of $C$.
To prove the equality of genera now it suffices
to show the following relation for the Euler characteristic
$$\chi(\overline{\tilde{V}_\infty})=2\chi(\overline{\Gamma}),$$
where $\overline{\tilde{V}_\infty}$ and $\overline{\Gamma}$
are results of 1-point compactifications of {\em each end} of
$\tilde{V}_\infty$ and $\Gamma$ respectively. (Note that
$\overline{\tilde{V}_\infty}$ is homeomorphic to the closure
of $h(\tilde{V}_\infty)\subset\tor$ in $\C T_{\Delta}\supset\tor$.)

We have $\chi(\overline{\Gamma})=x+v-e$, where $x$ is the number of ends
of the simple tropical curve $C$, $v$ is the number of its 3-valent
vertices and $e$ is the number of the edges of $\Gamma$.
From combinatorics of $\overline{\Gamma}$ we have
$2e=3v+x$ and therefore $$\chi(\overline{\Gamma})=\frac{x-v}{2}.$$
On the other hand the Euler characteristic of $(h\circ\Log)^{-1}(p)$
is $-1$ if $p\in C$ is a 3-valent vertex and zero otherwise.
Therefore $$\chi(\overline{\tilde{V}_\infty})=x-v.$$\end{proof}

\IGNORE{
\begin{prop}
The {\em genus of a simple parameterized tropical curve
$h:\tilde{V}_\infty\to\tor$} equals the
genus of $\tilde{V}_\infty$, the number of ends equals the number
of ends of $\tilde{V}$.
\end{prop}
\begin{proof}
Combining Propositions \ref{genusupdown} and \ref{simplegen}
we get one inequality for each of the two pairs of corresponding numbers.
W

\end{proof}
}

\subsection{Passing through a configuration
of points in $\tor$}
We start from an elementary enumerative observation.
Suppose that we have two points $q_1,q_2\in\tor$ such
that $\Log(q_1),\Log(q_2)\in\R^2$ are in tropical general position.
\begin{lem}\label{1cline}
There exists a unique complex tropical line $\Lambda\subset\tor$
such that $\Lambda\ni q_1,q_2$.
\end{lem}
\begin{proof}
We have a unique tropical line $C\subset\R^2$
such that $\Log(q_1),\Log(q_2)\in C$.
Acting on $\tor$ by an element of $SL_2(\Z)$ if needed
we may assume that $q_1$ sits on the horizontal edge of $C$
while $q_2$ sits on the vertical edge of $C$.
Recall that all complex tropical lines differ by a multiplicative
translation in $\tor$.
We have $(S^1)^2$ distinct complex tropical lines projecting
to $C$, $S^1\times\{\alpha_0\}$ of them pass through $q_2$
while $\{\beta_0\}\times S^1$ pass through $q_1$ for some $\alpha_0,\beta_0\in S^1$
(cf. Proposition \ref{cline-desc}).
Thus we have a unique tropical line through $q_1$ and $q_2$.
\end{proof}

Suppose that $C\subset\R^2$ is an image of a tropical line in $\R^2$
under a map $L_\Delta:\R^2\to\R^2$ (see \eqref{LDelta})
corresponding to a lattice triangle $\Delta\subset\R^2$.
Note that we have a natural way to equip the edges of $C$
with weights so that $C\subset\R^2$ is a tropical curve
with the help of Proposition \ref{obrazM}.
Indeed, a tropical line is the image of a complex tropical
line $\Lambda\subset\tor$ under the map $\Log$.
The image $M_{\Delta}(\Lambda)$ is a tropical curve and
thus $C=\Log(M_{\Delta})\subset\R^2$ is a tropical curve
by Proposition \ref{ctrop-trop}.

It is easy to see that the weights of the 3 edges of $C$
are the lattice lengths of the 3 sides of the triangle $\Delta$
(recall that a lattice length of a segment $[a,b]\subset\Z^2$,
$a,b\in\Z^2$ equals to the number of $\Z^2$-points on $[a,b]$
minus 1).
In the same time these weights are degrees of the covering
$\tilde{A}\to A$ of the holomorphic annuli in $\tor$
obtained as the corresponding restrictions of the map
$M_{\Delta}|_{\Lambda}:\Lambda\to\tor$. Let $w_1$ and $w_2$
be the weights of the edges of $C$ containing $q_1$ and $q_2$
respectively.

\begin{prop}\label{deltacov}
There exist $\frac{2\Area(\Delta)}{w_1w_2}$ distinct
simple complex tropical curves $N\subset\tor$ of degree $\Delta$
and genus $0$ with 3 ends at infinity such that
$N\ni q_1,q_2$.
\end{prop}
\begin{proof}
Recall that $M_{\Delta}:\tor\to\tor$ is a $2\Area(\Delta)$-fold covering.
We have $(2\Area(\Delta))^2$
distinct pairs of a point from $M_{\Delta}^{-1}(q_1)$
and a point from $M_{\Delta}^{-1}(q_2)$.
By Lemma \ref{1cline} this makes $(2\Area(\Delta))^2$
complex tropical lines in $\tor$ connecting such pairs.
However, not all of these lines are distinct.

Indeed, each complex tropical line passing through such a pair
contains $w_1$ points from $M_{\Delta}^{-1}(q_1)$ and
$w_2$ points from $M_{\Delta}^{-1}(q_2)$ so we only have
$\frac{(2\Area(\Delta))^2}{w_1w_2}$.
In addition, the $2\Area(\Delta)$-fold covering $M_{\Delta}$
identifies $2\Area(\Delta)$-tuples of such lines
(those that differ by the deck transformations).
Thus we have $\frac{2\Area(\Delta)}{w_1w_2}$
distinct simple parameterized complex tropical curves
passing through $q_1$ and $q_2$ of the form
$M_{\Delta}|_{\Lambda}:\Lambda\to\tor$ for some complex
tropical lines $\Lambda\subset\tor$.

To finish the proof we need to show that any
simple parameterized complex tropical curve
$h:\tilde{V}_\infty\to\tor$ of genus 0 with 3 ends and
passing via $q_1$ and $q_2$
is equivalent to $M_{\Delta}|_{\Lambda}:\Lambda\to\tor$
for a complex tropical line $\Lambda$.
Clearly, we have $\Log(h(\tilde{V}_\infty))=C$ as $C$
is the only tropical curve of genus 0 with 3 ends passing
via $\Log(q_1)$ and $\Log(q_2)$.

Let us note that $h:\tilde{V}_\infty\to\tor$
lifts under a covering $M_{\Delta}:\tor\to\tor$.
To establish this it is convenient to consider
the compactification map \ref{barMDelta}.
The fundamental group of the three times punctured
sphere $\tilde{V}_\infty$ is generated by two elements,
$(\Log\circ h)^{-1}(q_1)$ and $(\Log\circ h)^{-1}(q_2)$
as the three punctures correspond to the three sides
of the lattice triangle $\Delta$.
Both of these elements are in the kernel of
the induced homomorphism $M_{\Delta*}:\pi_1(\tor)\to\pi_1(\tor)$
as they wrap the holomorphic annuli of $h(\tilde{V}_\infty)$
$w_1$ and $w_2$ times respectively. The lift of $\tilde{V}_\infty$
is a complex tropical curve of degree $\Delta_1$ and, therefore
is a complex tropical line.
\end{proof}

In general case consider
a configuration $\qq\subset\tor$ of $n$ points such
that $\Log(\qq)\subset\R^2$ is in tropically general position.
Let $C\supset\Log(\qq)$ be a tropical curve in $\R^2$ of
genus $g$ and degree $\Delta$. Let $\mult$ be the multiplicity
of $C$ (see Definition \ref{multdim2}) and let $x\le s$ be the number
of ends of $C$. Suppose that $n=x+g-1$.

\IGNORE{
\begin{lem}
Any complex tropical curve $V_\infty$ of genus $g$ with $x$ ends that pass
via $\qq$ is a simple complex tropical curve.
\end{lem}
\begin{proof}

\end{proof}
}

Let $V_\infty$ be a simple complex tropical curve such that
$$\qq\subset V_\infty\gap \text{and}\gap \Log(V_\infty)=C\subset\R^2.$$
We define {\em the edge multiplicity $\mu_{\operatorname{edge}}(C,\ppp)$
of a simple tropical curve} $C\supset\ppp=\Log(\qq)$
to be the product of the weights
of all the edges of the parameterizing graph $\Gamma\to C\subset\R^2$
that are disjoint from $\ppp$ times
the product of the squares of the weights of all the edges of $\Gamma$
that are not disjoint from $\ppp$.
\IGNORE{
We define {\em the $\qq$-edge multiplicity $\mu(V_\infty,\qq)$
of a simple complex tropical curve} $V_\infty$
to be the product of the weights (i.e. the degrees of
the covering $\tilde{A}\to A$) of the elementary annuli
that contain points from $\qq$.
This multiplicity only depends on the projection $\Log(V_\infty)=C$
and $\Log(\qq)$
since it is equal to the product of the weights of all edges of $C$
containing the points from $\Log(\qq)$. Accordingly, we also
denote this multiplicity with $\mu(C,\Log(\qq))=\mu(V_\infty,\Log(\qq))$.
}
Note that the unbounded edges of $\Gamma$ are all of weight 1
by Proposition \ref{tropfinite} so they do not contribute
to the edge multiplicity.
\begin{prop}\label{ctropmain}
There are $\mult(C)/\mu_{\operatorname{edge}}(C,\ppp)$
simple complex tropical curves in $\tor$
of genus $g$ and degree $\Delta$ such that they project
to $C$ and pass via $\qq$.
\end{prop}

\begin{proof}
\IGNORE{
We start with a special case.
Suppose (as in Lemma \ref{cttri})
that $g=0$, $\Delta$ is a triangle and $C$ has 3 ends.
Then $\qq$ consists of two points.

By Lemma \ref{cttri} any simple complex tropical curve
over $C$ is an image of a complex tropical projective line
under $M_\Delta$. The set $(M_\Delta)^{-1}(\qq)$ consists
of $2\deg(M_\Delta)=4\Area(\Delta)$ points.
Each pair projecting to different points under $M_\Delta$
can be connected with a complex tropical projective line.
We have a total of $(\deg(M_\Delta))^2$ of such lines and
each of them projects to a simple complex tropical curve
of degree $\Delta$ under $M_\Delta$. Conversely, each
such simple complex tropical curve lifts to $d$ distinct lines.
Therefore, we have $\deg(M_\Delta)=2\Area(\Delta)$ distinct
simple complex tropical curves projecting to $C$ and passing via $\qq$.

Now we are ready to treat general case.
}
Recall that $\ppp=\Log(\qq)$ is a set of $x+g-1$ points in general position
As in Lemma \ref{descompl} let $K$ be a component of
the complement in $\Gamma$ of the inverse image of $\Log(\qq)$
under the parameterization $\Gamma\to C$.
By Lemma \ref{descompl} the component $K$
is a tree which contains one end at infinity.
Let $A,B\in\qq$ be two points that are the endpoints
of the edges of weights $w_A$ and $w_B$ adjacent to the same
3-valent vertex $\nu$ in $K$ (as in Figure \ref{komp-izm}).
Let $\Delta'\subset\Delta$ be the subpolygon dual to this 3-valent vertex of $C$.
By Proposition \ref{deltacov}
we have $\frac{2\Area(\Delta')}{w_Aw_B}$ simple parameterized
complex tropical curves that are locally distinct in $\Log^{-1}(U)$,
where $U\subset\R^2$ is a small neighborhood of the vertex $\nu$.
We proceed inductively for each component $K$ and then take
the product over all such components.
\end{proof}

\subsection{Polynomials to define complex tropical curves}
Let $V_K\subset (K^*)^2$ be a curve given by a non-Archimedean
polynomial $$f_K(z)=\sum\limits_j \alpha_jz^j,$$
$z\in (K^*)^2$, $j\in\Delta$,
$\alpha_j\in K^*$.
Let $V_\infty=W(V_K)\subset\tor$ be the corresponding complex tropical curve.
By Kapranov's theorem \cite{Ka} (cf. Proposition \ref{ctrop-trop})
$f(y)=``\sum\limits_j \val(\alpha_j)y^j "$ is a tropical polynomial
defining the tropical 1-cycle $C=\Log(V_\infty)=\Val(V_K)$.

\begin{prop}\label{ctrop-poly}
If $j$ and $j'$ are vertices of $\subdiv_f$ then
the ratio $\frac{w(\alpha_j)}{w(\alpha_{j'})}$
(where $w$ is defined by \eqref{mapw}) depends only
on $V_\infty$ and does not depend on the choice of $V_K$
as long as $W(V_K)=V_\infty$.
\end{prop}
\begin{proof}
Since $\Delta$ is connected we need to prove
the proposition only in the case when $j$ and $j'$ are
connected with an edge $\Delta'\in\subdiv_f$.
The complex tropical curve $V_\infty$ is $\Delta'$-compatible,
therefore the proposition follows from the special case when
$\Delta=\Delta'$ is 1-dimensional.

After an automorphism in $\tor$ we may assume that
$\Delta=\Delta'=[0,l]\times\{0\}\subset\R^2$,
$l\in\N$. Then $V_\infty$ is a collection of the holomorphic cylinders
$\{(\zeta,\eta)\in\tor\ |\ \zeta=z_k\}$, $z_k\in\C$, of weight $d_k$.
(Note that we have $|z_k|$ the same for all $k$ since
all these cylinders must project to the same edge of $C$.
Every cylinder is given by the equation $\zeta^{d_k}=z_k^{d_k}$.
This equation incorporates the weight of the cylinder
and is well-defined up to a multiplication by a constant
or a monomial.
Proceeding inductively we
get $\frac{w(\alpha_j)}{w(\alpha_{j'})}=\prod z_k^{d_k}$.
\end{proof}

Thus if $j\mapsto -\val(\alpha_j)$ is strictly convex then
to recover a complex tropical curve $w(V_K)$
it suffices to know only $a_j=w(\alpha_j)$.
\begin{prop}\label{ctrop-thm}
Let $f(z)=\sum\limits_{j\in\Delta} a_j z^j,$
be a formal sum of monomials, $a_j\in\C^*$, where $j$ runs over
all lattice points of $\Delta$.
Suppose that $j\mapsto -\log|a_j|$
is strictly convex. Then $f$ defines a complex tropical curve
$V_\infty\subset\tor$ of degree $\Delta$.
We have the following properties.
\begin{itemize}
\item
If two polynomials define the same complex tropical curve
then they differ by a multiplication by a constant.
\item
Complex tropical curves defined in this way form an open and dense
set in the space of all complex tropical curves of degree $\Delta$.
\end{itemize}
\end{prop}
\begin{proof}
Form a polynomial
$$f^K(z)=\sum\limits_{j\in\Delta\cap\Z^2} a_{j}t^{\log|a_{j}|}z^j,$$
and treat it as a polynomial over $K$, i.e. $z\in (K^*)^2$.
Define $V_\infty=W(\{z\in (K^*)^2\ |\
f^K(z)=0\}.$
To see that different polynomials define different curves it suffices
to note that for the corresponding tropical polynomial
$$f_{\operatorname{trop}}(y)=``\sum\limits_{j\in\Delta} \log|a_j| y^j "$$
the subdivision $\subdiv_{f_{\operatorname{trop}}}$ contains all lattice
points of $\Delta$ as its vertices. In particular, no lattice
point of $\Delta$ is contained in the interior of an edge of
$\subdiv_{f_{\operatorname{trop}}}$. For different
polynomials we have different $\Delta'$-compatible elements
for some edge $\Delta'\in\subdiv_{f_{\operatorname{trop}}}$.

On the other hand, any complex tropical curve can be obtained
from a non-Archimedean curve $V_K\subset(K^*)^2$ given by
a polynomial $f_K(z)=\sum\limits_j \alpha_jz^j$ such that
the function $j\mapsto -\val(\alpha_j)$ is convex.
The polynomials with strictly convex functions $j\mapsto -\val(\alpha_j)$
are open and dense among all such polynomials.
\end{proof}

\begin{rmk}
If $V_K\subset(K^*)^2$ is given by a polynomial
$f_K(z)=\sum\limits_j \alpha_jz^j$ such that
the function $j\mapsto -\val(\alpha_j)$ is not strictly convex
(even if it is convex non-strictly)
then the collection of complex numbers $a_j=w(\alpha_j)$ does
not necessarily determine the complex tropical curve $w(V_K)$.

E.g. the polynomials $f_K(\zeta,\eta)=\zeta^2+3\zeta+2$ and
$f'_K(\zeta,\eta)=\zeta^2+\zeta+1$
(treated as polynomials over $K$) produce the same collection
of three complex numbers, namely $a_{(0,0)}=a_{(1,0)}=a_{(2,0)}=1$.
However we have $W(V_K)=\{(\zeta,\eta)\in\tor\ |\ \zeta=-1\}$ and
$W(V'_K)=\{(\zeta,\eta)\in\tor\ |\ \zeta=-\frac12\pm\frac{\sqrt{3}}{2}\}.$
\end{rmk}

\ignore{
The following proposition follows from the second
description from Proposition \ref{jinfty}.
\begin{prop}\label{jinfty-param}
Any complex tropical curve $V_\infty$ with the Newton
polygon $\Delta$ can be obtained in this way from a polynomial
$$f(z,w)=\sum\limits_{(j,k)\in\Delta\cap\Z^2} a_{jk}z^jw^k$$ such
that $(j,k)\mapsto\log|a_{j,k}|$ defines a convex function
$\alpha:\Delta\cap\Z^2\to\R$. Furthermore, if $V$ is obtained from
a polynomial $f$ with a strictly convex $\alpha$ then the choice
of such $f$ is unique.
\end{prop}
}

\section{Statement of the main theorems}
\label{sectionmain}
\subsection{Enumeration of complex curves}
Recall that the numbers $\ntrop(g,\Delta)$ and $\ntropr(g,\Delta)$
were introduced in Definition \ref{tropnumbers} and {\em a priori}
they depend on a choice of a configuration $\ppp\subset\R^2$ of
$s+g-1$ points in general position.
\begin{thm}\label{main}
For any generic choice $\ppp$ we have
$\ntrop(g,\Delta)=\ncl(g,\Delta)$ and
$\ntropr(g,\Delta)=\nclr(g,\Delta)$.

Furthermore, there exists a configuration $\qq\subset\tor$
of $s+g-1$ points in general position such that
for every tropical curve $C$ of genus $g$ and degree $\Delta$ passing
through $\ppp$ we have $\mult(C)$ (see Definition \ref{multdim2})
distinct complex curves of genus $g$
and degree $\Delta$ passing through $\qq$. These curves
are distinct for distinct $C$ and are irreducible if
$C$ is irreducible.
\end{thm}

The following is an efficient way to compute $\ntropr(g,\Delta)$
(and therefore also $\nclr(g,\Delta)$).
Let us choose a configuration $\ppp$ to be contained
in an affine line ${L}\subset\R^2$.
Furthermore, we make sure that the order of $\ppp$ coincides
with the order on ${L}$ and that
the distance between $p_{j+1}$ and $p_j$ is much greater
than the distance between $p_j$ and $p_{j-1}$ inductively for each $j$.
These conditions ensure that $\ppp$ is in tropically general
position as long as the slope of ${L}$ is irrational.
Furthermore, these conditions specify the combinatorial types
of the tropical curves of
genus $g$ and degree $\Delta$ passing via $\ppp$
(see Definition \ref{defcombtype}).
We shall see that for this choice of $\ppp$
the forests corresponding to such combinatorial types are
paths connecting a pair of vertices of $\Delta$.
This observation can be used to compute $\nclr(g,\Delta)$
once an irrational slope of ${L}$ is chosen.

\subsection{Counting of complex curves by lattice paths}
\begin{defn}
A path $\gamma:[0,n]\to\R^2$, $n\in\N$, is called a {\em lattice path}
if $\gamma|_{[j-1,j]}$, $j=1,\dots,n$ is an affine-linear map
and $\gamma(j)\in\Z^2$, $j\in 0,\dots,n$.
\end{defn}
Clearly, a lattice path is determined by its values at
the integer points.
%
Let us choose an auxiliary linear map $$\lambda:\R^2\to\R$$ that is irrational,
i.e. such that $\lambda|_{\Z^2}$ is injective.
Note that an affine line ${L}\subset\R^2$ with
an irrational slope determines a choice of $\lambda$ as
we can take $\lambda$ to be the orthogonal projection onto ${L}\subset\R^2$.
Let $p,q\in\Delta$ be the vertices where $\alpha|_\Delta$ reaches
its minimum and maximum respectively.
A lattice path is called {\em $\lambda$-increasing}
if $\lambda\circ\gamma$ is increasing.
\begin{figure}[h]
\centerline{\psfig{figure=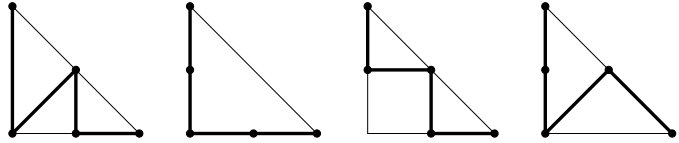,height=.71in,width=3.5in}}
\caption{\label{paths} All $\lambda$-increasing paths for a triangle with
vertices $(0,0)$, $(2,0)$ and $(0,2)$ where
$\lambda(x,y)=x-\epsilon y$ for a small $\epsilon>0$ and
$\delta=1$.}
\end{figure}


The points $p$ and $q$ divide the boundary $\dd\Delta$ into
two increasing lattice paths
$$\alpha^+:[0,n_+]\to\dd\Delta\ \  \text{and}\ \
\alpha^-:[0,n_-]\to\dd\Delta.$$
We have $\alpha_+(0)=\alpha_-(0)=p$,
$\alpha_+(n_+)=\alpha_-(n_-)=q$, $n_++n_-=m-l+3$.
To fix a convention we assume that $\alpha_+$ goes clockwise
around $\dd\Delta$ while $\alpha_-$ goes counterclockwise.

Let
$\gamma:[0,n]\to\Delta\subset\R^2$
be an increasing lattice path such that $\gamma(0)=p$
and $\gamma(n)=q$.
The path $\gamma$ divides $\Delta$ into two closed regions:
$\Delta_+$ enclosed by $\gamma$ and $\alpha_+$
and $\Delta_-$ enclosed by $\gamma$ and $\alpha_-$.
Note that the interiors of $\Delta_+$ and $\Delta_-$
do not have to be connected.

We define the positive (resp. negative)
multiplicity $\mu_\pm(\gamma)$
of the path $\gamma$ inductively.
We set $\mu_\pm(\alpha_\pm)=1$.
If $\gamma\neq\alpha_\pm$
then we take $1\le k\le n-1$ to be the smallest number
such that $\gamma(k)$ is a vertex of $\Delta_\pm$
with the angle less than $\pi$ (so that $\Delta_\pm$
is locally convex at $\gamma(k)$).

If such $k$ does not exist we set $\mu_\pm(\gamma)=0$.
If $k$ exist we consider two other increasing lattice paths connecting
$p$ and $q$,
$\gamma':[0,n-1]\to\Delta$ and $\gamma'':[0,n]\to\R^2$.
We define $\gamma'$ by $\gamma'(j)=\gamma(j)$ if $j<k$ and
$\gamma'(j)=\gamma(j+1)$ if $j\ge k$.
We define $\gamma''$ by $\gamma''(j)=\gamma(j)$ if $j\neq k$ and
$\gamma''(k)=\gamma(k-1)+\gamma(k+1)-\gamma(k)\in\Z^2$.
We set
\begin{equation}\label{pathmult}
\mu_\pm(\gamma)=2\operatorname{Area}(T)\mu_\pm(\gamma')+
\mu_\pm(\gamma''),
\end{equation}
where $T$ is the triangle with the vertices $\gamma(k-1)$,
$\gamma(k)$ and $\gamma(k+1)$.
The multiplicity is always integer since the area of a lattice
triangle is half-integer.

Note that it may happen that $\gamma''(k)\notin\Delta$.
In such case we use a convention  $\mu_\pm(\gamma'')=0$.
We may assume that $\mu_\pm(\gamma')$ and $\mu_\pm(\gamma'')$
is already defined since the area of $\Delta_\pm$ is smaller
for the new paths. Note that $\mu_\pm=0$ if $n<n_\pm$
as the paths $\gamma'$ and $\gamma''$ are not longer than $\gamma$.

We define {\em the multiplicity $\mu(\gamma)$ of the path $\gamma$}
as the product $\mu_+(\gamma)\mu_-(\gamma)$.
Note that the multiplicity of a path connecting two
vertices of $\Delta$ does not depend on $\lambda$.
We only need $\lambda$ to determine whether a path
is increasing.

\begin{exa}
Consider the path $\gamma:[0,8]\to\Delta_3$ depicted
on the extreme left of Figure \ref{mult+}.
This path is increasing
with respect to $\lambda(x,y)=x-\epsilon y$, where
$\epsilon>0$ is very small.

\begin{figure}[h]
\centerline{\psfig{figure=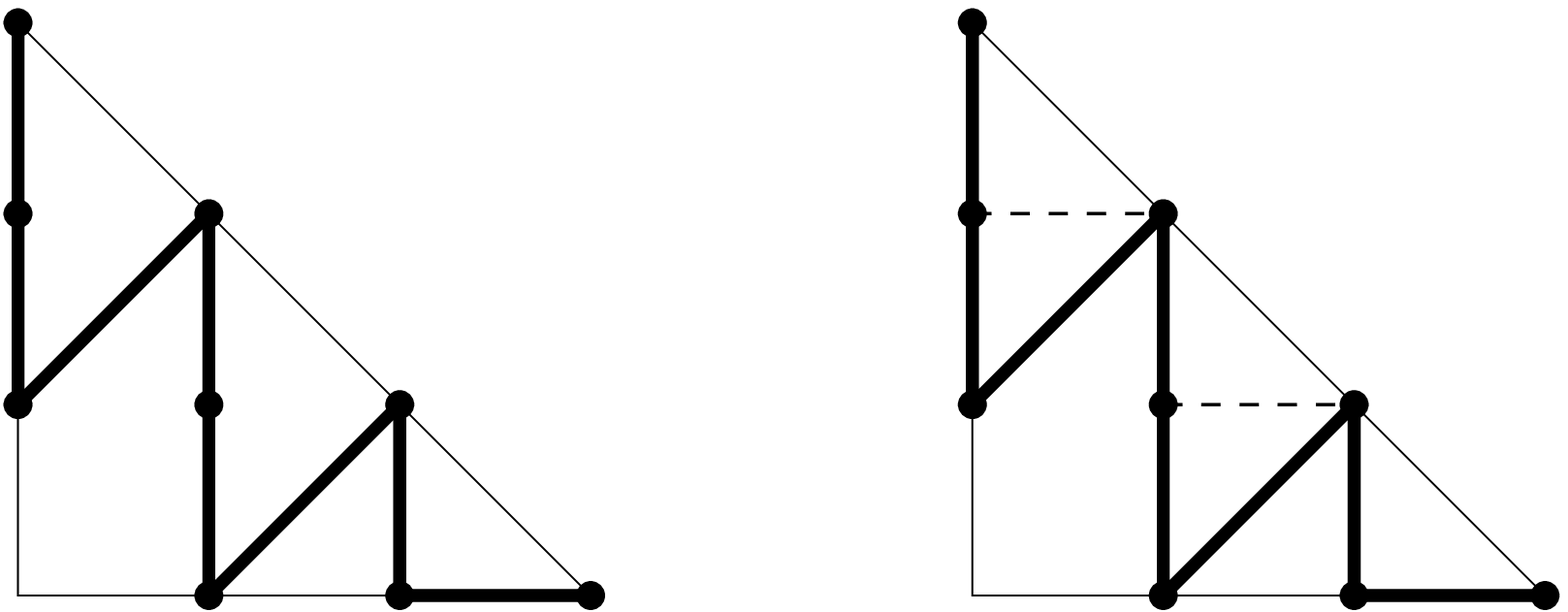,height=0.75in,width=1.95in}
\hspace{0.5in}
\psfig{figure=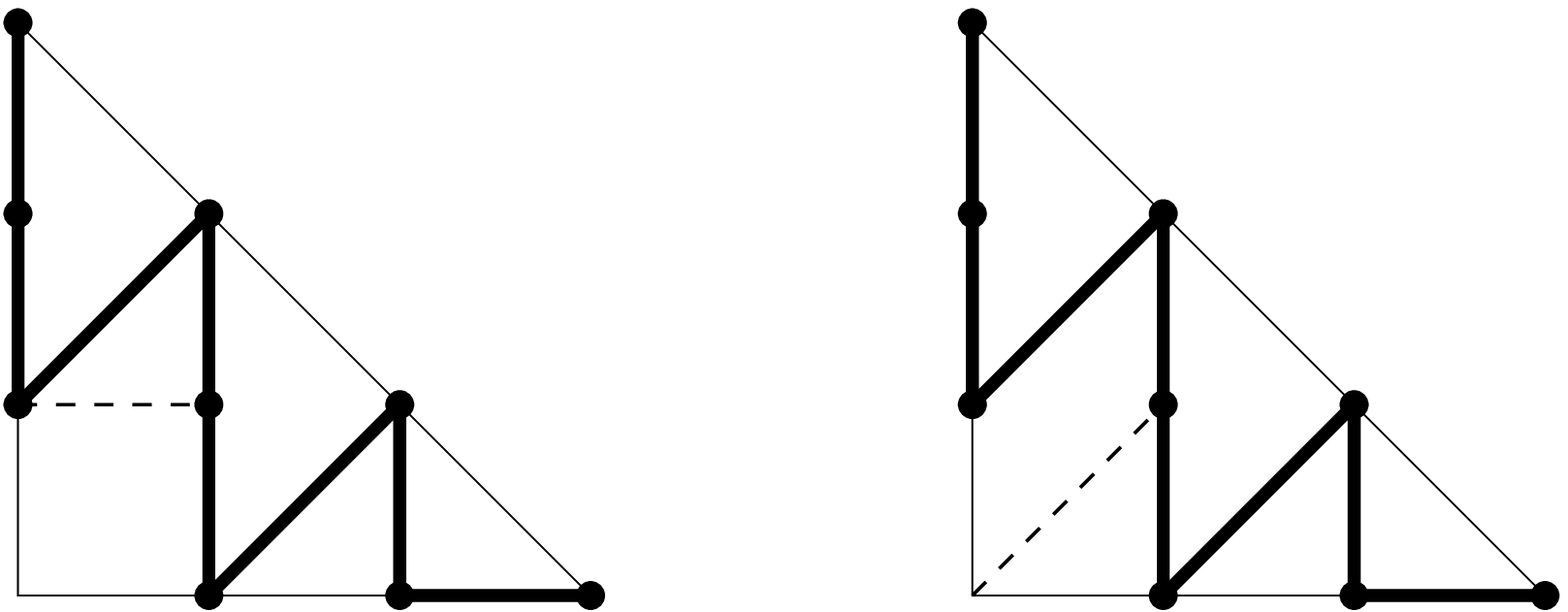,height=0.75in,width=1.95in}}
\caption{\label{mult+} A path $\gamma$ with $\mu_+(\gamma)=1$
and $\mu_-(\gamma)=2$.}
\end{figure}

Let us compute $\mu_+(\gamma)$. We have $k=2$ as $\gamma(2)=(0,1)$
is a locally convex vertex of $\Delta_+$.
We have $\gamma''(2)=(1,3)\notin\Delta_3$ and thus
$\mu_+(\gamma)=\mu_+(\gamma')$, since $\operatorname{Area}(T)=\frac12$.
Proceeding further we get
$\mu_+(\gamma)=\mu_+(\gamma')=\dots=\mu_+(\alpha_+)=1$.


Let us compute $\mu_-(\gamma)$. We have $k=3$ as $\gamma(3)=(1,2)$
is a locally convex vertex of $\Delta_-$.
We have $\gamma''(3)=(0,0)$ and $\mu_-(\gamma'')=1$.
To compute $\mu_-(\gamma')=1$ we note that $\mu_-((\gamma')')=0$
and $\mu_-((\gamma')'')=1$.
Thus the full multiplicity of $\gamma$ is 2.
\end{exa}

Recall that we fixed an (irrational) linear function
$\lambda:\R^2\to\R$ and this choice gives us a pair
of extremal vertices $p,q\in\Delta$.

\begin{thm}
\label{thm1}
The number $\ntropr(g,\Delta)$ is equal to the
number (counted with multiplicities)
of $\lambda$-increasing lattice paths
$[0,s+g-1]\to\Delta$ connecting $p$ and $q$.

Furthermore, there exists a configuration $\ppp\in\R^2$
of $s+g-1$ points in tropical general position such that
each $\lambda$-increasing lattice path encodes a number
of tropical curves of genus $g$ and degree $\Delta$
passing via $\ppp$ of total multiplicity $\mu(\gamma)$.
These curves are distinct for distinct paths.
\end{thm}

\begin{exa}\label{ecusp}
Let us compute $\nclr(0,\Delta)=5$ for the polygon $\Delta$
depicted on Figure \ref{cusp} in two different ways.
Using $\lambda(x,y)=-x+\epsilon y$ for a small $\epsilon>0$
we get the left two paths depicted on Figure \ref{cusp}.
Using $\lambda(x,y)=x+\epsilon y$ we get the three right paths.
The corresponding multiplicities are shown under the path.
All other $\lambda$-increasing paths have zero multiplicity.
\begin{figure}[h]
\centerline{\psfig{figure=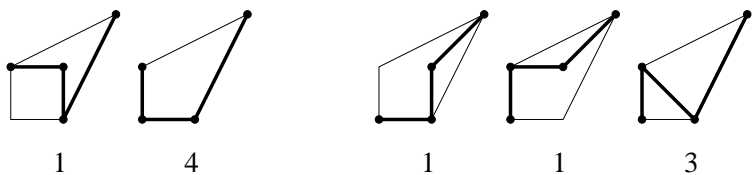,height=0.85in,width=3.9in}}
\caption{\label{cusp} Computing $\nclr(0,\Delta)=5$ in two different ways.}
\end{figure}
\end{exa}

In the next two examples we use $\lambda(x,y)=x-\epsilon y$
as the auxiliary linear function.

\begin{exa}\label{edeg3}
Figure \ref{deg3} shows a computation of
the well-known number $\nclr(0,\Delta_3)$.
This is the number of rational cubic curves through 8 generic
points in $\cp^2$.
\begin{figure}[h]
\centerline{\psfig{figure=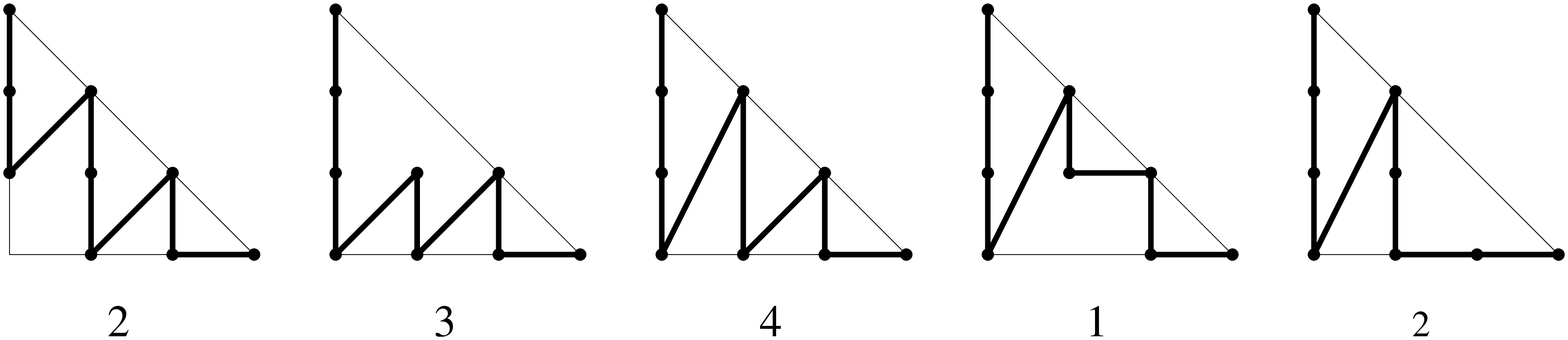,height=0.88in,width=4.1in}}
\caption{\label{deg3} Computing $\nclr(0,\Delta_3)=12$.}
\end{figure}
\end{exa}

%
\begin{exa}\label{edeg4}
Figure \ref{deg4} shows a computation
of a less well-known number $\nclr(1,\Delta_4)$.
This is the number of genus 1 quartic curves through
12 generic points in $\cp^2$.
\begin{figure}[h]
\centerline{\psfig{figure=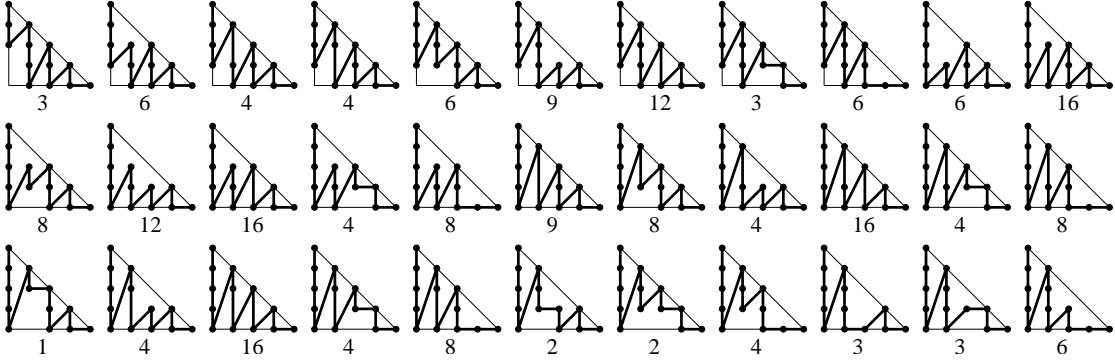,height=1.85in,width=5.8in}}
\caption{\label{deg4} Computing $\nclr(1,\Delta_4)=225$.}
\end{figure}
\end{exa}

\subsection{Enumeration of real curves}\label{enumrcurves}
Let $\qq=\{q_1,\dots,q_{s+g-1}\}\subset\rtor$ be a configuration
of points in general position. We have a total of $\nclr(g,\Delta)$
complex curves of genus $d$ and degree $\Delta$ in $\tor$ passing
through $\qq$. Some of these curves are defined over $\R$
while others come in complex conjugated pairs.
\begin{defn}\label{realenumerativeproblem}
We define the number $\ncl_{\R}(g,\Delta,\qq)$ to be the number of
irreducible real curves of genus $g$ and degree $\Delta$
passing via $\qq$.
Similarly we define the number $\nclr_{\R}(g,\Delta,\qq)$ to be
the number of all real curves of genus $g$ and degree $\Delta$
passing via $\qq$.
\end{defn}
Unlike the complex case the numbers $\ncl_{\R}(g,\Delta,\qq)$
and $\nclr_{\R}(g,\Delta,\qq)$ do depend on the choice of $\qq$.
We have
$$0\le\nclr_{\R}(g,\Delta,\qq)\le\nclr(g,\Delta)$$
while $\nclr_{\R}(g,\Delta,\qq)\equiv\nclr(g,\Delta)\pmod2$
and, similarly, $0\le\ncl_{\R}(g,\Delta,\qq)\le\ncl(g,\Delta)$
and $\ncl_{\R}(g,\Delta,\qq)\equiv\ncl(g,\Delta)\pmod2$.

Tropical geometry allows one to compute $\nclr_{\R}(g,\Delta,\qq)$
and $\nclr_{\R}(g,\Delta,\qq)$ for some configurations $\qq$.
Note that one can extract a {\em sign} sequence
$\{\sign(q_j)\}_{j=1}^{s-g+1}\subset\Z^2_2$
from $\qq=\{q_j\}_{j=1}^{s-g+1}$ by taking the coordinatewise sign.
Accordingly we can enhance the tropical configuration data by
adding a choice of signs that take values in
$\Z^2_2$.

\begin{defn}
A {\em signed} tropical configuration of points
$$\qqq=\{(r_1,\sigma_1),\dots,(r_{s+g-1},\sigma_{s+g-1})\}$$
is a collection of $s-g+1$ points $r_j$ in the tropical plane $\R^2$
together with a choice of signs $\sigma_j\in \Z^2_2$.

We denote with $|\qqq|=\{r_1,\dots,r_{s+g-1}\}\subset\R^2$ the resulting
configuration after forgetting the signs.
\end{defn}

Suppose that $C$ is a simple tropical curve and is given
by an immersion $h:\Gamma\to\R^2$ where $\Gamma$ is 3-valent.
Let $E$ be an edge of $\Gamma$ of weight $w_E$.
Suppose that $(x_E,y_E)$ is a primitive integer vector parallel to $E$.
%
We define the set $S_E$ as the quotient of $\Z^2_2$ by the equivalence
relation $(X,Y)\sim (X+w_Ex_E,Y+w_Ey_E)$, $X,Y\in\Z_2$.
Thus if the weight $w_E$ is even then $S_E=\Z^2_2$ but
if $w_E$ is odd then $S_E$ is a 2-element set.

An element $\sigma_E\in S_E$ is called {\em the phase} of $E$.
Every vertex of $\Gamma$ is adjacent to three edges $E_1,E_2,E_3$.
Their phases $\sigma_{E_1},\sigma_{E_2}, \sigma_{E_3}$ are called
{\em compatible} if there exist representatives $\rho^j_k\in\Z^2_2$,
$j,k=1,2,3$, $j\neq k$, such that
\begin{itemize}
\item $\rho^j_k\in \sigma_{E_k}$,
\item $\rho^{j}_k\neq\rho^{j'}_k$ if $j\neq j'$ and the equivalence
class $\sigma_{E_k}$ contains more than one element,
\item $\rho^k_j=\rho^j_k$.
\end{itemize}
In other words, once we associate to every phase $\sigma_{E_k}$
its two possible representatives in $\Z^2_2$
(which we choose distinct if there are two
elements in $\sigma_{E_k}$ and coinciding if $\sigma_{E_k}$
is a one-point set) the 6 resulting elements of $\Z^2_2$
should divide in three pairs of coinciding elements from
distinct edges.


\begin{defn}
A simple {\em real tropical curve} (or {\em signed tropical curve})
$(C,\{\sigma_E\})$
is a simple tropical curve $C$ parameterized by
$h:\Gamma\to\R^2$ with a choice of {\em the phase}
$\sigma_E\in S_E$ for every edge $E\subset\Gamma$ that is compatible
at every vertex of $\Gamma$.
We say that $(C,\{\sigma_E\})$ passes
through a signed configuration $\qqq$ if for every $r_j\in|\qqq|$
there exists a (closed) edge $E\subset C$ such that
$r_j\in E$ and $\sigma_j\in\sigma_E$.
\end{defn}

\begin{rmk}
Signed tropical configuration and real tropical curves
both can be viewed as 0- and 1-dimensional real tropical varieties
in $\R^2$ respectively.
\end{rmk}

Suppose that $\qqq$ is such that $|\qqq|\subset\R^2$ is
a configuration of $n$ points in general position.
Let $(C,\{\sigma_E\})$ be a real tropical curve passing through
a signed configuration $\qqq$. Suppose that $x+g-1=n$, where
$x$ is the number of ends of $C$ and $g$ is the genus of $C$.

\begin{defn}
The multiplicity $\mu(C,\{\sigma_E\},\qqq)$
is $2^{N^{\operatorname{even}}}$, where $N^{\operatorname{even}}$ is twice
the number of edges of $\Gamma$ of even weight that contain
points from $|\qqq|$ plus the number of the remaining edges
of $\Gamma$ of even weight.
\end{defn}

A real tropical curve $(C,\{\sigma_E\})$ is called irreducible
if the tropical curve $C$ is irreducible. Otherwise, it is called
reducible.

\begin{defn}
We define the number $\ntropR(g,\Delta,\qqq)$ to be the number of
real irreducible tropical curves $(C,\{\sigma_E\})$ of genus $g$ and degree $\Delta$
passing via $\qqq$ counted with multiplicities $\mu(C,\{\sigma_E\},\qqq)$.
Similarly we define the number $\ntropRr(g,\Delta,\qqq)$ to be
the number of all real tropical curves of genus $g$ and degree $\Delta$
passing via $\qqq$ (again, counted with multiplicities $\mu(C,\{\sigma_E\},\qqq)$).
\end{defn}

\begin{thm}\label{realmain}
Suppose that $\qqq$ is a signed configuration of $s+g-1$ points
in tropically general position. Then there exists a configuration
$\qq\subset\rtor$ of $s+g-1$ real points in general position
such that $\ncl_{\R}(g,\Delta,\qq)=\ntropR(g,\Delta,\qqq)$
and $\nclr_{\R}(g,\Delta,\qq)=\ntropRr(g,\Delta,\qqq)$.

Furthermore,
for every real tropical curve $(C,\{\sigma_E\})$ of genus $g$
and degree $\Delta$ passing
through $\qqq$ we have $\mu(C,\{\sigma_E\},\qqq)$
distinct real curves of genus $g$
and degree $\Delta$ passing through $\qq$. These curves
are distinct for distinct $C$ and are irreducible if
$C$ is irreducible.
\end{thm}


Our next goal is to define the {\em real multiplicity} $\mu_{\R}(C,\qqq)$
of a tropical curve $C$ passing through a signed configuration $\qqq$
in a way to include all signed tropical curves
$(C,\{\sigma_E\})$ with the same $C$ with corresponding multiplicities.
The multiplicity $\mu_{\R}(C,\qqq)$
depends only on $C$, $|\qqq|$ and the equivalence class $\sigma_E\in S_E$
of the sign $\sigma_j$ for every point $r_j\in\qqq$ contained in the edge $E$.
(Recall that since $|\qqq|\subset\R^2$ is tropically in general position,
the points $r_j$ are disjoint from the vertices of $C$.)

Let $h:\Gamma\to\R^2$ be a tropical curve of genus $g$ and degree $\Delta$ passing
through $|\qqq|$.
Recall that by Lemma \ref{descompl} each component of
$\Gamma\setminus h^{-1}(|\qqq|)$ is a tree with a single end
at infinity (see Figure \ref{typkomp}).

We define the real tropical multiplicity $\mu_\R(T,\qqq)$
of each component $T$ of $\Gamma\setminus h^{-1}(|\qqq|)$ inductively.
Let $A$ and $B$ be two 1-valent vertices of $T$
corresponding to marked points $r_a$ and $r_b$
such that the edges adjacent to $A$ and $B$ meet at
a 3-valent vertex $C$ (see Figure \ref{komp-izm}).
\begin{figure}[h]
\centerline{\psfig{figure=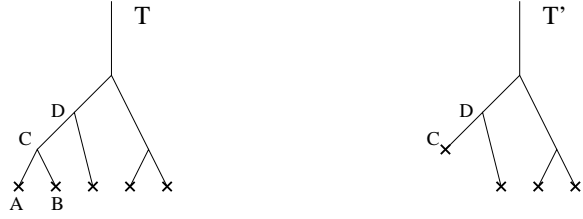,height=1.1in,width=3in}}
\caption{\label{komp-izm} Inductive reduction of the components
of $\Gamma\setminus h^{-1}(|\qqq|)$ in the definition of real
multiplicity}
\end{figure}

Form a new tree $T'$ by removing the edges $[A,C]$ and $[B,C]$
from $T$. The number of 1-valent vertices of $T'$ is less by one
($C$ becomes a new 1-valent vertex while $A$ and $B$ disappear).
By the induction assumption the real multiplicity of $T'$ is
already defined for any choice of signs. All the finite 1-valent
vertices of $T'$ except for $C$ have their signs induced
from the signs of $T$. To completely equip $T'$ with the
signs we have to define the sign $\sigma_d$ at the edge $[C,D]$.

Suppose that $(x_a,y_a)$ and $(x_b,y_b)$ are primitive integer
vectors parallel to $[A,C]$ and $[B,C]$ respectively.
Suppose that the signs of $[A,C]$ and $[B,C]$
are $\sigma_a$ and $\sigma_b$.
Suppose that $[C,D]$ is the third edge adjacent to $C$.
Let $S_a$, $S_b$, $S_d$ be the set of equivalence classes of signs corresponding
to $[A,C]$, $[B,C]$ and $[C,D]$.
Let $w_a, w_b, w_d$ be their weights.

\begin{defn}\label{realmult}
The sign at $C$ and the real multiplicity
of $T$ are defined according to the following inductive rules.
If $T$ does not have a 3-valent vertex (i.e. $T$ is homeomorphic
to $\R$) then $\mu_R(T)=1$.
\begin{itemize}
\item Suppose that $w_a\equiv w_b\equiv 1\pmod2$
and $(x_a,y_a)\equiv (x_b,y_b)\pmod2$.
In this case we have $S_a=S_b$ so the
signs $\sigma_a$ and $\sigma_b$ take values in the same set.
Note that $w_d\equiv 0\pmod2$ in this case.
The sign $\sigma_d$ on such edge takes values in $S_d=\Z^2_2$.
If $\sigma_a=\sigma_b$ then this sign can be presented by
two distinct equivalent elements $\sigma^+_d,\sigma^-_d\in\Z^2_2$.
Let $T'_+$ and $T'_-$ be the trees equipped with the corresponding
signs. We set
\begin{equation}\label{case1}
\mu_{\R}(T)=\mu_{\R}(T'_+)+\mu_{\R}(T'_-).
\end{equation}
If $\sigma_a\neq\sigma_b$ we set $\mu_{\R}(T)=0$.

\item Suppose that $w_a\equiv w_b\equiv 1\pmod2$
and $(x_a,y_a)\not\equiv (x_b,y_b)\pmod2$.
In this case $w_d\equiv1\pmod2$ and
the three sets $S_a, S_b, S_d$ are all distinct.
The sign $\sigma_d\in S_d$ is uniquely determined by
the condition that its equivalence class has common
elements both with the equivalence class $\sigma_a$
and with the equivalence class $\sigma_b$.
Let $T'$ be the tree equipped with this sign.
We set
\begin{equation}\label{case2}
\mu_{\R}(T)=\mu_{\R}(T').
\end{equation}

\item Suppose that one of the weights $w_a$ and $w_b$ is odd
and the other is even.
We may suppose without the loss of generality
that $w_a\equiv 1\pmod2$ and $w_b\equiv 0\pmod2$.
In this case we have $w_d\equiv 1\pmod2$ and $S_a=S_d$ while $S_b=\Z^2_2$.
If the equivalence class $\sigma_a$ contains $\sigma_b$
we set $\sigma_d=\sigma_a$ and
\begin{equation}\label{case3}
\mu_{\R}(T)=2\mu_{\R}(T'),
\end{equation}
where $T'$ is equipped with the sign $\sigma_d$ at $[C,D]$.
If the equivalence class $\sigma_a$ does not contain $\sigma_b$
we set $\mu_{\R}(T)=0$.

\item Suppose that $w_a\equiv w_b\equiv 0\pmod2$. Then $w_d$ is even
and $S_a=S_b=S_d=\Z^2_2$. If $\sigma_a=\sigma_b$ we set
$\sigma_d=\sigma_a$ and
\begin{equation}\label{case4}
\mu_{\R}(T)=4\mu_{\R}(T'),
\end{equation}
where $T'$ is equipped with the sign $\sigma_d$ at $[C,D]$.
If $\sigma_a\neq\sigma_b$ we set $\mu_{\R}(T)=0$.
\end{itemize}
Let $\qqq$ be a signed configuration of $s+g-1$ points and
$h:\Gamma\to\R^2$ be a tropical curve $C=h(\Gamma)$
of genus $g$ and degree $\Delta$ passing via $|\qqq|$.
The real multiplicity of a tropical curve passing through
the signed configuration $\qqq$ is the product
$$\mu_{\R}(C,\qqq)=\prod\limits_{T}\mu_{\R}(T),$$
where $T$ runs over all the components of $\Gamma\setminus h^{-1}(|\qqq|)$.
\end{defn}

\begin{prop}
The real multiplicity $\mu_{\R}(C,\qqq)$ is never greater than and has the
same parity as the multiplicity of $C$ from Definition \ref{multdim2}.
\end{prop}
\begin{proof}
The proposition follows directly from Definition \ref{realmult} by induction.
The multiplicity of a three-valent vertex is odd if and only if all
three adjacent edges have odd weights. This multiplicity is at least 2
if one of the adjacent edges has even weight. This multiplicity
is at least 4 if all three adjacent edges have even weights.
\end{proof}

\IGNORE{
Similarly we have a real tropical enumerative number once the genus $g$,
the degree $\Delta$ and a signed configuration $\qqq$ of $s+g-1$
points is fixed.
\begin{defn}
We define the number $\ntropR(g,\Delta,\qqq)$ to be the number of
irreducible tropical curves of genus $g$ and degree $\Delta$
passing via $|\qqq|$ counted with real multiplicities $\mu_{\R}(C,\qqq)$.
Similarly we define the number $\ntropRr(g,\Delta,\qqq)$ to be
the number of all tropical curves of genus $g$ and degree $\Delta$
passing via $|\qqq|$ counted with real multiplicities $\mu_{\R}(C,\qqq)$.
\end{defn}

\begin{thm}\label{realmain}
Suppose that $\qqq$ is a signed configuration of $s+g-1$ points
in tropically general position. Then there exists a configuration
$\ppp\subset\rtor$ of $s+g-1$ real points in general position
such that $\ncl_{\R}(g,\Delta,\ppp)=\ntropR(g,\Delta,\qqq)$
and $\nclr_{\R}(g,\Delta,\ppp)=\ntropRr(g,\Delta,\qqq)$.

Furthermore,
for every tropical curve $C$ of genus $g$ and degree $\Delta$ passing
through $|\qqq|$ we have $\mu_{\R}(C,\qqq)$ distinct real curves of genus $g$
and degree $\Delta$ passing through $\ppp$. These curves
are distinct for distinct $C$ and are irreducible if
$C$ is irreducible.
\end{thm}
}

\begin{prop}
The multiplicity $\mu_{\R}(C,\qqq)$ equals to the sum of
multiplicities $\mu(C,\{\sigma_E\},\qqq)$ for all real
tropical curves $(C,\{\sigma_E\})$ with the same $C$.
\end{prop}
\begin{proof}
The proposition follows from induction.
In all except for the first case of Definition \ref{realmult}
we have a unique choice for the phase of $[C,D]$ compatible at $C$.
In the first case we have two choices that are accounted in \eqref{case1}.
\end{proof}

\begin{coro}
The number $\ntropR(g,\Delta,\qqq)$ is equal to the sum
of $\mu_{\R}(C,\qqq)$ over all irreducible tropical curves
$C$ passing through $|\qqq|$.

The number $\ntropRr(g,\Delta,\qqq)$ is equal to the sum
of $\mu_{\R}(C,\qqq)$ over all tropical curves
$C$ passing through $|\qqq|$.
\end{coro}

\begin{exa}
Let us choose the signs of $\qqq$ so that every $\sigma_j$
contains $(+,+)\in\Z^2_2$ in its equivalence class.
Let $g=0$ and $\Delta$ be the quadrilateral whose
vertices are $(0,0)$, $(1,0)$, $(0,1)$ and $(2,2)$
as in Example \ref{exa-cusp}.
Then we have $\ntropRr(g,\Delta,\qqq)=3$ for the configuration of
3 points from Figure \ref{tren1} and $\ntropRr(g,\Delta,\qqq)=5$
for the configuration of 3 points from Figure \ref{tren2}.

For other choices of signs of $\qqq$ we can get
$\ntropRr(g,\Delta,\qqq)=1$ for Figure \ref{tren2}
while $\ntropRr(g,\Delta,\qqq)=3$ for Figure \ref{tren1}
for any sign choices.
\end{exa}

\subsection{Counting of real curves by lattice paths}
Theorem \ref{thm1} can be modified to give the relevant count
of real curves. In order to do this we need to
define the real multiplicity of a lattice path
$\gamma:[0,n]\to\Delta$ connecting the vertices $p$ and $q$
once $\gamma$ is equipped with signs.

Suppose $\gamma(j)-\gamma(j-1)=(y_j,x_j)\in\Z^2$, $j=1,\dots,n$.
Let $w_j\in\N$ be the GCD of $y_j$ and $x_j$.
Similarly to the previous subsection we define $S_j$
to be the set obtained from $\Z^2_2$ by taking the quotient
under the equivalence relation $(X,Y)\sim (X+x_j,Y+y_j)$, $X,Y\in\Z_2$.
Let $$\sigma=\{\sigma_j\}_{j=1}^n,\ \sigma_j\in S_j$$
be any choice of signs.

We set
\begin{equation}
\label{remult}
\mu^{\R}_\pm(\gamma,\sigma)=a(T)\mu^{\R}_\pm(\gamma',\sigma')+
\mu^{\R}_\pm(\gamma'',\sigma'').
\end{equation}
The definition of the new paths $\gamma'$, $\gamma''$
and the triangle $T$ is the same as in subsection \ref{paths}.
The sign sequence for $\gamma''$ is $\sigma''_j=\sigma_j,j\neq k,k+1$,
$\sigma''_k=\sigma_{k+1}$, $\sigma''_{k+1}=\sigma_k$.
The sign sequence for $\gamma'$ is $\sigma'_j=\sigma_j,j< k$,
$\sigma'_j=\sigma_{j+1}, j>k$.
We define the sign $\sigma'_{k}$ and
the function $a(T)$ (in a way similar to Definition \ref{realmult}) as follows.

\begin{itemize}
\item
If all sides of $T$ are odd
we set $a(T)=1$ and define the sign $\sigma'_k$ (up to the equivalence)
by the condition that
the three equivalence classes of $\sigma_k$, $\sigma_{k+1}$
and $\sigma'_k$ do not share a common element.
\item
If all sides of $T$ are even we set $a(T)=0$ if $\sigma_{k-1}\neq\sigma_k$.
In this case we can ignore $\gamma'$ (and its sequence of signs).
We set $a(T)=4$ if $\sigma_{k}=\sigma_{k+1}$. In this case we
define $\sigma'_k=\sigma_k=\sigma_{k+1}$.
\item
Otherwise
we set $a(T)=0$ if the equivalence classes of $\sigma_k$ and
$\sigma_{k+1}$ do not have a common element.
We set $a(T)=2$ if they do. In the latter case we define
the equivalence class of $\sigma'_k$ by the condition that
$\sigma_k$, $\sigma_{k+1}$
and $\sigma'_k$ have a common element.
There is one exception to this rule.
If the even side is
$\gamma(k+1)-\gamma(k-1)$
then there are two choices for
$\sigma'_k$ satisfying the above condition.
In this case we replace $a(T)\mu^{\R}_\pm(\gamma')$ in \eqref{remult}
by the sum of the two multiplicities of $\gamma'$ equipped with
the two allowable choices for $\sigma'_k$ (note that this agrees
with $a(T)=2$ in this case).
\end{itemize}
Similar to subsection \ref{paths} we define $\mu^{\R}_\pm(\alpha_\pm)=1$ and
$\mu^{\R}(\gamma,\sigma)=\mu^{\R}_+(\gamma,\sigma)\mu^{\R}_-(\gamma,\sigma)$.
As before $\lambda:\R^2\to\R$ is a linear map injective on $\Z^2$
and $p$ and $q$ are the extrema of $\lambda|_\Delta$.
\begin{thm} \label{thm2}
For any choice of $\sigma_j\in\Z^2_2$, $j=1,\dots,s+g-1$
there exists a configuration of $s+g-1$ of generic points in
the respective quadrants
such that the number of real curves among the $N(g,\Delta)$
relevant complex curves is equal to the number of $\lambda$-increasing
lattice paths $\gamma:[0,s+g-1]\to\Delta$ connecting $p$ and $q$
counted with multiplicities $\mu^{\R}$.

Furthermore,
each $\lambda$-increasing lattice path encodes a number
of tropical curves of genus $g$ and degree $\Delta$
passing via $\qqq$ of total real multiplicity $\mu^{\R}(\gamma,\sigma)$.
These curves are distinct for distinct paths.
\end{thm}


\begin{exa} Here we use the choice $\sigma_j=(+,+)$ so all
the points $z_j$ are in the positive quadrant
$(\R_{>0})^2\subset\rtor$.
The first count of the number $N(0,\Delta)$ from Example \ref{ecusp}
gives a configuration of 3 real points with
5 real curves. The second count gives a configuration
with 3 real curves as the real multiplicity of the last path is 1.
Note also that the second path on Figure \ref{cusp} changes its
real multiplicity if we reverse its direction.

Example \ref{edeg3} gives a configuration of 9 generic points in $\rp^2$
with all 12 nodal cubics through them real.
Example \ref{edeg4} gives a configuration of 12 generic points in $\rp^2$
with 217 out of the 225 quartics of genus 1 real. The path in the
middle of Figure \ref{deg4} has multiplicity 9 but real multiplicity 1.
A similar computation shows that there exists a configuration
of 11 generic points in $\rp^2$ such that 564 out of the 620
irreducible quartic through them are real.
\end{exa}

\subsection{Different types of real nodes and the Welschinger invariant}
\label{rgw}
Let $V\subset\tor$ be a curve defined over $\R$.
In other words it is  a curve
invariant with respect to the involution of complex conjugation.
Suppose that $V$ is nodal, i.e. all singularities of $V$ are
ordinary double points (nodes).

There are three types of nodes of $V$:
\begin{itemize}
\item Hyperbolic. These are the real nodes that locally
are intersections of a pair of real branches. Such nodes
are given by equation $z^2-w^2=0$ for a choice of local
real coordinates $(z,w)$.
\item Elliptic. These are the real nodes that locally
are intersections of a pair of imaginary branches. Such nodes
are given by equation $z^2+w^2=0$ for a choice of local
real coordinates $(z,w)$.
\item Imaginary. These are nodes at non-real points of $V$.
Such nodes come in complex conjugate pairs.
\end{itemize}

\ignore{
Theorems \ref{realmain} and \ref{thm2} can be refined
to account for such distinction.
Namely we can count the number of nodal curves with different
number of nodes separately.
Let us introduce the multiplicity
$$\mu_{\R}^{h,e}(V)=h^{n_h}e^{n_e}\in\Z[h,e]$$
for a real nodal curve $V$.
If $\ppp\subset\rtor$ is a configuration of $s+g-1$ real points
we set
$$N_{\R}^{\operatorname{irr},h,e}(g,\Delta,\ppp)=\sum\limits_V\mu_{\R}^{h,e}(V)$$
where $V$ runs over all irreducible real curves of genus $g$ and degree $\Delta$
passing via $\ppp$.
If $V$ runs over all real curves of genus $g$ and degree $\Delta$
passing via $\ppp$ the resulting sum is denoted with
$N_{\R}^{\operatorname{irr},h,e}(g,\Delta,\ppp)$.

Let $h:\Gamma\to\R^2$ be a tropical curve of genus $g$ and degree $\Delta$
passing via a signed configuration $\qqq$ of $s+g-1$ points.
By Lemma \ref{descompl} each component of
$\Gamma\setminus h^{-1}(|\qqq|)$ is a tree with a single end
at infinity.
We can define the real tropical multiplicity
that takes into account the types of real nodes.
For that we modify
Definition \ref{realmult}
by replacing equations \eqref{case1}, \eqref{case2}, \eqref{case3}
and \eqref{case4} with
$$\mu_{\R}^{h,e}(T)=\frac12(h+e)(\mu_{\R}(T'_+)+\mu_{\R}(T'_-)),$$
$$\mu_{\R}(T)=e^{\frac{m_c-1}{2}}\mu_{\R}(T')$$
$$\mu_{\R}(T)=2\mu_{\R}(T'),\ \text{and}$$
$$\mu_{\R}(T)=4\mu_{\R}(T').$$

Then again we set
$\mu_{\R}^{h,e}(C,\qqq)=\prod\limits_{T}\mu_{\R}^{h,e}(T).$
}

This distinction was used in \cite{W} in order to get a real
curve counting {\em invariant} with respect to the initial
configuration $\ppp\subset\rtor$.
Indeed, let us modify the
real enumerative problem from Definition \ref{realenumerativeproblem}
in the following way. Let $V\subset\tor$ be a real nodal curve.
Let $e(V)$ be the number
of real elliptic nodes of $V$. We prescribe to $V$ the sign
equal to $(-1)^{e(V)}$.
As usual we fix a genus $g$, a degree $\Delta$ and a configuration
$\ppp\in\rtor$ of $s+g-1$ points in general position.
\begin{defn}[see \cite{W}]
We define the number $\ncl_{\R,W}(g,\Delta,\ppp)$ to be the number of
irreducible real curves of genus $g$ and degree $\Delta$
passing via $\ppp$ counted with signs.
Similarly we define the number $\nclr_{\R,W}(g,\Delta,\ppp)$ to be
the number of all real curves of genus $g$ and degree $\Delta$
passing via $\ppp$ counted with signs.
\end{defn}
\begin{thm}[Welschinger \cite{W}]\label{thmW}
If $g=0$ and $\C T_\Delta$ is smooth then the number
$\ncl_{\R,W}(g,\Delta,\ppp)$ does not depend on the choice of $\ppp$.
\end{thm}
Theorem \ref{realmain} can be modified to compute
$\ncl_{\R,W}(g,\Delta,\ppp)$ and $\nclr_{\R,W}(g,\Delta,\ppp)$
for an arbitrary $g$, in particular, in the invariant situation $g=0$.

Let $C\subset\R^2$ be a simple tropical curve.
Recall that Definition \ref{multvert} assigns
a multiplicity $\operatorname{mult}_{V}(C)$
to every 3-valent vertex $V\in C$.
\begin{defn}\label{mult-W}
We define
$$\operatorname{mult}_{V}^{\R,W}(C)=
(-1)^{\frac{\operatorname{mult}_{V}(C)-1}{2}}$$
if $\operatorname{mult}_{V}(C)$ is odd and
$\operatorname{mult}_{V}^{\R,W}(C)=0$
if $\operatorname{mult}_{V}(C)$ is even.

The {\em tropical Welschinger sign $\operatorname{mult}^{\R,W}(C)$}
is the product of $\operatorname{mult}_{V}^{\R,W}(C)$
over all 3-valent vertices of $C$.
\end{defn}

As usual let us fix a genus $g$, a degree $\Delta$
and a configuration $\qqq\subset\R^2$ of $s+g-1$ points
in tropically general position. Define
$\ntropRW(g,\Delta,\qqq)$ to be the number of irreducible tropical
curves of genus $g$ and degree $\Delta$ passing via $\qqq$
counted with the Welschinger sign.
In a similar way define
$\ntropRWr(g,\Delta,\qqq)$ to be the number of all tropical
curves of genus $g$ and degree $\Delta$ passing via $\qqq$
counted with the Welschinger sign.

\begin{thm}\label{main-W}
Suppose that $\qqq\subset\R^2$ is a configuration of $s+g-1$ points
in tropically general position. Then there exists a configuration
$\ppp\subset\rtor$ of $s+g-1$ real points in general position
such that $$\nclr^{\irr}_{\R,W}(g,\Delta,\ppp)=\ntropRW(g,\Delta,\qqq)$$
and $$\nclr_{\R,W}(g,\Delta,\ppp)=\ntropRWr(g,\Delta,\qqq).$$

Furthermore,
for every tropical curve $C$ of genus $g$ and degree $\Delta$ passing
through $\qqq$ we have a number of distinct real curves of genus $g$
and degree $\Delta$ passing through $\ppp$ with the total sum equal
to the Welschinger sign of $C$. These curves
are distinct for distinct $C$ and are irreducible if
$C$ is irreducible.
\end{thm}

\begin{exa}
Let $g=0$ and $\Delta$ be the quadrilateral with vertices
$(0,0)$, $(1,0)$, $(0,1)$ and $(2,2)$.
In Figure \ref{tren1} we have 3 real curves, two of them have the sign $+1$
and one has the sign $-1$.
In Figure \ref{tren2} we have 2 real curves, one of them has the sign $+1$
and one has the sign $0$.
In both cases we have $\ntropRWr(g,\Delta,\qqq)=1$.
\end{exa}

Theorem \ref{thm1} can be adjusted to compute $\ntropRWr(g,\Delta,\qqq)$
by lattice paths. Let $\gamma:[0,n]\to\Delta$ be a lattice path
connecting the vertices $p$ and $q$ of $\Delta$.
Let us introduce
the multiplicity $\nu^{\R}$ inductively, in a manner similar
to our definition of the multiplicity $\mu(\gamma)$.
Namely we set $\nu^{\R}(\gamma)=\nu_+^{\R}(\gamma)\nu_-^{\R}(\gamma)$.
To define $\nu_\pm^{\R}(\gamma)$ we repeat the definition
of $\mu_\pm(\gamma)$ but replace \eqref{pathmult} with
$$\nu^{\R}_\pm(\gamma)=b(T)\nu^{\R}_\pm(\gamma')+
\nu^{\R}_\pm(\gamma'').$$
Here we define $b(T)=0$ if at least one side of $T$ is even
and $b(T)=(-1)^{\#(\Int T\cap\Z^2)}$ otherwise. The paths $\gamma'$,
$\gamma''$ and the triangle $T$ are the same as in the inductive
definition of $\mu_\pm$.

\begin{thm}
\label{pathsW}
For any choice of an irrational linear map $\lambda:\R^2\to\R$
there exists a configuration $\ppp$ of $s+g-1$ of generic points in $\rtor$
such that the number of (not necessarily irreducible)
real curves of genus $g$ and degree $\Delta$ passing
through $\ppp$ counted with the tropical Welschinger sign
is equal to the number of $\lambda$-increasing
lattice paths $\gamma:[0,s+g-1]\to\Delta$ connecting $p$ and $q$
counted with multiplicities $\nu^{\R}$.

Furthermore, there exists a configuration $\qqq\in\R^2$
of $s+g-1$ points in tropical general position such that
each $\lambda$-increasing lattice path encodes a number
of tropical curves of genus $g$ and degree $\Delta$
passing via $\qqq$ with the sum of signs equal to $\nu^{\R}(\gamma)$.
These curves are distinct for distinct paths.
\end{thm}

\begin{exa}
In example \ref{ecusp} we
have $\nclr_{\R,W}(0,\Delta,\ppp)=1$ for some choice of $\ppp$.
In example \ref{edeg3} we have
$\nclr_{\R,W}(0,\Delta_3,\ppp)=8$ for some choice of $\ppp$.
In example \ref{edeg4} we have
$\nclr_{\R,W}(1,\Delta_4,\ppp)=93$ for some choice of $\ppp$.
In all these examples we have $\ncl_{\R,W}=\nclr_{\R,W}$
as there are no reducible curves of these genera and degree.
Note that by Theorem \ref{thmW}
the number $\ncl_{\R,W}(0,\Delta_3,\ppp)=8$
does not depend on the choice of $\ppp$.
\end{exa}

The following observation is due to Itenberg, Kharlamov
and Shustin \cite{IKS}. If $\lambda(x,y)=x-\epsilon y$ and
$\Delta=\Delta_d$ or $\Delta=[0,d_1]\times [0,d_2]$
then $\nu^{\R}(\gamma)\ge 0$
for any $\lambda$-increasing path $\gamma$ and, furthermore,
any tropical curve encoded by $\gamma$ by Theorem \ref{pathsW}
has a non-negative tropical Welschinger sign. It is easy to show
that there exist $\lambda$-increasing paths that encode
irreducible tropical curves of non-zero tropical Welschinger signs.
We get the following corollary for any $d,d_1,d_2\in\N$.
\begin{cor}
For any generic configuration $\ppp\subset\rp^2$ of $3d-1$ points
there exists an irreducible rational curve $V\subset\rp^2$ of degree $d$
passing through $\ppp$.

For any generic configuration $\ppp\subset\rp^1\times\rp^1$ of
$2d_1+2d_2-1$ points
there exists an irreducible rational curve $V\subset\rp^1\times\rp^1$
of bidegree $(d_1,d_2)$ passing through $\ppp$.
\end{cor}

With the help of Theorem \ref{pathsW}
Itenberg, Kharlamov and Shustin in \cite{IKS} have obtained a non-trivial lower
bound for the number of such rational curves.
In particular,
they have shown that for any generic configuration $\ppp$ of $3d-1$ points
in $\rp^2$ there exists at least $\frac{d!}{2}$ rational curves of degree $d$
passing via $\ppp$.

\section{Proof of the main theorems}\label{prma}
\subsection{Complex amoebas in $\R^2$ and the key lemma}
Let $\qq=\{q_1,\dots,q_n\}\subset\tor$ be generic points in the
sense of Proposition \ref{cgeneric} (i.e. generic classically).
Suppose in addition that
the points $$p_1=\Log(q_1),\dots,p_{n}=\Log(q_{n})\in\R^2$$
are in general position tropically.
Denote $\ppp=\{p_1,\dots,p_n\}\subset\R^2$.
Recall that we fix a Newton polygon $\Delta\subset\R^2$ and a genus
$g\in\N$. There are $\nclr(g,\Delta)$ holomorphic curves
passing through $\qq$ as long as $n=s+g-1$.
By Proposition \ref{tropfinite} there are finitely many tropical curves
\begin{equation}\label{gammy}
C_1,\dots,C_m\subset\R^2
\end{equation}
of genus $g$ with the Newton polygon $\Delta$ and passing
through $\ppp$. Note that $m$ depends
on the choice of the points $x_j$ (unlike the number
$\ntropr(g,\Delta)\ge m$ for which we are proving invariance).

\begin{prop}\label{jt-hol}
For generic $t$ we have $\nclr(g,\Delta)$ $J_t$-holomorphic
curves passing through $\qq$.
\end{prop}
\begin{proof}
Indeed, this number
equals to the number of holomorphic
(i.e. $J_e$-holomorphic) curves through
the points $H_t(p_1),\dots,H_t(p_n)$. (These points are
in general position for generic $t$ since they are for $t=e$.)
\end{proof}

Let $\NN_{\epsilon}(C_j)$ be the $\epsilon$-neighborhood
(in the sense of the standard Euclidean metric on $\R^2$)
of $C_j$ for $\epsilon>0$.
Recall (see \cite{GKZ}) that the {\em amoeba} of a curve
$V\subset\tor$ is its image $$\am=\Log(V)\subset\R^2.$$
Note that if $V$ is a $J_t$-holomorphic
curve then $\Log(V)$ can be obtained from the amoeba of some
holomorphic curve by the $\log(t)$-contraction. This allows us to speak
of the Newton polygons of $J_t$-holomorphic curves.

\begin{prop}\label{trvnutri}
If $V$ is a $J_t$-holomorphic curve whose Newton polygon is
$\Delta$ then its amoeba $\am=\Log(V)$ contains a tropical curve
$C$ with the same Newton polygon $\Delta$.
\end{prop}
\begin{proof}
If $t=e$, i.e. $C$ is holomorphic with respect to the standard
holomorphic structure, then the statement follows from the theorem
of Passare and Rullg{\aa}rd \cite{PR}.
Recall (see \cite{PR}) that if $V$ is a complex curve
defined by a polynomial $f(z,w)=\sum\limits_{j,k} a_{j,k}z^jw^k$
then the {\em spine} of its amoeba $\am$ is a tropical
curve defined by a tropical polynomial
$N^{\infty}_f(x,y)=``\sum\limits_{j,k} b_{j,k}x^jy^k ",$
where $$b_{j,k}=\frac{1}{(2\pi i)^2}\int\limits_{\Log^{-1}(r)}
\log|f(z,w)|\frac{dz}{z}\frac{dw}{w}$$
and $r\in \R^2\setminus\am$ is any point such that its {\em index} is $(j,k)$.
If there are no points in $\R^2\setminus\am$ of index $(j,k)$ then
the monomial $x^jy^k$ is omitted from $N^{\infty}_f$.

It is shown in \cite{PR} that the tropical hypersurface defined
by $N^{\infty}_f$ is contained in $\am$. Clearly, the Newton
polygon of $N^{\infty}_f$ is $\Delta$.
To finish the proof we note that the image of a tropical
curve under a homothety is tropical with the same Newton polygon.
\end{proof}

With the help of Proposition \ref{jt-hol}
Theorem \ref{main} follows from the following two lemmas.

\ignore{
As before, we fix a number $g\in\Z$ as the genus, a lattice polygon $\Delta\in\R^2$
as the degree and choose a configuration $\qq\in\tor$ of $s+g-1$ points in general
position.
}
\begin{lem}\label{keylem-a}
For any $\epsilon>0$ there exists $T>1$
such that if $t>T$ and $V$ is a $J_t$-holomorphic curve of genus $g$,
degree $\Delta$ and passing through
$\qq$ then its amoeba $\Log(V)$ is contained in
the $\epsilon$-neighborhood $\NN_{\epsilon}(C_j)$
of $C_j$ for some $j=1,\dots,m$.
\end{lem}
\begin{lem}\label{keylem-b}
For a sufficiently small $\epsilon>0$ and a sufficiently
large $t>0$
the multiplicity $\mult(C_j)$ of each $C_j$ from \eqref{gammy}
(see Definition \ref{multdim2}) is equal to the
number of the $J_t$-holomorphic curves $V$ of genus $g$ and
degree $\Delta$ passing through $\qq$
and such that $\Log(V)$ is contained in
$\NN_{\epsilon}(C_j)$.
Furthermore, if $C_j$ is irreducible then any
$J_t$-holomorphic curve $V$ of genus $g$ and degree $\Delta$
passing through $\qq$ with $\Log(V)\subset\NN_{\epsilon}(C_j)$
is irreducible while if $C_j$ is reducible than any such curve $V$
is reducible.
\end{lem}

\subsection{Proof of Lemma \ref{keylem-a}}
A holomorphic curve $V\subset\tor$ is given by a polynomial
$$F(z_1,z_2)=\sum\limits_{(j,k)\in\Delta} a_{j,k}z_1^jz_2^k.$$
To a curve $V\subset\tor$
we associate its {\em tropicalization}
$V^{\trop}\subset\R^2$ given by the tropical polynomial
$$F^{\operatorname{trop}}(y_1,y_2)=
\max\limits_{(j,k)\in\Delta} (jy_1+ky_2+\log|a_{j,k}|).$$

\begin{lem}\label{lemeps}
The amoeba $\Log(V)$ is contained in the
$\delta$-neighborhood
of $V^{\trop}$ (with respect to the Euclidean metric in $\R^2$),
where $$\delta=\log(\#(\Delta\cap\Z^2)-1).$$
\end{lem}
\IGNORE{
We define $\Psi_{\Delta}(y_1,y_2)=0$ if $(y_1,y_2)\in V^{\trop}_t$.
If $(y_1,y_2)\notin V^{\trop}$ then there exists
$(j',k')\in\Delta\cap\Z^2$ such that for any
$(j,k)\in\Delta\cap\Z^2$, $(j,k)\neq(j',k')$
we have
$$j'y_1+k'y_2+\log|a_{j',k'}| > jy_1+ky_2+\log|a_{j,k}|.$$
The index $(j',k')$ depends only on the component
of $\R^2\setminus V^{\trop}$ that contains $(y_1,y_2)$.
In this component we set
$$\Psi_{\Delta}(y_1,y_2)=j'y_1+k'y_2
-\max\limits_{(j,k)\neq (j',k')} (jy_1+ky_2).$$
Clearly, $\Psi_{\Delta}$ is continuous and such that
$\Psi_{\Delta}^{-1}(0)=V^{\trop}$. To check that
$\Psi_{\Delta}(\Log_t(V))\le\log_t(\#(\Delta\cap\Z^2))$
we note that if $(y_1,y_2)\in\Log_t(V)$
then we have a triangle inequality for
$(z_1,z_2)\in V\subset\tor$ such that $y_1=\log_t|z_1|$, $y_2=\log_t|z_2|$:
$$|a_{j',k'}z_1^{j'}z_2^{k'}|\le
|\sum \limits_{(j,k)\neq (j',k')}
a_{j,k}z_1^{j}z_2^{k}|.$$
Let us apply $\log_t$ to both sides of this inequality:
$$j'y_1+k'y_2+\log_t|a_{j',k'}|=\log_t|a_{j',k'}z_1^{j'}z_2^{k'}|\le$$
$$\log_t|\sum \limits_{(j,k)\neq (j',k')}
a_{j,k}z_1^{j}z_2^{k}|\le
\log_t((\#(\Delta\cap\Z^2)-1)\times\max \limits_{(j,k)\neq (j',k')}
|a_{j,k}z_1^{j}z_2^{k}|)=$$
$$\log_t|\#(\Delta\cap\Z^2)-1)|+
\max\limits_{(j,k)\neq (j',k')} (jy_1+ky_2+\log_t|a_{j,k}|).
$$
}
\begin{proof}
Suppose that $(y_1,y_2)$ is not contained in the
$\delta$-neighborhood of $V^{\trop}$. Then there exists $(j',k')$
such that
\begin{equation}\label{logeq}
j'y_1+k'y_2+\log|a_{j',k'}| > jy_1+ky_2+\log|a_{j,k}|+\delta
\end{equation}
for any $(j,k)\neq(j',k')$.
Indeed, the distance from
$(y_1,y_2)$ to the line
$j'y_1+k'y_2+\log|a_{j',k'}| = jy_1+ky_2+\log|a_{j,k}|$
is greater than $\delta$ by the hypothesis and the norm of the gradient of the
function $(j-j')y_1+(k-k')y_2$ is at least 1 (as $j,j',
k,k'$ are all integers).

Suppose that
$(z_1,z_2)\in V\subset\tor$ is such that $\Log(z_1,z_2)=(y_1,y_2)$.
Since $\sum\limits_{(j,k)\in\Delta} a_{j,k}z_1^jz_2^k=0$
the triangle inequality implies that
$$|a_{j',k'}z_1^{j'}z_2^{k'}|\le
|\sum \limits_{(j,k)\neq (j',k')}
a_{j,k}z_1^{j}z_2^{k}|.$$
Let us apply $\log$ to both sides of this inequality:
$$j'y_1+k'y_2+\log|a_{j',k'}|=\log|a_{j',k'}z_1^{j'}z_2^{k'}|\le$$
$$\log|\sum \limits_{(j,k)\neq (j',k')}
a_{j,k}z_1^{j}z_2^{k}|\le
\log((\#(\Delta\cap\Z^2)-1)\times\max \limits_{(j,k)\neq (j',k')}
|a_{j,k}z_1^{j}z_2^{k}|)=$$
$$\delta+
\max\limits_{(j,k)\neq (j',k')} (jy_1+ky_2+\log|a_{j,k}|).
$$
Thus a point from the amoeba has to be in the $\delta$-neighborhood
of $V^{\trop}$.
\end{proof}
\begin{coro}\label{coroeps}
The amoeba $\Log(V_t)$ of a $J_t$-holomorphic curve $V_t=H_t(V)$
is contained in the $\delta$-neighborhood of some tropical
curve in $\R^2$, where
$\delta=\log_t (\#(\Delta\cap\Z^2))$.
\end{coro}
\begin{proof}
The corollary is obtained by applying the $\log t$-contraction
to Lemma \ref{lemeps} since $\Log(V_t)=\Log_t(V)=\frac{\Log(V)}{\log t}$.
Clearly,
$\frac{V^{\trop}}{\log t}\subset\R^2$ is a tropical curve.
\end{proof}

Let $V_k\subset\tor$, $k\in\N$, be a sequence of curves passing through
$\qq$ and such that $V_k$ is a
$J_{t_k}$-holomorphic curve for some $t_k>0$, where $t_k\to
\infty$, $k\to\infty$. As in the previous subsection we assume
that the holomorphic curve $H_{t_k}^{-1}(V_k)$ is of genus
$g$ and has the Newton polygon $\Delta$ for each $k$.
Denote with $\am_k=\Log(V_k)$ the amoeba of $V_k$.
Proposition \ref{trvnutri} ensures (after applying
the $\log t_k$-contraction) that there exists
a tropical curve $\ce_k\subset\am_k$.

\begin{prop}\label{1lem}
There is a subsequence $V_{k_\alpha}$, $\alpha\in\N$, such that
the sets $\am_{k_\alpha}\subset\R^2$ converge in the Hausdorff
metric in $\R^2$ to some tropical curve $C_j$ from \eqref{gammy}.
\end{prop}

\begin{proof}
By Proposition
\ref{compact} we can extract a subsequence from $\ce_k$
which converges to a tropical curve $C$.
To prove the proposition it suffices to show that
$C$ is a tropical curve passing
through $\ppp$ of genus $g$ whose Newton
polygon is $\Delta$.
Proposition \ref{compact} and
Corollary \ref{coroeps} ensure convergence in the Hausdorff metric
in $\R^2$.

We have $C\supset\ppp$ since $V_k\supset\qq$ and thus $\am_k\supset\ppp$.
The degree of $C$ is a subpolygon $\Delta'\subset\Delta$
since $C$ is the limit of curves of degree $\Delta$.
We want to prove that $\Delta'=\Delta$.

Choose a disk
$D_R\subset\R^2$ of radius $R$ so large that $D_R$ contains all
vertices of $C$. Furthermore, making $R$ larger if needed
we may assume that the extension of the exterior edges of
$C\cap D_R$ beyond $D_R$ do not intersect.

\ignore{
Then $C\setminus D_R$ consists
of $s'=\#(\Delta'\cap\Z^2)$ rays.
}
The Newton polygon of $\ce_k$ is $\Delta$. Therefore, it has
$s$ ends. By Proposition \ref{compact} the intersection $\ce_k\cap D_R$
is an approximation of $C\cap D_R$. Therefore, for a large $k$
we have $\ce_k\supset \ppp'$ where $\ppp'$ is a configuration
of $s+g-1$ points in tropically general position obtained
by a small deformation of $\ppp$.
We have $\Delta'=\Delta$ if and only if
$\ce_k^{\trop}\setminus D_R$ is a disjoint union of rays (each going to $\infty$).
If not, $\ce_k\setminus D_R$ has a bounded edge connecting a point
of $\ce_k\cap\dd D_R$ with a vertex of $\ce_k$ outside of $D_R$.
A change of the length of this edge produces a deformation of $\ce_k$
such that all curves in the family pass via $\ppp'$.
This contradicts to the tropical general position of $\ppp'$.

\ignore{
We can treat $V_k^{\trop}\setminus D_R$ as
a cobordism which connects $s'$ points at $V_k^{\trop}\cap\dd D_R$ to
$s$ points at infinity. We have to prove that this
cobordism is trivial. But if not then at least two points from
$V_k^{\trop}\cap\dd D_R$ belong to the same component of the cobordism.
}

Note that the genus of $C$ cannot be smaller than $g$,
otherwise the configuration $\ppp$ is not in general
position. The genus of curves $\ce_k$ may be larger than
$g$ even though the genus of $\ce_k$ is $g$. However, the genus
of their limit $C$ cannot be larger than $g$
by Proposition \ref{genusupdown} as $C$ can be presented as
the image under $\Log$ of the limit of a subsequence of $V_k$.
Therefore, the genus of $C$ is $g$.
Thus $C$ has to be one of $C_k$ from \eqref{gammy}.
\end{proof}

Since Proposition \ref{1lem} can be applied to any infinite subsequence
of $V_k$ we have the following corollary.
\begin{coro}\label{cor1lem}
The sequence $\am_{k}\subset\R^2$, $k\in\N$
can be presented as a union of $m$
converging subsequences $\am^j_k\subset\R^2$
such that either the number of terms in $\am^j_k$ is finite
or $\am^j_k$ converges in the Hausdorff metric in $\R^2$
with $\lim\limits_{k\to+\infty}\am^j_k=C_j$, $j=1,\dots,m$.
\end{coro}
\begin{proof}
By Proposition \ref{1lem} for any $\delta>0$
the number of indices $k\in\N$ such that the Hausdorff
distance from $\am_k$ to $C_j$ is greater than $\delta$
for every $j=1,\dots,m$ is finite.
\end{proof}

To deduce Lemma \ref{keylem-a} from Corollary \ref{cor1lem}
Suppose that the amoebas $\am_k$ of $J_{t_k}$-holomorphic
curves $V_k$ converge to a tropical curve $C=C_j$, $j=1,\dots,m$ for $k\to\infty$.
Then the spines $\ce_k\subset\am_k$ also converge to $C$.

Note that
the number of edges of $C$ is not greater than that of $\ce_k$.
Some edges of the 1-cycle $\ce_k$ tend to the ``corresponding"
edges of $C$. The remaining edges of $\ce_k$ are vanishing:
their length tends to zero when $k\to\infty$.
Lemma \ref{keylem-a} follows from the following Proposition.
\begin{prop}
There exists $T>0$ and a function $\tilde\delta:[T,+\infty)\to\R$
such that for every edge $E_k$ of $\ce_k$ whose length is higher than $\delta$
there exists an edge $E$ of $C$ parallel to $E_k$ and within
$\tilde{\delta}(t_k)$-distance (in the Hausdorff metric in $\R^2$) from $E$.
\end{prop}
\begin{proof}
Since $\Delta$ is a bounded polygon we have a finite number
of possibilities for the slopes of the edges of $\ce_k$.
Thus each edge of $C$ gets approximated by a parallel edge $E_k$ of
$\ce_k$ when $t_k\to\infty$.

Suppose that $E\subset C$
is an edge containing a point $p\in\ppp$. The distance
between the parallel lines containing $E$ and $E_k$ cannot
be more than $\delta_k=\log_{t_k}(\#(\Delta\cap\Z^2))$ by Corollary \ref{coroeps}
(assuming that $t_k$ is so large that $d(p,C\setminus E)>\delta_k$).

Recall that $C\setminus\ppp$ is a disjoint union of trees with only one end
at infinity.
Since the number of edges of $\ce_k$ is bounded from above
it suffices to prove that the length of vanishing edges of $\ce_k$
is uniformly bounded from above by a quantity tending to zero when
$k\to\infty$. Every vanishing edge is contained in a small
neighborhood $U\subset\R^2$ of a vertex $v$
of the 1-cycle $C$ for large $t_k$.
Recall that since $\ppp$ is in general position
the 1-cycle $C$ is simple.

Suppose $v\in C$ is a 3-valent vertex and $U\ni v$ is a small open disk
around $v$. As in the proof of Proposition \ref{genusupdown}
it is easy to see that $\Log^{-1}(U)\cap V_k$ is connected
and homeomorphic to a pair-of-pants, i.e.
a sphere punctured three times. Indeed, every
component of $\Log^{-1}(U)\cap V_k$ has at least three ends by the
maximum principle. Thus the Euler characteristic of each such component
is at most $-1$ and strictly less than $-1$ unless $\Log^{-1}(U)\cap V_k$
is a pair-of-pants.
Our claim follows since the genus of $V_k$ is $g$
(which coincides with the genus of $C$).

Unless $U\cap\ce_k$ is a union of three rays emanating from $v$ it
must contain a cycle. Since $\ce_k$ is a deformational retract
of $\am_k$ the intersection $U\cap\ce_k$ then also contains a cycle.
Let $w\in U\setminus\am_k$ be a point inside of this cycle and let
$w'\in U\setminus\am_k$ be any other point such that $[w,w']\cap\am_k$
is non-empty and connected. Note that the line $L\subset\R^2$
containing $[w,w']$ also must intersect $U\cap\am_k$ outside
the interval $[w,w']$ by topological reasons. The pull-back
$\Log^{-1}(U\setminus [w,w'])\cap V_k$ cannot be disconnected.
Otherwise this would contradict to the maximum principle for
the harmonic function $\pi_L\circ\Log|_{V_k}$, where $\pi_L:\R^2\to\R$
is the linear projection in the direction of $L$ (recall that
$V_k\cap\Log^{-1}(U)$ has three ends and $\am_k$ is contained
in the $\delta_k$-neighborhood of $\ce_k$).
But if $\Log^{-1}(U\setminus [w,w'])\cap V_k$ is connected then
the genus of $V_k\cap\Log^{-1}(U)$ is positive which is also
a contradiction. Thus $\ce_k$ cannot have vanishing edges
near 3-valent vertices of $C$.

Suppose now that $v\in C$ is a 4-valent vertex of the 1-cycle $C$.
Let $U\ni v$ is a small open disk around $v$. The pull-back
$V_k\cap\Log^{-1}(U)$ is a union of two components $V^U_1$
and $V^U_2$.
As in the 3-valent case the intersection
$\ce_k\cap U$ cannot have cycles as
it would lead to a contradiction with the maximum principle.
Therefore, $U\cap\ce_k$ is a tree with 4 ends
and thus may contain not more than one vanishing edge $E$.

Suppose that the length of $E$ is
greater that $2\delta_k/s$, where $s$ is the minimal value
of the sine of the angle between two distinct slopes of the
edges of $\ce_k$.
Let $L\subset\R^2$ be the line passing through the midpoint
of $E$ and parallel to one of the two edges of $C$ passing through $V$.
Corollary \ref{coroeps} combined with
our assumption on the length of $E$ implies
that $L\cap\am_k$ is compact and therefore $\Log(V^U_1)$
cannot be locally concave which contradicts to the maximum
principle for $\pi\circ\Log|_{V^U_1}$ for some linear projection
$\pi:\R^2\to\R$ (cf. Lemma 1 from \cite{Mi}).
\end{proof}

\subsection{Proof of Lemma \ref{keylem-b}}
Let $C\subset\R^2$ be one of the tropical curves $C_j$ from \eqref{gammy}
and $\mult(C)$ be its multiplicity from Definition \ref{multdim2}.
Denote with $\subdiv_C$ the lattice subdivision of $\Delta$ dual to $C$.
Let
\begin{equation}\label{sootvtrop}
f_{\trop}(y)=\sum\limits_{j\in\Delta\cap\Z^2} \beta_j y^j
\end{equation}
be a tropical polynomial that defines $C$
and such that $f_{\trop}$ includes
all monomials of indices $j\in\Delta\cap\Z^2$ so that
$j\mapsto -\beta_j$ is a (non-strictly) convex function
defined on $\Delta\cap\Z^2$.
To get rid of the ambiguity in the choice of $f_{\trop}$
we choose a ``reference" index $j_0\in\Delta\cap\Z^2$
among the vertices of $\Delta$
and assume that $\beta_{j_0}=0$.

By Proposition \ref{ctropmain} there are
$\mult(C)/\mu_{\operatorname{edge}}(C,\Log(\qq))$
simple complex tropical curves
projecting to $C$ and passing via $\qq$.
We shall see that each of them gives rise to $\mu_{\operatorname{edge}}(C,\Log(\qq))$
distinct nearby $J_{t_k}$-holomorphic curves of degree $\Delta$, genus $g$
and passing via $\qq$ for large $t_k>1$ and this exhausts all $J_{t_k}$
holomorphic curves with this property in a small neighborhood of
$\Log^{-1}(C)\subset\tor$.

Let $V_\infty$ be any complex tropical curve passing via $\qq$
and such that $\Log(V_\infty)=C$.
By Proposition \ref{ctrop-poly} the curve $V_\infty$
defines the coefficients $a_j\in\C$ for the vertices
$j$ of $\subdiv_C$ once we set $a_{j_0}=1$.
Note that since $C\supset\qqq$ and $\qqq$ is in tropically
general position the number of ends of $C$ is $s$ and therefore
$$\dd\Delta\cap\Z^2\subset\operatorname{Vert}(\subdiv_C).$$
We have $\log|a_j|=\beta_j.$
It turns out that these coefficients $a_j$, $j\in\operatorname{Vert}(\subdiv_C)$,
depend only on $C$ and $\qq$ as the next proposition shows. Note that
unlike the situation in Proposition \ref{ctrop-thm} not all
points of $\Delta\cap\Z^2$ are vertices of $\subdiv_C$
and thus we may have distinct complex tropical curves with a given set
of the coefficients $a_j$, $j\in\operatorname{Vert}(\subdiv_C)$.
\begin{prop}\label{coeff-aj}
The coefficients $a_j$, $j\in \operatorname{Vert}(\subdiv_C)$,
do not depend on the choice of $V_\infty$ as long as we set $a_{j_0}=1$.
Thus these coefficients depend only on $C$ and $\qq$.
\end{prop}
\begin{proof}
By Lemma \ref{descompl} each component $K$ of $\Gamma\setminus h^{-1}(\qq)$
is a tree which contains one end at infinity.
As in the proof of Proposition \ref{ctropmain} we proceed
inductively cutting the branches of the tree $K$.

Let $A,B\in\qq$ be two points that connect to the same
3-valent vertex in $K$ as in Figure \ref{komp-izm}.
Then the points $A,B$ are contained in the edges of $C$ dual
to the edges $[j,j']$, $[j',j'']$ such that the triangle with
the vertices $j,j',j''$ belongs to $\subdiv_C$.
The ratios $\frac{a_j}{a_{j'}}$
and $\frac{a_{j'}}{a_{j''}}$ are determined since the
points $A$ and $B$ have to be contained in the curve.
Therefore, we know the ratio $\frac{a_j}{a_{j''}}$ and we can proceed
further by induction.
\end{proof}

Consider the polynomial
$$f_t(z)=\sum\limits_{j\in\operatorname{Vert}(\subdiv_C)} \arg(a_j)t^{\log|a_j|}z^j,$$
$z\in\tor$.
The sum here is taken only over the vertices $j$ of $\subdiv_C$.
Define $V_t\subset\tor$ by
$$V_t=H_t(\{z\in\tor\ |\ f_t(z)=0\}).$$
For each $t>1$ the curve $V_t$ is $J_t$-holomorphic.
For large values of $t$ we have $\Log(V_t)\subset\NN_\epsilon(C)$

For large values of $t$ we may consider the curve $V_t$ as
an {\em approximate} solution to the problem of finding
the $J_t$ holomorphic curves of genus $g$ and degree $\Delta$
via $\qq$. Indeed, generically we expect $V_t$ to be smooth
(and therefore of a ``wrong" genus $\#(\Delta\cap\Z^2)\ge g$) and
not to contain $\qq$.
However it is close to a (singular) curve of genus $g$
and is very close to the configuration $\qq$.
We need to find a genuine solution near this approximate one
for large values $t>>1$.

Recall that the amoebas of the ($J_t$-holomorphic)
curves we are looking for have to be contained
in a small neighborhood of $C$ while the curves themselves
have to contain $\qq$.
For large $t>>1$ this implies that such curve can be presented
in the form $V_t^\zeta=H_t(\{f_t^\zeta(z)=0\})$, where
$$f_t^\zeta(z)=\sum\limits_{j\in\Delta\cap\Z^2}
\arg(\zeta_{j})t^{\log|\zeta_{j}|}z^j$$
and $\zeta\in\C^n$, $n=\#(\Delta\cap\Z^2)-1$, are such
that $|\zeta_j-a_j|<\epsilon'_j$ for
$j\in\operatorname{Vert}(\subdiv_C)$ while
$\log|\zeta_j|-\beta_j<\epsilon'_j$
for $j\notin\operatorname{Vert}(\subdiv_C)$.
Here $\epsilon'_j>0$ is some collection of small numbers.
All $\zeta$ that comply to these conditions form
a polydisc $\DD\subset\C^n$.

\begin{prop}
If $t$ is sufficiently large and ${\mathcal V}_t\subset\tor$
is a $J_t$-holomorphic curve of genus $g$ such that
${\mathcal V}_t\supset\qq$ and $\Log_t({\mathcal V}_t)\subset\NN_\epsilon(C)$
then there exists $\zeta\in\DD$ such that ${\mathcal V}_t=H_t(V^\zeta_t)$.
\end{prop}
Recall that by Lemma \ref{keylem-a} the image $\Log_t({\mathcal V}_t)$
has to be contained in the $\epsilon$-neighborhood of one of $C_j$
from \eqref{gammy}.
\begin{proof}
Suppose that $\log|\zeta_j|-\beta_j\ge\epsilon'_j$ for some
$j\notin\operatorname{Vert}(\subdiv_C)$ and
this holds for a sequence of values of $t$ going to $+\infty$.
Then there exists a subsequential limit of $\Log_t({\mathcal V}_t)$
of genus higher than $g$ (as the limiting tropical curve has an extra
complement component corresponding to $j$). By Proposition \ref{genusupdown}
the genus of ${\mathcal V}_t$ is also higher than $g$.

Suppose that $|\zeta_j-a_j|\ge\epsilon'_j$ for some $j$.
Tracing back the proof of Proposition \ref{coeff-aj} we see
that then there exists an edge $[j',j'']\in\Xi$ corresponding
to a point $q\in\qq$ such that
$|\frac{a_{j'}}{a_{j''}}-\frac{\zeta_{j'}}{\zeta_{j''}}|$
is bounded from below by a positive constant not depending on $t$.
But then the curve ${\mathcal V}_t$ is disjoint from $q$ for sufficiently
large $t$.
\end{proof}

We cover $\R^2$ with open sets $U(\Delta')$
corresponding to the polygons $\Delta'$
from the subdivision $\subdiv_C$ in the following way.
\begin{itemize}
\item
If $\Delta'\in\subdiv_C$ is a 2-polygon then it is dual
to a vertex $p_{\Delta'}\in C$. We choose $U(\Delta')$
to be a small open disk centered at $p_{\Delta'}$.
\item
If $\Delta'$ is an edge of $\subdiv_C$ then it is
dual to an edge $e_{\Delta'}\subset C$ connecting
two vertices $p_{\Delta_1},p_{\Delta_2}\in C$,
where $\Delta_1,\Delta_2\in\subdiv_C$ are the two 2-polygons
adjacent to the edge $e_{\Delta'}$. We choose
$U(\Delta')$ to be a small regular open neighborhood
of $e_{\Delta'}\setminus(U(\Delta_1)\cup U(\Delta_2))$.
\item
According to our previous choice of $U(\Delta')$
the curve $C$ is a deformational retract of the union
${\mathcal U}=\bigcup\limits_{\Delta'} U(\Delta')$ where $\Delta'$ runs
over positive-dimensional polygons from $\subdiv_C$.
Therefore there is a natural bijection between the components
of $\R^2\setminus{\mathcal U}$ and $\operatorname{Vert}(\subdiv_C)$.
We choose $U(\Delta')$ to be a small open neighborhood of the
component of $\R^2\setminus{\mathcal U}$ corresponding to $\Delta'$
if $\Delta'\in\operatorname{Vert}(\subdiv_C)$.
\end{itemize}

The following {\em patchworking principle} (due to Viro \cite{Vi})
can be used to localize the problem, i.e. to reduce it to an individual
problem in every $U(\Delta')$:
{\em
deformation of $\zeta_j$ with $j\notin\Delta'$
has little effect on $V^\zeta_t\cap\Log^{-1}_t(U(\Delta'))$ for large $t$.
}
To state this principle formally in this tropical set-up let us consider
for every 2-dimensional $\Delta'\in\subdiv_C$
a translation $$\phi_{\Delta'}:\R^2\to\R,$$
$\phi_{\Delta'}(y_1,y_2)=(y_1+b_1,y_2+b_2)$, $b_1,b_2\in\R$,
such that $\phi_{\Delta'}(p_{\Delta'})=0\in\R^2$.
Similarly, if $\Delta'\in\subdiv_C$ is 1-dimensional
we choose a translation $\phi_{\Delta'}:\R^2\to\R^2$
so that there exists $p_{\Delta'}\in C\cap U(\Delta')$
with $\phi_{\Delta'}(p_{\Delta'})=0$. If $\Delta'$ is a vertex
of $\subdiv_C$ we choose the translation $\phi_{\Delta'}:\R^2\to\R^2$
so that $\phi^{-1}_{\Delta'}\in U(\Delta')$.
Consider the lifting $$\Phi_{\Delta',t}:\tor\to\tor$$
defined by $\Phi_{\Delta',t}(z_1,z_2)=(t^{b_1}z_1,t^{b_2}z_2)$,
$t>1$.
We have $$\phi_{\Delta'}\circ\Log=\Log_t\circ\Phi_{\Delta',t}.$$

Note that (since $j\mapsto -\log|a_j|$ is convex on
$\operatorname{Vert}(\subdiv_C)$) we have
$$f^\zeta_t\circ\Phi^{-1}_{\Delta',t}(z)=t^c\sum\limits_j \arg(\zeta_j)t^{\eta_j}z^j,$$
where $\eta_j<\delta_{\Delta'}$ if $j\notin\Delta'$ and $\zeta\in\DD$
for some constant $\delta_{\Delta'}<0$ while $\eta_j$ for $j\in\Delta'$
is sufficiently close to zero (as $\DD$ is small).
Here $c\in\R$ is the constant term of the affine-linear
map $j\mapsto -\beta_j$ on $\Delta'\cap\Z^2$.

If $j\notin\Delta'$ and $\Log_t(z)\in\phi_{\Delta'}(U(\Delta'))$
we have $|\arg(\zeta_j)t^{\eta_j}z^j|\le t^{\delta'(\Delta')}$,
for some constant $\delta'(\Delta')>0$
since $U(\Delta')$ is chosen so that the tropical monomials
in $f^{\trop}$ corresponding to lattice points of $\Delta'$ dominate
the monomials corresponding to lattice points of $\Delta\setminus\Delta'$
by some definite amount.
Therefore we have a uniform upper bound for
$|\sum\limits_{j\in\Delta\setminus\Delta'} \arg(\zeta_j)t^{\eta_j}z^j|$
independent of $\zeta\in\DD$ which tends to zero when $t\to+\infty$.
Let $$V^\zeta_{\Delta',\infty}\subset\tor$$ be the image under
$\Phi_{\Delta',t}$ of the complex tropical
curve (with the Newton polygon $\Delta'$) given by complex
tropical coefficients $\zeta_j$, where $j\in\Delta'$ is a vertex
of the subdivision of $\Delta'$ defined by $j\mapsto-\log|\zeta_j|$,
as in Proposition \ref{ctrop-thm}.
Then $$\Phi_{\Delta',t}(V^\zeta_t\cap\Log^{-1}_t(U(\Delta')))\subset\tor$$
tends (as a sequence of subsets of $\tor$ with respect to the Hausdorff metric)
to $V^\zeta_{\Delta',\infty}\cap\Log^{-1}_t(U(\Delta'))$.
Furthermore, we have the following proposition.
\begin{prop}\label{pp}
For any $\epsilon'>0$ there exists $T>1$ such that
for every $t\ge T$ and $\Delta'\in\subdiv_C$
the image
$\Phi_{\Delta',t}(V^{\zeta}_t\cap\Log^{-1}_t(U(\Delta')))$
is contained in the $\epsilon'$-neighborhood (with respect to
the product metric in $\tor\approx\R^2\times (S^1)^2$) of the
complex tropical curve $V^\zeta_{\Delta',\infty}\subset\tor$.
\end{prop}
\begin{proof}
Suppose that $z\in\tor$ is outside of the $\epsilon$-neighborhood of
$V^\zeta_{\Delta',\infty}$. Then for sufficiently large $t$
the absolute value of the sum
$|\sum\limits_{j\in\Delta'} \arg(\zeta_j)t^{\eta_j}z^j|$
is larger than zero by a definite amount that can be made large
than $|\sum\limits_{j\in\Delta\setminus\Delta'} \arg(\zeta_j)t^{\eta_j}z^j|$ so
such a point cannot be contained in
$\Phi_{\Delta',t}(V^{\zeta}_t\cap\Log^{-1}(U(\Delta')))$.
\end{proof}
This proposition can be considered as a tropical manifestation of the
patchworking principle mentioned above as it states that the geometry of
$V^{\zeta}_t\cap\Log^{-1}_t(U(\Delta'))$ is close to the geometry
of $\Phi_{\Delta',t}^{-1}(V^\zeta_{\Delta',\infty})$ no matter what are the values
of $\zeta_j$ for $j\notin\Delta'$ (as long as $\zeta\in\DD$).

\begin{coro}\label{coro-pp}
Suppose that $\Delta'\in\subdiv_C$ is a 2-dimensional polygon and
$V^\zeta_{\Delta',\infty}\subset\tor$ is a complex tropical curve
of genus $g'$.
If $t$ is sufficiently large then $V^\zeta_t\cap U(\Delta')$ is of genus
not less than $g'$.

Suppose that $\Delta'\in\subdiv_C$ is a 1-dimensional polygon and
$V^\zeta_{\Delta',\infty}\subset\tor$ has $k'$
connected components.
If $t$ is sufficiently large then $V^\zeta_t\cap U(\Delta')$
has at least $k'$
connected components.
\end{coro}
\begin{proof}
The first statement follows from Proposition \ref{pp}
since by Definition \ref{genusinfty} a complex tropical curve
cannot be approximated by curves of smaller genus.
The second statement follows from uppersemicontinuity of
the number of connected components.
\end{proof}

\IGNORE{
For any $\epsilon'>0$ there exists $T>1$ such that for any $t\ge T$
we can find neighborhoods $U'_t(\Delta')$
with the following properties.
\begin{itemize}
\item
If $\Delta'$ is a 2-dimensional polygon then
$U'_t(\Delta')$ is a neighborhood of the vertex of $C$
dual to $\Delta'$. Furthermore, $U'_t(\Delta')$ is disjoint
from all other vertices of $C$
\item
If $\Delta'$ is a 1-dimensional polygon then
the intersection $U'_t(\Delta')\cap C$
is contained in the edge of $C$ dual to $\Delta'$.
\item
We have
$$\bigcup\limits_{\Delta'
}
U'_t(\Delta')\supset \Log_t(V_t^\zeta)$$
for any $\zeta\in\DD$.
\item
$\phi_{\Delta',t_1}(U'_{t_1}(\Delta'))\subset
\phi_{\Delta',t_2}(U'_{t_2}(\Delta'))$ if $t_1<t_2.$
\item
$$\bigcup\limits_t\phi_{\Delta',t}(U'_{t}(\Delta'))=\R^2.$$
\item
There exists $\epsilon'(t)>0$ such that
$\lim\limits_{t\to+\infty}\epsilon'(t)=0$ and
$$V_t^\zeta\cap\Log^{-1}(\phi_{\Delta',t}(U'_{t}(\Delta')))$$ is
contained in the $\epsilon'(t)$-neighborhood of $V'^\zeta_{\Delta'}\subset\tor$
(with respect to the product metric in $\tor\approx\R^2\times S^1\times S^1$).
\end{itemize}
\end{prop}
}


The problem of passing from an approximate solution to
the exact solution for large $t$ is ready to be localized.
We need to find the number of such choices for $\zeta$
that $V_t^\zeta\supset\qq$ and $V_t^\zeta$ has genus $g$.
Recall that the amoeba $\Log_t(V_t)\subset\R^2$ is contained in
a small neighborhood $\NN_\epsilon(C)\supset C$
and thus, by Lemma \ref{keylem-a}, we have
$$\Log_t(V^\zeta_t)\subset\bigcup\limits_{\Delta'\ :\ \dim(\Delta')>0}U(\Delta')$$
for large $t$.
The curve $C\subset\R^2$ is simple
and thus (if considered as a subspace of $\R^2$) is a graph
with 3- and 4-valent vertices.
Every edge or vertex of $C$ corresponds to a subpolygon
$\Delta'\in\subdiv_C$. This subpolygon is a triangle,
parallelogram or an edge for a 3-valent vertex, a 4-valent vertex
or an edge (respectively).

\IGNORE{
Our goal is to cover $\NN_\epsilon(C)$ with open sets $U'_t(\Delta')$
corresponding to these subpolygons and keep track of
$V_t^\zeta\cap\Log^{-1}(U'_t(\Delta'))$ separately.

We require that our choice of the neighborhoods $U'_t(\Delta')$
is made so that the following holds.
\begin{itemize}
\item
If $\Delta'\subset\subdiv_V$ is a 2-dimensional subpolygon
corresponding to a vertex $v'\in C$ (we write
$\Delta'\subset\operatorname{Faces}(\subdiv_C)$ in this case)
then $U'_t(\Delta')\subset\R^2$ is
a small open neighborhood of $v'$ such that $U'(\Delta')\cap C$
is connected.
\item
In the case if $\Delta'\subset\subdiv_V$ is an edge
dual to an edge $E'\in C$
then we require that $U'_t(\Delta')\subset\R^2$ is
a small open neighborhood of
$$E'\setminus
\bigcup\limits_{\Delta'\subset\operatorname{Faces}(\subdiv_C)} U(\Delta')$$
such that $U(\Delta')\cap C$ is an interval.
\end{itemize}
}
\IGNORE{
Suppose that the neighborhoods $U(\Delta')$
are chosen so that
\begin{equation}\label{okrest}
\NN_\epsilon(C)\subset
\bigcup\limits_{\Delta'\in\subdiv_C\setminus\operatorname{Vert}(\subdiv_C)}
U(\Delta')
\end{equation}
(recall that $t$ is chosen so that $V^\zeta_t\subset\NN_\epsilon(C)$
for any $\zeta\in\DD$).
}

\begin{prop}\label{localize}
Suppose that $\zeta\in\DD$,
$V^\zeta_t\supset\qq$
and $t$ is large.
The curve $V^\zeta_t$ is a curve of genus $g$
if and only if all the following conditions hold.
\begin{itemize}
\item
If $\Delta'$ is a parallelogram with vertices $k_0,k_1,k_2,k_3\in\Z^2$,
$k_3-k_2=k_1-k_0$,
then $V^\zeta_t\cap\Log^{-1}_t(U(\Delta'))$
is a union of two (not necessarily connected) curves,
one in a small neighborhood of a complex tropical curve
with the Newton polygon $[k_0,k_1]$
and one in a small neighborhood of a complex tropical curve
with the Newton polygon $[k_2,k_3]$.
\item
If $\Delta'$ is an edge then $V^\zeta_t\cap\Log^{-1}_t(U(\Delta'))$
is homeomorphic to an immersed annulus (and, therefore, connected).
\item
If $\Delta'$ is a triangle then $V^\zeta_t\cap\Log^{-1}_t(U(\Delta'))$
has genus 0.
\end{itemize}
\IGNORE{
Conversely, suppose that
and $V^\zeta_t$, $\zeta\in\DD$, satisfies to the conditions above.
Then the genus of $V^t_\zeta$ is $g$.
In addition, if $t$ is large then
we have $\am_t^\zeta=\Log(V_t^\zeta)\subset\NN_\epsilon(C)$.
}
\end{prop}
\begin{proof}
Recall that each point of $\dd\Delta\cap\Z^2$ is a vertex of $\subdiv_C$.
If $V_t^\zeta$ satisfies the conditions of Proposition \ref{localize}
then no component of $V^\zeta_t\cap\Log^{-1}_t(U(\Delta'))$ has positive
genus. The last condition (when $\Delta'$ is an edge)
guarantees that the genus of $V_t^\zeta$ coincides with the genus of $C$.
Any other choice gives a higher genus.
\end{proof}

The following corollary takes care the last statement of Lemma \ref{keylem-b}.

\begin{coro}
If a curve $C=C_j$ from \eqref{gammy} is irreducible then any
$J_t$-holomorphic curve $V$ of genus $g$ and degree $\Delta$
with large $t$
passing through $\qq$ with $\Log(V)\subset\NN_{\epsilon}(C_j)$
is irreducible while if $C_j$ is reducible than any such curve $V$
is reducible.
\end{coro}
\begin{proof}
For large $t$ the $J_t$-holomorphic curve $V$ must appear as a curve $V_t^\zeta$
for some $\zeta\in\DD$ since $\Log_t(V)\subset\NN_{\epsilon}(C_j)$.
Recall that $C\subset\R^2$ is a simple tropical curve and therefore
admits a unique simple parameterization by a 3-valent graph which
is connected if and only if $C$ is irreducible.
By Proposition \ref{localize} each component of the parameterizing
3-valent graph corresponds to a component of $V_t^\zeta$.
\end{proof}



\IGNORE{
\begin{lem}\label{delta-vertex}

\end{lem}

\begin{lem}\label{delta-par}
Let $\Delta'\subset\R^2$ be a parallelogram with vertices
$k_0,k_1,k_2,k_3\in\Z^2$, $k_1-k_0=k_3-k_2$.
For any choice $b_{k_0},b_{k_1},b_{k_2}\in\C^*$
there exist
$l_{1}l_{2}=(\#([k_0,k_1]\cap\Z^2)-1)(\#([k_0,k_2]\cap\Z^2)-1)$
choices of coefficients
$\{b_j\}$, $j\in (\Delta'\cap\Z^2)\setminus \{k_0,k_1,k_2\}$,
such that the curve
$$\{z\in\tor\ |\ \sum\limits_{j\in\Delta'\cap\Z^2} b_jz^j=0\}$$
is a reducible curve of genus $-1$ and degree $\Delta'$
that splits to two irreducible components of degrees $[k_0,k_1]$
and $[k_0,k_2]$ respectively.

Suppose that the asymptotic directions of the curve corresponding to the sides
$[k_0,k_1]$ and $[k_0,k_2]$ are chosen.
There is a unique choice of coefficients $\{b_j\}$ above
that agrees with the choice of these directions.
\end{lem}
\begin{proof}
Any such curve is a union of two holomorphic cylinders taken
with multiplicities $l_1$ and $l_2$ respectively.
There are $l_1$ choices for the first cylinder and $l_2$
choices for the second.
\end{proof}
}



Recall that $C\supset\ppp\in\R^2$. Let $\XX\subset\Delta$
be the tree given by Proposition \ref{tree}.
Without loss of generality we may assume that $j_0\in\Xi$.
The number $\zeta_{j_0}=1$ is already determined.
Let us orient $\XX$ so that $j_0$ is its only source.

We partite the points of $(\Delta\subset\Z^2)\setminus\{j_0\}$
into sets $G_{\Delta'}\subset\Delta\cap\Z^2$,
$\Delta'\in\subdiv_C$ using the following rules.
\begin{itemize}
\item Suppose that $\Delta'\subset\subdiv_C$ is a parallelogram then one
of its diagonal, say $[j',j]$, is contained in the tree $\XX$.
Suppose that $[j',j]$ is oriented positively with respect
to the chosen orientation of $\XX$.
Let $k_1$ and $k_2$ be the other two vertices of the parallelograms $\Delta'$.
We let $$G_{\Delta'}=(\Delta'\cap\Z^2)\setminus ([j',k_1]\cup [j',k_2]).$$
\item Suppose that $\{j\}$, $j\neq j_0$ is a vertex of $\subdiv_C$
that is not contained in $G_{\Delta'}$ for any parallelogram $\Delta'$.
We let $$G_{\{j\}}=\{j\}.$$
If $j$ is contained in
$G_{\Delta'}$ for some parallelogram $\Delta'$ then we let
$G_{\{j\}}=\emptyset$. We also let $G_{\{j_0\}}=\emptyset$.
\item Suppose that $[j_1,j_2]\in\subdiv_C$,
is an edge such that $[j_1,j_2]\setminus\{j_1,j_2\}$ is disjoint from
$G_{\Delta'}$ for every parallelograms $\Delta'\subset\subdiv_C$.
We let $$G_{[j_1,j_2]}=([j_1,j_2]\cap\Z^2)\setminus\{j_1,j_2\}.$$
If $[j_1,j_2]\setminus\{j_1,j_2\}$ is not disjoint from
$G_{\Delta'}$ for some parallelogram $\Delta'$ then we let
$G_{[j_1,j_2]}=\emptyset$.
\item Suppose that $\Delta'\in\subdiv_C$ is a triangle with vertices
$k_0,k_1,k_2$.
We let $$G_{\Delta'}=(\Int(\Delta')\cap\Z^2).$$
\end{itemize}
Clearly, the sets $G_{\Delta'}$ are disjoint.
Furthermore, since $\XX$ contains all vertices of $\subdiv_C$
we have
$$\bigcup\limits_{\Delta'\in\subdiv_C} G_{\Delta'}=
(\Delta\subset\Z^2)\setminus\{j_0\}.$$

For a given large $t$
we say that $\zeta\in\DD$ is {\em $\Delta'$-compatible} in the following cases.
\begin{itemize}
\item Suppose that $\Delta'\in\subdiv_C$ is a positive-dimensional
polygon with $G_{\Delta'}\neq\emptyset$.
We say that $\zeta\in\DD$ is $\Delta'$-compatible
if the condition of Proposition \ref{localize} corresponding to $\Delta'$
holds for $V^\zeta_t$.
\item Suppose that $\Delta'=\{j\}$ is a vertex of $\subdiv_C$ with
$G_{\Delta'}\neq\emptyset$. Then $j$ is the endpoint of
a unique oriented edge from $\Xi$ and thus corresponds to a point
$p\in\ppp$. We say that $\zeta\in\DD$ is $\Delta'$-compatible
if $V^\zeta_t\ni q$, where $q\in\qq$ and $\Log(q)=p$.
\IGNORE{
\item If $\Delta'=[j',j]$ is an edge of the forest $\Xi$ we require
that $V_t^\zeta\cap\Log^{-1}(U'_t(\Delta'))$ is an immersed cylinder
that contains a point from $\qq$.
\item If $\Delta'=[j_1,j_2]$ is an edge
such that $[j_1,j_2]\setminus\{j_1,j_2\}$ is disjoint from
$G_{\Delta'}$ for all parallelograms $\Delta'\subset\subdiv_C$
we require that $V_t^\zeta\cap\Log^{-1}(U'_t(\Delta'))$ is an immersed cylinder.
\item If $\Delta'$ is a triangle we require that
$V_t^\zeta\cap\Log^{-1}(U'_t(\Delta'))$ is a rational curve with 3 ends.
}
\item If $G_{\Delta'}=\emptyset$ then any $\zeta\in\DD$ is
by default $\Delta'$-compatible.
\end{itemize}
By Proposition \ref{localize} the curve $V^\zeta_t\subset\tor$ with $\zeta\in\DD$
contains $\qq$ and has genus $g$ if and only if $\zeta$
is $\Delta'$-compatible for every $\Delta'\in\subdiv_C$.

Lemmas \ref{delta-vertex}, \ref{delta-par} and \ref{delta-triangle}
compute the number
of $\Delta'$-compatible choices for individual polygons $\Delta'$.

\begin{lem}\label{delta-vertex}
Let $[k',k'']$ be the edge of $\Xi$ and $q\in\tor$ be any point.
For any choice of $b_j\in\C$, $j\in[k',k'']\setminus\{k''\}$
there exists a unique choice of $b_{k''}$ such that $q$ is
a point of
$$\{z\in\tor\ |\ \sum\limits_{j\in[k',k'']}b_jz^j=0\}.$$
\end{lem}
\begin{proof}
The equation $\sum\limits_{j\in[k',k'']}b_jz^j$ is linear on $b_{k''}$
and thus has a unique solution.
\end{proof}

\begin{lem}\label{delta-par}
Let $\Delta'\subset\R^2$ be a parallelogram with vertices
$k_0,k_1,k_2,k_3\in\Z^2$, $k_1-k_0=k_3-k_2$.
For any choice $b_{j}\in\C^*$, $j\in [k_0,k_1]\cup [k_0,k_2]$,
there exists a unique
choice of coefficients
$\{b_j\}$, $j\in (\Delta'\cap\Z^2)\setminus([k_0,k_1]\cup [k_0,k_2]),$
such that the curve
$$\{z\in\tor\ |\ \sum\limits_{j\in\Delta'\cap\Z^2} b_jz^j=0\}$$
is a union of two curves with Newton polygons $[k_0,k_1]$
and $[k_0,k_2]$ respectively.
%
\end{lem}
\begin{proof}
The coefficients $b_j$ should be such that the corresponding
curve is the union of two curves
$\{z\in\tor\ |\ \sum\limits_{j\in[k_0,k_1]}b_jz^j=0\}$
and $\{z\in\tor\ |\ \sum\limits_{j\in[k_0,k_2]}b_jz^j=0\}$.
\end{proof}


Our next task is to deal with curves corresponding
to the triangles in $\subdiv_C$.
\begin{rmk}\label{proends}
If $V\subset\tor$ is an algebraic curve of degree $\Delta$
and if $\Delta\subset\R^2$ is a lattice polygon with $j$ sides then $V$ has at
least $j$ ends and at most $\#(\dd\Delta\cap\Z^2)$ ends.
Indeed, any such end corresponds to an intersection of the closure
of $V$ in the toric surface $\C T_\Delta$ with one of the boundary
divisor of $\C T_\Delta$ (which in turn corresponds to a side of $\Delta$).
\end{rmk}
\begin{lem}\label{3up}
Let $V\subset\tor$ be an algebraic curve in $\tor$ homeomorphic
to a sphere punctured 3 times. Then there exists a (multiplicative)
group endomorphism
$M:\tor\to\tor$ and a complex line
$$\C\Lambda=\{(z_1,z_2)\in\tor\ |\ b_1z_1+b_2z_2+b_0=0\},$$
$b_0,b_1,b_2\in\C^*$, such that $V=M(\C\Lambda)$.
\end{lem}
\begin{proof}
By Remark \ref{proends} the Newton polygon $\Delta'\subset\R^2$
of $V$ is a triangle.
Let $M=M_{\Delta'}$, where $M_{\Delta'}$ is defined by \eqref{MDelta}.
Its compactification is $\bar{M}_{\Delta'}:\cp^2\to\C T_{\Delta'}$
defined by \eqref{barMDelta}. Any side $\Delta''\subset\dd\Delta'$
corresponds to a boundary divisor $\C T_{\Delta''}$ and near
a general point of $\C T_{\Delta''}$ (away from intersection points
of boundary divisors) the map $\bar{M}_{\Delta'}$
is a branched covering with branching locus over $\C T_{\Delta''}$
and branching index $\#(\Delta''\cap\Z^2)-1$.

We claim that $V$ lifts under the covering $M_\Delta$.
Topologically our curve $V$ is a sphere with punctures and thus
its fundamental group is generated by loops going around the punctures.
Each loop goes around a boundary divisor $\C T_{\Delta''}$
$\#(\Delta''\cap\Z^2)-1$ times as this is the intersection number
of $\C T_{\Delta''}$ and the closure of $V$ in $\C T_{\Delta'}$.
Thus the closure of $V$ lifts to a closed surface in $\cp^2$
that is holomorphic and intersects each boundary divisor
of $\cp^2$ (that is $\cp^1\subset\cp^2$) once. Thus the lift
is a line disjoint from the intersection of the coordinate axes.
\end{proof}
\begin{coro}
Suppose that $V\subset\tor$ is a rational curve
of degree $\Delta'$,
where $\Delta'\subset\R^2$ is a lattice triangle
with no lattice points on its boundary except its vertices
(i.e. such that $\#(\Delta'\cap\Z^2)=3$).
Then $V=M_{\Delta'}(\C\Lambda)$ for a complex line
$$\C\Lambda=\{(z_1,z_2)\in\tor\ |\ b_1z_1+b_2z_2+b_0=0\},$$
$b_0,b_1,b_2\in\C^*$.
\end{coro}

By the {\em asymptotic direction of $V\subset\tor$
corresponding to a side $\Delta'\subset\dd\Delta$} we mean
the intersection point of the closure $\bar{V}\subset\C T_{\Delta}$
with the divisor $\C T_{\Delta'}$ assuming there is only
one such point (perhaps not transverse).

\begin{lem}\label{delta-triangle}
Let $\Delta'\subset\R^2$ be a triangle with vertices
$k_0,k_1,k_2\in\Z^2$. For any choice $b_{k_0},b_{k_1},b_{k_2}\in\C^*$
there exist $2\area(\Delta')$ distinct choices of coefficients
$\{b_j\}$, $j\in (\Delta\cap\Z^2)\setminus\operatorname{Vert}(\Delta')$,
such that the curve
\begin{equation}\label{Vb}
V^b=\{z\in\tor\ |\ \sum\limits_{j\in\Delta'\cap\Z^2} b_jz^j=0\}
\end{equation}
is a rational (i.e. genus 0) curve of degree $\Delta'$ with
3 ends at infinity.

Furthermore,
for any choice of the asymptotic directions of $V^b$ corresponding
to the sides $[k_0,k_1]$ and $[k_0,k_2]$
we have $\frac{2\area(\Delta')}{l_1 l_2}$ choices of
coefficients $\{b_j\}$, $j\in (\Delta'\cap\Z^2)$,
such that the curve
$V^b$ defined by \eqref{Vb}
is a rational curve of degree $\Delta'$ with
3 ends at infinity and with the given choice of the asymptotic
directions corresponding to $[k_0,k_1]$ and $[k_0,k_2]$.
Here $l_{1}=\#([k_0,k_1]\cap\Z^2)-1$ and $l_{2}=\#([k_0,k_2]\cap\Z^2)-1$.


\ignore{
Suppose that the asymptotic directions of the curve $V^b$ corresponding
to all three sides of $\Delta'$ are chosen then we have one or none
such coefficient choices. Out of the total of $l_0l_1l_2$ choices
of the asymptotic directions, $2\area(\Delta')$ have one such choice
of coefficients, $l_0=\#([k_1,k_2]\cap\Z^2)-1$.
}
\end{lem}
\begin{proof}
Consider the (singular) covering
$\bar{M}_{\Delta'}:\cp^2\to\C T_{\Delta'}$
of degree $2\Area(\Delta')$ defined in \eqref{barMDelta}.
By Lemma \ref{3up} any rational curve $V$ of degree $\Delta'$
with 3 ends is an image of a line in $\tor$.

Consider the closure $\bar{V}$ of $V$ in
$\C T_{\Delta'}$.
Since $b_{k_0},b_{k_1},b_{k_2}$ are fixed we have
$l_{1}$ possibilities for
the (unique) intersection point $p_1=\bar{V}\cap\C T_{[k_0,k_1]}$
and $l_{2}$ possibilities for
the point $p_2=\bar{V}\cap\C T_{[k_0,k_2]}$.
The points $p_1$ and $p_2$ have $\frac{2\area(\Delta')}{l_{1}}$ and
$\frac{2\area(\Delta')}{l_{2}}$ inverse images under
the map $\bar{M}_{\Delta'}$.
Connecting different liftings for different
choices of $p_1$ and $p_2$ we get $(2\area(\Delta'))^2$ different
lines in $\cp^2$ that project to $2\area(\Delta')$ different
rational curves in $\C T_{\Delta'}$.
\end{proof}


\begin{lem}\label{lem-order}
There exists an order on the polygons $\Delta'\in\subdiv_C$
such that if $\Delta'$ is greater than
$\Delta''$ then $G_{\Delta'}$ is disjoint from $\Delta''\cap\Z^2$.
\end{lem}
\begin{proof}
Clearly we can ignore the polygons $\Delta'\in\subdiv_C$ with
$G_{\Delta'}=\emptyset$ by assigning to them the highest possible
weight. We can do the same for the triangles as if $\Delta'\in\subdiv_C$
is a triangle then $G_{\Delta'}$ is disjoint from any other polygon
in $\subdiv_C$.

To sort out the remaining $\Delta'\in\subdiv_C$
we choose an order for the vertices and edges of the tree $\XX$
so that it agrees with the already chosen orientation of $\XX$
(recall that this is the orientation such that the only sink is $j_0$).
This means that $\Delta'$ must have a higher order than $\Delta''$ if
to connect $\Delta'$ to $j_0$ we have to pass through $\Delta''$.

Recall that each edge $E$ of $\XX$ is either an edge $\Delta'$ of $\Xi$
or a diagonal of a parallelogram $\Delta'\in\subdiv_C$.
Thus the order on vertices and edges of $\XX$ yields
the required order on $\subdiv_C$.
%
\end{proof}

Let $\Delta'_1,\dots,\Delta'_N$ be the polygons
from $\subdiv_C$ enumerated according to an order given
by Lemma \ref{lem-order}.

\IGNORE{
\begin{defn}
The {\em inner multiplicity} $\mu'_k$ of a polygon $\Delta'_k$
is the following number.
\begin{itemize}
\item If $\Delta'_k$ is an edge with $G_{\Delta'_k}\neq\emptyset$ then
$\mu'_k=(l_k)^2$, where the integer length $l_k$
is defined as $\#(\Delta'_k\cap\Z^2)-1$.
\item If $\Delta'_k$ is a triangle then
$\mu'_k=\frac{2\area(\Delta'_k)}{l'_k}$,
where $l'_k$ is the product of the integer lengths
of the 3 sides of $\Delta'_k$. Note that in this case $\mu'_k$
is often non-integer.
\item If $\Delta'_k$ is a parallelogram or a vertex or
any polygon with $G_{\Delta'_k}=\emptyset$ then $\mu'_k=1$.
\end{itemize}
\end{defn}

For any edge $E$ of $\subdiv_C$ contained entirely in $\dd\Delta$
we have $G_E=\emptyset$ since $\qq$ is in tropically general
position (and thus $x=s$ by Proposition \ref{tropfinite}).
Therefore $\prod\limits_{u=1}^N \mu'_u=\mult(C)$.
Lemma \ref{keylem-b} follows inductively from the following
proposition.
}

Recall that $V_{\infty}\subset\tor$ is one of
the $\mult(C)/\mu_{\operatorname{edge}}(C,\ppp)$ complex tropical
curves of genus $g$ passing via the configuration $\qq$
(by Proposition \ref{ctropmain}).
Let $\NN(V_\infty)$ be the $\epsilon$-neighborhood of $V_\infty$ in $\tor$
(recall that $\epsilon>0$ is chosen to be small).
Let $\mu'(\Delta'_k)=\#(\Delta'_k\cap\Z^2)-1$ if $\Delta'_k$
is an edge disjoint from $\ppp=\Log(\qq)$,
$\mu'(\Delta'_k)=(\#(\Delta'_k\cap\Z^2)-1)^2$ if $\Delta'_k$
is an edge not disjoint from $\ppp=\Log(\qq)$,
and $\mu'(\Delta'_k)=1$ otherwise.
Lemma \ref{keylem-b} and Theorem \ref{main} inductively follow
from the next proposition.

\begin{prop}\label{ind-step}
%
Suppose that $t$ is large and $\zeta\in\DD$ is chosen compatible with
$\Delta'_1,\dots,\Delta'_{k-1}$.
There exist $\prod\limits_{u=1}^k\mu'(\Delta'_u)$ choices
of $\zeta'\in\DD$
with the following properties.
\begin{itemize}
\item
The parameter $\zeta'$ is compatible with
$\Delta'_1,\dots,\Delta'_{k}$.
\item
We have $\zeta'_j=\zeta_j$ if $j\in G_{\Delta'_u}$, $u>k$.
\item
$V^{\zeta'}_t\subset\NN(V_\infty)$.
\end{itemize}

\IGNORE{
If $\Delta'_k$ is a triangle then we have
either one or none such
choices of $\zeta'\in\DD$ depending on the choice of asymptotic directions
corresponding to the sides $\Delta'_u$ of $\Delta'_k$,
(note that $u<k$).
We have a total of ${l'_k}$
of possibilities for this choice
of asymptotic direction.
For $2\area(\Delta')$ of these $l'_k$ choices
we have a unique choice of coefficients.
}
\end{prop}
\begin{proof}
We have the coefficient $\zeta'_j=\zeta_j$ already chosen
for $j\in G_{\Delta'_u}$, $u>k$.
Suppose that $j\in G_{\Delta'_k}$. Let us vary the corresponding
coefficients $\zeta'_j$ within $\DD$, i.e. within the disk
$|\zeta'_j-a_j|<\epsilon'_j$ for
$j\in\operatorname{Vert}(\subdiv_C)\cap G_{\Delta'_k}$ while
$\log|\zeta'_j|-\beta_j<\epsilon_j$
for $j\notin\operatorname{Vert}(\subdiv_C)\cap G_{\Delta'_k}$.
Denote the corresponding $\#(G_{\Delta'_k})$-dimensional polydisc
with $\DD_k$.

\IGNORE{
Inductively, Proposition \ref{ind-step} gives
$\prod\limits_{u=1}^{k-1} \mu'_u$ choices of coefficients
$\zeta'_j$, ${j\in G_{\Delta'_u}}$, $u<k$.
By Proposition \ref{patch-principle} and Lemma \ref{lem-order}
we can choose $T>1$ so large that
for any choice of
$$\eta=\{\zeta'_j\}_{j\in G_{\Delta'_k}}\in\DD_k$$
the curve
$V^{\zeta'}_t\cap\Log^{-1}(U'_t(\Delta'_u))$
is contained in a small neighborhood of
$V^\zeta_t\cap\Log^{-1}(U'_t(\Delta'_u))$.

This yields a map $$\Phi_k:\DD_k\to \bm_{l,s},$$
where $l=\#(\Int\Delta\cap\Z^2)$ and
$\bm_{l,s}\supset\M_{l,s}$ is the Deligne-Mumford
compactification (see e.g. \cite{HM})
$$\Phi_k:\eta\mapsto V_t^{\zeta'}.$$
The map $\Phi_k$ is {\em a priori} only partially defined
(for $\eta$ such that $V_t^{\zeta'}$ has singularities
other than nodes).
}
We have a map $$\Phi_k:\DD_k\to \M(\Delta),$$
where $\M(\Delta)$ is the space of all curves in $\tor$
given by polynomials whose Newton polygon is contained in $\Delta$.
Consider also the map $\Psi_k:\DD_k\to \M(\Delta),$
$\Psi_k(\zeta)=V_{\Delta',t}^\zeta,$
where $V_{\Delta',t}^\zeta$ is the zero locus of
$f_{\Delta',t}^\zeta=\sum\limits_{j\in\Delta'}\arg(\zeta_j)t^{|\zeta_j|}$.

The polydisc $\DD$ is a product of disks in $\C$.
Recall that the {rim} $\delta\DD_k$ of $\DD_k$ is the product
of the boundaries of the disks from this product.
By Corollary \ref{coro-pp} if $t$ is large then
both $\Phi_k(\delta\DD_k)$ and $\Psi_k(\delta\DD_k)$ consists of nonsingular curves.
Furthermore, by Proposition \ref{pp} and Corollary \ref{coro-pp}
the curves from $\Phi_k(\delta\DD_k)$ and $\Psi_k(\delta\DD_k)$
are close to the corresponding complex tropical curves
that are ``non-singular" in $\Log_t^{-1}(U(\Delta'_k))$,
i.e. have the highest possible genus or the number of connected
components for the given Newton polygon $\Delta'$, so that
$\Phi_k$ and $\Psi_k$ can be connected by a homotopy whose
restriction to $\delta\DD_k$ stays within the curves that have non-singular
intersection with $\Log^{-1}_t(U(\Delta'_k))$.
Therefore, the intersection number of $\Phi_k(\DD_k)$ and $\Psi_k(\DD_k)$
with the locus of $\Delta'_k$-compatible curves is the same.
Similarly, $\Phi_u$ and $\Psi_u$ can also be connected by a homotopy whose
restriction to $\delta\DD_u$ stays within the curves that have non-singular
intersection with $\Log^{-1}_t(U(\Delta'_u))$ because of Proposition \ref{pp}
and since the order of $\Delta'_u$ is taken as in Lemma \ref{lem-order}.

Suppose that $\subdiv_C$ foes not have multiple edges,
i.e. edges $E$ with $G_E\neq\emptyset$.
The map $\Psi_k:\DD_k\to\M(\Delta)$
intersects the stratum of $\Delta'_k$-compatible
curves in a single point by Lemmas \ref{delta-vertex}, \ref{delta-par} and
\ref{delta-triangle}.
Thus we have the same for the map $\Phi_k$ and this gives us
the values $\zeta'_j$ for $j\in\Delta'_k$ such that
the collection $\zeta'_j$, $j\in\Delta'_k$, and $\zeta_j$,
$j\notin\Delta'_k$, is $\Delta'_k$-compatible. Note that even
though this collection does no longer have to be $\Delta'_u$-compatible
for $u<k$, it is close to a compatible configuration by Proposition \ref{pp}
since $G_{\Delta'_k}\cap\Delta'_u=\emptyset$.
Thus, modifying this construction inductively at earlier steps with the
new choice of $\zeta'_j$, $j\in\Delta'_k$ we can find
the values $\zeta'_j$, $j\in\Delta_u$, $u\le k$, such that together with $\zeta_j$,
$j\in\Delta_u$, $u>k$, they are $\Delta'_u$-compatible for any $u\le k$.

This finishes the proof of the proposition (and thus the proof
of Lemma \ref{keylem-b} and Theorem \ref{main}) in the case when
$\subdiv_C$ does not have multiple edges. If there are such edges
in $\subdiv_C$ then the proposition follows from Lemma \ref{delta-edge} below.
\end{proof}

\begin{lem}\label{delta-edge}
Let $\Delta'=[k',k'']\in\subdiv_C$ be an edge, $k',k''\in\Z^2$,
such that $\Delta'\not\subset\dd\Delta$.

If $\Delta'\not\subset\Xi$ and $t$ is sufficiently large then
there exist $l'$ different choices
of coefficients
$\zeta'\in\DD$
such that $\zeta'_j=\zeta_j$ if $j=k',k''$ or
$j\in \Delta\setminus\Delta'$, $V^{\zeta'}_t\subset\NN(V_\infty)$,
and the intersection
$$V^{\zeta'}_t\cap\Log_t^{-1}(U(\Delta'))$$
is an immersed cylinder. Here $l'=\#(\Delta'\cap\Z^2)-1$ is the integer length
of $\Delta'$.

If $\Delta'\subset\Xi$ then there exist $(l')^2$ different choices
of coefficients $\zeta'\in\DD$
such that $\zeta'_j=\zeta_j$ if $j=k'$ or
$j\in \Delta\setminus\Delta'$, $V^{\zeta'}_t\subset\NN(V_\infty)$ and
the intersection
$$V^{\zeta'}_t\cap\Log_t^{-1}(U(\Delta'))$$
is an immersed cylinder which contains a point from $\qq$.
(Note that in the second case we also vary the coefficient
corresponding to one of the endpoints of the interval $\Delta'$.)
\end{lem}
\ignore{
\begin{lem}\label{delta-inedge}
Let $\Delta'=[k',k'']\in\Xi$ be an edge.
If $\Delta$t$ is sufficiently large then
there exist $l'=\#(\Delta'\cap\Z^2)-1$ different choices
of coefficients
$\zeta'\in\DD$
such that $\zeta'_j=\zeta_j$ if $j=k',k''$ or
$j\in \Delta\setminus\Delta'$
and the curve
$$\{z\in\Log^{-1}(U'_t)\ |\ \sum\limits_{j\in (\Delta\cap\Z^2)} b_jz^j=0\}$$
is an immersed cylinder that contains a point from $\qq$.
\end{lem}
}

\begin{proof}
First we note that it suffices to check
this lemma for a particular model of $\Delta$
and $\subdiv_C{\Delta}$
as long as $\Delta'\in\subdiv_\Delta$ is an edge not contained
in the boundary of $\Delta$.
This can be deduced from the patchworking principle, Proposition \ref{pp}.
Indeed, suppose that
$$f^{(1)}(z)=\sum\limits_{j\in\Delta_1\cap\Z^2}a^{(1)}_jz^j\ \ \text{and}\ \
f^{(2)}(z)=\sum\limits_{j\in\Delta_2\cap\Z^2}a^{(2)}_jz^j,$$
$a^{(1)}_j,a^{(2)}_j\in\C^*$,
are polynomials with Newton polygons $\Delta_1$ and $\Delta_2$
such that  $\Delta'\in\subdiv_{f^{(1)}},\subdiv_{f^{(2)}}$,
$\Delta'\not\subset\dd\Delta_1,\dd\Delta_2$ and $a^{(1)}_j=a^{(2)}_j$
if $j\in\Delta'$.
Embedding $\Delta_1$ and $\Delta_2$ into a larger polygon
and applying Proposition \ref{pp} in the same fashion
as in the proof of Proposition \ref{ind-step} we may assume that
$\Delta_1=\Delta_2$.
The polynomials $f^{(1)}$ and $f^{(2)}$ can be deformed to each other
by $f^{(s)}$, $1\le s\le 2$ with $a^{(s)}_j=a^{(1)}_j$ for $j\in\Delta'$.

Let $$\DD'=\{\zeta\in\C^{\#(\Delta'\cap\Z^2)-2}\ |\
\log|\zeta_j^{(s)}|-\log|a_j^{(s)}|<\epsilon'_j,
j\in\Delta'\setminus\{k',k''\}\}.$$
be the polydisc, where $\epsilon'_j>0$ are some small numbers.
For each $\zeta\in\C^{\#(\Delta'\cap\Z^2)-2}$ we form
$$f^{\zeta^{(s)}}_t(z)=\sum\limits_{j\in\Delta\cap\Z^2}
\arg(\zeta^{(s)}_{j})t^{\log|\zeta^{(s)}_{j}|}z^j,$$
where $\zeta_j^{(s)}=a_j^{(s)}$, if $j=k',k''$ or
if $j\in\Delta\setminus\Delta'$.
Denote the zero set of $f^{\zeta^{(s)}}_t$ in $\tor$
with $V^{\zeta^{(s)}}_t$.


As $\zeta$ runs over the rim $\delta\DD'$
the curve $V^{\zeta^{(s)}}_t$ never develops
a singularity within $\Log^{-1}_t (U(\Delta'))$
by Corollary \ref{coro-pp}. Furthermore,
the intersection
$V^{\zeta^{(s)}}_t\cap\Log^{-1}_t (U(\Delta'))$
is a union of $\#(\Delta'\cap\Z^2)-1$ disjoint cylinders
whose mutual position changes when $\zeta$ changes within $\DD'$
so any $\zeta\in\delta\DD'$ is ``maximally $\Delta'$-incompatible"
for any $1\le s\le 2$.
Note that the locus of $\Delta'$-compatible
values of $\zeta$ is locally given as an intersection of hypersurfaces
and for every $\zeta$ from those hypersurfaces
we have $V^{\zeta^{(s)}}_t\cap\Log^{-1}_t (U(\Delta'))$
consisting of at most $\#(\Delta'\cap\Z^2)-2$ components, so they
miss the rim $\delta\DD'$.
Thus the number of $\Delta'$-compatible values of $\zeta\in\DD'$ does
not depend on $s$ and we can use any model for $\Delta\supset\Delta'$
and $\subdiv_C\ni\Delta'$ as long as $\Delta'\not\subset\dd\Delta$.

First we treat the case when $\Delta'\not\subset\Xi$.
Let $\Delta'=[(0,0),(l',0)]$.
Suppose that $l'$ is odd.
In this case we take $\Delta$ to be the parallelogram with vertices
$(0,-1)$, $(0,0)$, $(l',0)$ and $(l',1)$. Let $g=0$ and $j_0=(0,0)$.
We have $s+g-1=3$. We may choose $\qq=\{q_1,q_2,q_3\}$
so that the only tropical curve $C$ passing via $\Log(\qq)$
has $\subdiv_C\ni [(0,0),(l',0)]$ or so that
$\subdiv_C\ni [(0,-1),(l',1)]$.
Both choices can be made so that the forest $\Xi$ from Proposition \ref{forest}
consists of the edges $[(0,-1),(0,0)]$,
$[(0,-1),(l',0)]$ and $[(l',0),(l',1)]$.
The points $\qq$ determine the coefficients
$\zeta_{(0,0)}, \zeta_{(0,-1)}, \zeta_{(l',0)},
\zeta_{(l',1)}$ since $\dd\Delta\cap\Z^2$ consists only of the
vertices of $\Delta$.
We need to determine the number of compatible choices
for coefficients at the points
$j\in\Int(\Delta)$.
Note that there are no lattice points inside $[(0,-1),(l',1)]$.

To compute $\ncl(0,\Delta)$ we may use both configurations.
For the first choice of $\qq$ we have $\ncl(g,\Delta)=N$,
where $N$ is the sum of the numbers
of choices of $\zeta'$ we need to find for all possible tropical curves
$V^{(j)}_\infty$ of degree $\Delta$ and genus zero passing via $\qq$.
All of them have are mapped by $\Log$ to the same tropical curve $C$,
since $C$ is the only tropical curve of degree $\Delta$ and genus zero
passing via $\Log(\qq)$. There are $l'$ of them by Proposition \ref{ctropmain}.
Thus, we have $N/l'$ choices for $\zeta'$
so that $V^{\zeta'}_t$ is contained in the neighborhood of any individual
curve by symmetry.
%


For the second choice of $\qq$ we do not have multiple edges
in $\subdiv_C$. Thus Lemma \ref{keylem-b} and Theorem \ref{main}
are already established for this choice of configuration.
We have $\ncl(0,\Delta)=(l')^2$ from the two triangles
of $\subdiv_C$ and, therefore, $N=(l')^2$ and $N/l'=l'$.

If $l'$ is even we choose $\Delta$ to be the quadrilateral
with vertices $(0,-1)$, $(0,0)$, $(l',0)$ and $(l'-1,1)$.
The rest of the proof is the same.

If $\Delta'\in\Xi$ we make same choices for our model of
the ambient polygon $\Delta$ (depending on the parity of $l'$).
However, we choose $\qq$ so that $\Xi$ consists of the edges $[(0,-1),(0,0)]$,
$[(0,0),(l',0)]$ and $[(l',0),(l',1)]$ if $l'$ is odd or
$[(0,-1),(0,0)]$, $[(0,0),(l',0)]$ and $[(l',0),(l'-1,1)]$ if $l'$ is even.
By Proposition \ref{ctropmain} there is a unique complex tropical
curve of degree $\Delta$ and genus 0 passing via $\qq$, so it must coincide
with $V_\infty$ and the number of compatible choices for $\zeta'\in\DD$
is $N(0,\Delta)=(l')^2$ (as this number was already computed above).
\end{proof}

This finishes the proof of Lemma \ref{keylem-b} and Theorem \ref{main}
in the general case.

\begin{rmk}\label{phases}
Let $E$ be the edge of $C$ dual to $\Delta'$ in Lemma \ref{delta-edge}.
Let $A,B\in C$ be the endpoints of the corresponding edge of $\Gamma$
(recall that $C$ is a simple tropical curve parameterized by
$h:\Gamma\to\R^2$). Both $A$ and $B$ are three-valent vertices
corresponding to triangles $\Delta_A,\Delta_B\subset\Delta$.
By the patchworking principle (cf. Proposition \ref{pp})
the intersection $V^\zeta_t\cap\Log^{-1}_t(U(\Delta_A))$
approximates a curve
$$V^\zeta_{\Delta_A}=\{z\in\tor\ |\
\sum\limits_{j\in\Delta_a}\arg(\zeta_j)t^{\log|\zeta_j|}z^j\}$$
which, in turn, approximates
a curve $V_A$ in $\C T_{\Delta_A}$ with the Newton polygon $\Delta_A$
that has tangency of order $l'=\#(\Delta'\cap\Z^2)-1$ with the toric
divisor corresponding to $\Delta'$.
Let $b_A\approx S^1$ be the link of this tangency,
$b_A\subset V^\zeta_t\cap\Log^{-1}_t(U(\Delta_A))$.
The map $\arg|_{b_A}:b_A\to S^1\times S^1$
approximates the $\Z_{l'}$-covering $$\beta_A:b_A\to S^1,$$
where the base $S^1$ is the geodesic circle in $S^1\times S^1$
corresponding to the phase of the holomorphic annulus of $V_{\infty}$
over the edge $E$.
Similarly, we get the $\Z_{l'}$-covering $$\beta_B:b_B\to S^1$$
for the other endpoint $B$.
There are $l'$ ways to match the arguments of $b_A$ and $b_B$ in the
corresponding curves $V_A\subset\C T_{\Delta_A}$ and $V_B\subset\C T_{\Delta_B}$.
By the $\Z_{l'}$-symmetry (in either $\C T_{\Delta_A}$ or
$\C T_{\Delta_B}$) we have equal number of choices for $\zeta'$
in Lemma \ref{delta-edge} for any of this phase matching.
Thus we have a unique choice of $\zeta'$ if $\Delta'\not\subset\Xi$ and
$l'$ choices otherwise. Similarly, the $l'$ choices in the second case
is distinguished by the $l'$ points of $\beta_A^{-1}(\alpha)$,
where $\alpha$ is the argument of the point $q_j\in\qq$ with $\Log(q_j)\in E$.
\end{rmk}

\subsection{Real curves: proof of Theorem \ref{realmain} and \ref{main-W}}
Remark \ref{phases} is useful for detecting real curves. Suppose that
$V_A$ and $V_B$ are real, i.e. invariant with respect to $\conj$.
In that case we can choose the circles $b_A$ and $b_B$
$\conj$-invariant as well. Let $p^+_A,p^-_A\in b_A$ and $p^+_B,p^-_B\in b_B$
be the points fixed by $\conj$. To get a real curve $V^{\zeta'}_t$
we have to match the real points of $b_A$ with the real points of $b_B$
with the same value of argument.

If $l'$ is odd than the intersection number of $V_A$ with the
toric divisor in $\C T_{\Delta_A}$ corresponding to $\Delta'$ is odd
as well and $p^+_A$ and $p^-_A$ belong to distinct quadrants in $\rtor$
and, similarly, $p^+_B$ and $p^-_B$ also must have distinct arguments.
Thus, only one out of $l'$ ways of matching the phase can give
a real curve $V^{\zeta'}_t$.
On the other hand, if $\qq$ is real then such $V^{\zeta'}_t$ must be real:
if not then $\conj(V^{\zeta'}_t)$ is another curve of the same genus and degree
passing through $\qq$ with the same pattern of phase matching at $\Delta'$.
If $\Delta'\in\Xi$ then we must match $q_j$ to either $p^+_A$ or $p^-_A$
(according to its quadrant) and, again, the real curve is unique.

If $l'$ is even then $p^+_A$ and $p^-_A$ are both from the same quadrant.
Also, $p^+_B$ and $p^-_B$ are both from the same quadrants.
If the quadrants of these two pairs do not coincide then no real
matching is possible and not a single curve in Lemma \ref{delta-edge} is real.
If these two quadrants coincide then we can match
$p^+_A$ to $p^+_B$ or $p^-_B$. Each of these two patterns gives
a real curve, otherwise we'd have two distinct (conjugate)
curves with the same matching pattern.
If $\Delta'\in\Xi$ then there are two ways of matching $q_j$:
to $p^+_A$ and $p^-_A$.
Again, each of these two ways has to give a real curve, so we
have a total of 4 out of $(l')^2$ curves real in this case.
This finishes the proof of Theorem \ref{realmain}.


To prove Theorem \ref{main-W} we note that if
$\subdiv_C$ contains an edge of an even integer length then
$C$ contributes zero to the Welschinger invariant.
Indeed, for each such edge $E$ we have two real branches
of $V^{\zeta'}_t\cap\Log^{-1}(U(\Delta'))$ in the same quadrant,
where $\Delta'\subset\Delta$ is the edge dual to $E$.
We have two ways of matching of the real points of $b_A$ and $b_B$
both leading to real nodal curves as in Remark \ref{phases}.
By topological reasons these two choices must have different
parity of the number of hyperbolic nodes in that quadrant.
Therefore, the total contribution of such $C$ to the Welschinger
invariant is zero.

If all edges of $\subdiv_C$ have odd integer length then
elliptic nodes can appear only from
$V^{\zeta'}_t\cap\Log^{-1}(U(\Delta'))$
where $\Delta'$ is a triangle.
This part of $V^{\zeta'}_t$ has a total of $\#(\Delta'\cap\Z^2)$
nodes. None of these nodes can be real hyperbolic since
the restriction $M_{\Delta'}|_{\rtor}:\rtor\to\rtor$ (see \eqref{MDelta})
is injective if all sides of $\Delta'$ have odd integer length.
Therefore the multiplicity $\operatorname{mult}_{V}^{\R,W}(C)$
from Definition \ref{mult-W} agrees with Welschinger's signs and
gives the right count for Theorem \ref{main-W}.

\subsection{Counting by lattice paths:
proof of Theorems \ref{thm1}, \ref{thm2} and \ref{pathsW}}
Recall that we have a linear map $\lambda:\R^2\to\R$
injective on $\Z^2$.
Let $L\subset\R^2$ be an affine line (in the classical sense)
orthogonal to $\lambda$.

We choose a configuration $\ppp=\{p_1,\dots,p_{s+g-1}\}\subset L$
so that the order of $p_k$ agrees with a linear order on $L$.
Furthermore, we choose each $p_k$ so that the distance
from $p_k$ to $p_{k-1}$ is much larger than the distance
from $p_{k-1}$ to $p_{k-2}$. Such configuration can be chosen
in a tropically general position since the slope of $L$ is irrational
(and therefore intersects any tropical curve in $\R^2$ in a finite
number of points).

Let $C\subset\R^2$ be a tropical curve of genus $g$ and degree $\Delta$
passing via $\ppp$.
Let $\Xi$ be the forest $\Xi$ from Proposition \ref{forest}.

\begin{lem}\label{noextraint}
We have $C\cap L=\ppp$.
\end{lem}
\begin{proof}
Let $K$ be a component of $\Gamma\setminus h^{-1}(\ppp)$.
Suppose that $h(K)$ intersects $L$ at a point not from $\ppp$.
One of the components of $K\setminus h^{-1}(L)$ would
yield a bounded graph with edges at $L$ contained in a half-plane.
Clearly such graph can not be balanced.
\end{proof}

\begin{coro}
The forest $\Xi\subset\Delta$ is a $\lambda$-increasing
path that connects vertices
$p$ and $q$ as in Theorem \ref{thm1}.
\end{coro}
\begin{proof}
By Lemma \ref{noextraint} the vertices of $\Xi$
correspond to the components of $L\setminus\ppp$.
Therefore,
the forest $\Xi$ is a path that connects vertices of $\Delta$.
Note that the path is $\lambda$-increasing since the linear order
on $L$ is consistent with $\lambda$.
\end{proof}

This corollary allows one to enumerate all tropical curves
of genus $g$ and degree $\Delta$ passing via $\qq$
by the corresponding paths.
Suppose that such a path $\gamma$ is chosen.

The path $\gamma$ determines the slopes of the edges that contain
points from $\ppp$ for any tropical curve with $\Xi=\gamma([0,s+g-1])$.
Let $\Line(p_j)\subset\R^2$ be the line through $p_j$ in the corresponding
direction.
We need to find all tropical curves of genus $g$ and degree $\Delta$
that pass through $\ppp$ along $\Line(p_j)$.
The sum of multiplicities of such curves should coincide
with the multiplicity of $\gamma$. We need to compute
$\mu(\gamma)$, $\mu^{\R}(\gamma,\sigma)$ or $\nu^{\R}(\gamma)$
for Theorem \ref{thm1}, \ref{thm2} or \ref{pathsW} respectively.

Let $H_+$ and $H_-$ be the two half-planes bounded by $L$.
Lemmas \ref{descompl} and \ref{noextraint} imply that
for any tropical curve $C$ of genus $g$ and degree $\Delta$
that pass through $\ppp$
$\Gamma\cap h^{-1}(H_{\pm})$ is a tree with one end at infinity,
where $h:\Gamma\to\R^2$ is a parameterization of $C$.

Recall that our definition of the multiplicities
$\mu_{\pm}(\gamma)$, $\mu_{\pm}^{\R}(\gamma,\sigma)$
and $\nu_{\pm}^{\R}(\gamma)$ was inductive.
If $k$ is such that $\Delta_\pm$
is locally strictly convex at $\gamma(k)$ then the intersection
point $p_{k,k+1}$ of $\Line(p_k)$ and $\Line(p_{k+1})$ is contained
in $H_\pm$. If $k$ is the smallest number with this property then
$p_{k,k+1}$ is the closest to $L$ intersection point of the lines
$\Line(p_j)$ in $H_\pm$.

Let $\hat{L}$ be the line parallel to $L$ and such that
the strip $A(L,\hat{L})$ between $L$ and $\hat{L}$ contains $p_{k,k+1}$
and does not contain any other intersection point of
lines $\Line(p_j)$.
We have two possible cases for $C\cap A(L,\hat{L})$.

In the first case we have a 3-valent vertex at $p_{k,k+1}$.
Then $C\cap A(L,\hat{L})$ has a new interval emanating at
$p_{k,k+1}$ and intersecting $\hat{L}$ at a point $p'_k$.
This case corresponds to the path $\gamma'$ from \eqref{pathmult}.
We set $\{p'_j\}=\Line(p_j)\cap\hat{L}$ if $j<k$,
$\{p'_j\}=\Line(p_{j+1})\cap\hat{L}$ if $j>k$ and proceed
inductively by incorporating possibilities for
$C\cap (H_{\pm}\setminus A(L,\hat{L}))$.

In the second case $p_{k,k+1}$ is a point of self-intersection
of $h:\Gamma\to\R^2$. This case corresponds to the path
$\gamma''$ from \eqref{pathmult}.
We set $\{p''_j\}=\Line(p_j)\cap\hat{L}$ if $j\neq k,k+1$,
$\{p''_k\}=\Line(p_{k+1})\cap\hat{L}$ and
$\{p''_{k+1}\}=\Line(p_{k})\cap\hat{L}$
Again, we proceed inductively with the new, smaller half-plane.
$H_{\pm}\setminus A(L,\hat{L})$.

All multiplicities $\mult(C)$, $\mu(C,\{\sigma_E\},\qqq)$ and
$\operatorname{mult}^{\R,W}(C)$ are multiplicative and therefore
we can compute them by taking the products of the corresponding
numbers in every annulus $A(L,\hat{L}))$ from the induction.
\IGNORE{
We compute these numbers inductively for each $H_{\pm}$.
Let $R_j\subset\R^2$ be the ray

 containing the edge $E_j$.
Suppose that the intersection point $\{r_{12}\}=\Line(E_1)\cap\Line(E_2)$
is contained in $H_+$. Then
...
Since the distance between $r_1$ and $r_2$ is the smallest
if these extensions intersect in $H_+$ they do so before intersecting
the extension of any other edge $E_u$.
Each such intersection can give
a 3-valent vertex of $C$. Otherwise these extension pass through
the intersection point without an interaction. The first case
is dual to a triangle in the dual subdivision. The second
case is dual to a parallelogram. Then we consider the intersection
of the result with the extension of $E_3$ and so on.
After incorporating the extension of $E_{s+g-1}$ the resulting
curve has to have the ends orthogonal to the positive part
of $\dd\Delta$, i.e. given by the path $\alpha_+$.

Note that our inductive procedure agrees with the definition
of the positive multiplicity of $\gamma$.
Similarly, the negative multiplicity of $\gamma$ agrees with
a choice of extensions to $H_-$.
Note that our procedure necessarily gives a tropical curve $C$
of degree $\Delta$ parameterized by $h:\Gamma\to C$ so that
each component of $\Gamma\setminus h^{-1}(\ppp)$ is a tree with
one end at infinity.
Therefore, the genus of $C$ is $g$.
}
Theorems \ref{thm1}, \ref{thm2} and \ref{pathsW} follow
from Theorems \ref{main}, \ref{realmain} and \ref{main-W} respectively.

\end{document}